\documentclass[twoside,12pt]{article}% {amsart}
\usepackage{amsmath,amsfonts,amscd,bezier}
 \usepackage[usenames,dvipsnames]{pstricks}
 \usepackage{epsfig}
 \usepackage{pst-grad} % For gradients
 \usepackage{pst-plot} % For axes
\usepackage{amsthm,amsmath,amsfonts,amssymb}
\usepackage{graphicx}
\usepackage{hyperref}
\usepackage{ulem}
\usepackage[latin1]{inputenc}             % Para traduzir acentos
\usepackage{color}

\newcommand{\R}{{\mathbb R}}
\newcommand{\N}{{\mathbb N}}
\newcommand{\C}{{\mathbb C}}

\renewcommand{\Re}{{\mbox{Re\,}}}

\numberwithin{equation}{section}

\newcommand{\bPto}{{\stackrel{{\cal{P}}^\beta}{\longrightarrow}}}

\newcommand{\wto}{{\stackrel{ }{-{\hspace{-2mm}}\rightharpoonup}}}

\newcommand{\dwto}{{\stackrel{E}{-{\hspace{-2mm}}\rightharpoonup}}}

\newcommand{\dto}{{\stackrel{E}{\longrightarrow}}}
\newcommand{\Pto}{{\stackrel{{\cal{P}}}{\longrightarrow}}}
\newcommand{\Eto}{{\stackrel{{{E}}}{\longrightarrow}}}
\newcommand{\ddto}{{\stackrel{EE}{\longrightarrow}}}
\newcommand{\PPto}{{\stackrel{{\cal{PP}}}{\longrightarrow}}}

\newcommand{\ccto}{{\stackrel{CC}{\longrightarrow}}}
\newcommand{\cqd}{\hfill $\rule{1.5mm}{2.5mm}$ \vspace*{0.5cm}}

\def \ts {\textstyle}
\def \ds {\displaystyle}

\def \half{\frac{1}{2}}

\newcommand{\trip}[1]{{|\kern -1pt\|#1\|\kern -1pt|}}

\def \demb {{{ \bigcap\kern -12.2pt{_\big\downarrow}}}}
\def \uemb {{{ \bigcup\kern -5.8pt{^\big\uparrow}}}}

\def \Re{{I \kern -.3em R}}
\def \Rn {{{\bf I\kern -1.6pt{\bf R}}}^{\rm n}}

\def \half {{1\over 2}}
\def \ts {\textstyle}
\def \eps {\epsilon}

\newtheorem{theorem}{Theorem}[section]
\newtheorem{lemma}[theorem]{Lemma}
\newtheorem{corollary}[theorem]{Corollary}
\newtheorem{remark}[theorem]{Remark}
\newtheorem{definition}[theorem]{Definition}
\newtheorem{proposition}[theorem]{Proposition}

\begin{document}
%\color{red}
\title{Continuity of attractors of parabolic equations with nonlinear boundary conditions and rapidly varying boundaries. The case of a Lipschitz deformation}

\author{Gleiciane S. Arag\~ao\footnote{Departamento de Ci\^encias Exatas e da Terra, Universidade Federal de S\~ao Paulo, Rua S\~ao Nicolau, 210, Diadema-SP, Cep 09913-030, Brazil. E-mail: gleiciane.aragao@unifesp.br}, \; Jos\'e M. Arrieta \footnote{Departamento de An\'alisis y Matem\'atica Aplicada, Universidad Complutense de Madrid, 28040 Madrid and Instituto de Ciencias Matem\'aticas
CSIC-UAM-UC3M-UCM, C/Nicol\'as Cabrera 13-15, Cantoblanco, 28049 Madrid, Spain. E-mail: arrieta@mat.ucm.es}\,  and \; Simone M. Bruschi \footnote{Departamento de Matem\'atica, Universidade de Bras\'ilia, Campus Universit\'ario Darcy Ribeiro, ICC Centro, Bloco A, Asa Norte, Cep 70910-900, Brazil and Department of Mathematical Sciences, George Mason University, 4400 University Dr Fairfax-VA, 22030, USA. E-mail: sbruschi@unb.br, sbruschi@gmu.edu}}

\maketitle

\begin{abstract}
In this paper we obtain the continuity of attractors for 
 nonlinear parabolic equations with nonlinear boundary conditions when the boundary of the domain varies very rapidly as a parameter $\epsilon$ goes to zero. We consider the case
 where the boundary of the domain presents a highly oscillatory behavior as the parameter $\epsilon$ goes to zero. For the case where we have a Lipschitz deformation of the boundary with the Lipschitz constant uniformly bounded in $\epsilon$ but the boundaries do not approach in a Lipschitz sense, the solutions of these equations converge in certain sense to the solution of a limit parabolic equation of the same type but where the boundary condition has a factor that captures the oscillations of the boundary. 
 To address this problem, it is necessary to consider the notion of convergence of functions defined in varying domains and the convergence of a family of operators defined in different Banach spaces. 
Moreover, since we consider problems with nonlinear boundary  conditions, it is necessary to extend these concepts to the case of spaces with negative exponents and to operators defined between these spaces. 

 \vspace{0.2cm}

\noindent \textit{keywords}: parabolic equation; nonlinear boundary conditions; varying boundary; oscillations; Lipschitz deformation; continuity of attractors. 

\vspace{0.2cm}


\end{abstract}

%\tableofcontents

%%%%%%%%%%%%%%%%%%%%%%%%%%%%%%%%%%%%%%%%%%%%%%%%%%%%%%%%%%%%%%%%%%%%%%%%%%%%%%%%%%%%%%%%%%%%%%%%%%%%%%%%%%%%%%%%%%%%%%%%%%%%%%%%%%%%%%%%%%%%%%%%%%%%%%%%%%%%%%%%%%%

\section{Introduction}
\label{introd}

\par

In this paper we are concerned with the asymptotic dynamics of  parabolic problem with nonlinear Neumann boundary conditions of the type
\begin{equation}  
\label{nbc}
\left\{
\begin{array}{lll}
\frac{\ts\partial u_{\epsilon}}{\ts\partial t}-\Delta u_{\epsilon}+u_{\epsilon}=\displaystyle  f(x,u_{\epsilon}), &\hbox{in} \ \Omega_\eps  \times (0,\infty) \\
\frac{\ts\partial u_{\epsilon}}{\ts\partial n_{\epsilon}}=g(x,u_{\epsilon}), & \hbox{on} \ \partial \Omega_\eps \times (0,\infty) \\
u_{\epsilon}(0)=u^{0}_{\epsilon}(x), &\hbox{in} \ \Omega_\eps
\end{array} \right.
\end{equation}
when the boundary of the domain varies very rapidly as a parameter $\epsilon \to 0$.
To describe the problem, we consider a family of uniformly bounded smooth domains $\Omega_\epsilon \subset \mathbb{R}^N$, with $N\geq 2$ and $0\leq\epsilon\leq \epsilon_0$, for some $\epsilon_0>0$ fixed, and we will look at this problem from the perturbation of the domain point of view. We will refer to $\Omega \equiv\Omega_0$ as the unperturbed domain and $\Omega_\eps$ as the perturbed domains. We will assume that  $\Omega_\eps\to \Omega$ and $\partial\Omega_\eps\to \partial\Omega$ in the sense of Hausdorff. Although the domains behave continuously as $\eps \to 0$, the way in which boundary $\partial \Omega_\eps$ approaches $\partial \Omega$ may not be smooth. In particular, this setting includes the case
where boundary $\partial \Omega_\eps$ presents a highly oscillatory behavior as $\eps \to 0$. In this work, we will assume that the boundary $\partial\Omega_\eps$ is expressed in local charts as a Lipschitz deformation of $\partial \Omega$ with the Lipschitz constant uniformly bounded in $\eps$.  

\par\medskip 

It is reasonable to expect that the family of solutions $\{u_{\epsilon}\}_{\epsilon \in (0,\epsilon_0]}$ of (\ref{nbc}) will converge to the solution of an equation  with a nonlinear boundary condition on $\partial \Omega$ that  inherits the information about the behavior  of the measure of the deformation of $\partial \Omega_\eps$ with respect to $\partial \Omega$. More precisely, under certain conditions, the solutions of (\ref{nbc}) converge, in some sense that we will define later,  to the solution of the parabolic problem with nonlinear Neumann boundary conditions
\begin{equation}  
\label{nbc_limite_gamma_F}
\left\{
\begin{array}{lll}
\frac{\ts\partial u_{0}}{\ts\partial t} -\Delta u_{0}+u_{0}=f(x,u_{0}), &\hbox{in} \ \Omega \times (0,\infty) \\
\frac{\ts\partial u_{0}}{\ts\partial n}=\gamma(x)g(x,u_{0}), & \hbox{on} \ \partial \Omega \times (0,\infty) \\
u_{0}(0)=u^{0}_{0}(x), &\hbox{in} \ \Omega
\end{array} \right.
\end{equation} 
whenever the initial conditions $u^{0}_{\eps}$ ``converge'' to $u^{0}_{0}$. The function $\gamma\in L^\infty(\partial\Omega)$  is related to the behavior of the $(N-1)$-dimensional measure of  $\partial\Omega_\eps$ as the parameter $\eps\to 0$ and captures the oscillatory behavior of the boundary. We refer to Subsection \ref{perdominio} for a complete and thorough definition of the domains and in particular the function $\gamma(\cdot)$, see Definition \ref{definition-of-gamma}.   Indeed, we will prove the existence and continuity of the family of attractors of (\ref{nbc}) and (\ref{nbc_limite_gamma_F}) in $H^1(\Omega_\eps)$ at $\epsilon =0$. 

\par\medskip 

One of the main difficulties when treating problems which are posed in different domains, like \eqref{nbc} and \eqref{nbc_limite_gamma_F}, is that
the solutions live in different spaces, say $H^1(\Omega_\eps)$  and $H^1(\Omega)$ or $L^p(\Omega_\eps)$ and $L^p(\Omega)$, see for instance \cite{AB0,AB1,AC,ACL1,ACL2,ACL3}.  So it is necessary to devise a tool to compare functions which are defined in different spaces and give a meaning to statements like $u_\eps\in H^1(\Omega_\eps)$ ``converges'' to $u\in H^1(\Omega)$. 

\par\medskip 

For this, consider the linear operator $E_\eps:H^1(\Omega)\to H^1(\Omega_\eps)$, which is defined as 
\begin{equation} \label{definitionE_epsilon}
E_\epsilon = R_\epsilon \circ P,
\end{equation} 
where $P:H^1(\Omega) \to H^1(\R^N)$ is a linear and continuous operator that extends a function $u$ defined in $\Omega$ to a function defined in $\R^N$, see for instance \cite{brezis},  and $R_\epsilon$ is the restriction operator from functions defined in $\R^N$ to functions defined in $\Omega_\eps$, that is, $R_\epsilon w=w_{|\Omega_\epsilon}$. Observe that we also have $E_\eps:L^p(\Omega)\to L^p(\Omega_\eps)$ and  $E_\eps:W^{1,p}(\Omega)\to W^{1,p}(\Omega_\eps)$, for all $1\leq p\leq \infty$. As a matter of fact, denoting by   $Z_\eps=H^1(\Omega_\eps)$ or $L^p(\Omega_\eps)$ or $W^{1,p}(\Omega_\eps)$, for $\eps \geq 0$, we have $E_\eps:Z_0\to Z_\eps$. By  \cite{AB0} we obtain
\begin{equation}
\label{extensaobase} 
\|E_\epsilon u\|_{Z_\epsilon} \to \|u\|_{Z_0}, \quad \mbox{as $\epsilon \to 0$,} \quad \mbox{and} \quad \|E_\epsilon\| \leq \|R_\eps\| \|P\|\leq \|P\|, \quad \mbox{independent of $\epsilon$}.
\end{equation}

With this operator we can define the following concept of convergence:

\begin{definition} 
\label{Def_E_convergence_A_S}
A family of elements $\{u_\eps\}_{\eps\in (0,\eps_0]}$, $u_\eps\in H^{1}(\Omega_\eps)$, is said to be $E$-convergent to $u\in H^{1}(\Omega)$   if  $\|u_\eps-E_\eps u\|_{H^1(\Omega_{\epsilon})}\to 0$ as $\eps\to 0$. We write this as $u_\eps\Eto u$ in $H^1(\Omega_\epsilon)$.  (Similarly we have the definition of $u_\eps$ E-converges to $u$ in $L^p(\Omega_\epsilon)$, $W^{1,p}(\Omega_\epsilon)$, etc).
\end{definition}

The $E$-convergence of the solutions of nonlinear elliptic equations with nonlinear boundary conditions and rapidly varying boundaries  
was studied in \cite{AB0} for this case of uniformly Lipschitz deformation of the boundary. The authors proved the $E$-convergence of the solutions in $H^1(\Omega_\epsilon)$. %and $C^{\alpha}(\Omega_\epsilon)$, for some $0<\alpha\leq 1$. 
In particular, if we regard these nonlinear elliptic equations as stationary equations of the parabolic evolutionary equations (\ref{nbc}) and (\ref{nbc_limite_gamma_F}), then\cite{AB0} proves the $E$-continuity of the set of equilibria of  (\ref{nbc}) in $H^1(\Omega_\epsilon)$ at $\epsilon =0$. Moreover, the authors proved the convergence of the eigenvalues and the $E$-convergence of the eigenfunctions of the linearization around the equilibrium points.  Also, they showed that if a equilibrium of the limit equation is hyperbolic, then the perturbed equation has one and only one equilibrium nearby. The goal of this work is to continue the analysis initiated in \cite{AB0} and address the complete analysis of the asymptotic dynamics, showing that, under certain conditions of the nonlinearities and the domains,  the attractor of \eqref{nbc} $E$-converges in $H^1(\Omega_\epsilon)$ to the attractor of \eqref{nbc_limite_gamma_F}. 

\par\medskip 

 As a matter of fact, with the concept of $E$-convergence we consider the following definition of $E$-continuity of a family of sets
\begin{definition} \label{def_family_E_continuity}
Let $\mathcal{A}_{\epsilon}\subset H^1(\Omega_{\epsilon})$, $\epsilon\in (0,\epsilon_0]$, and  $\mathcal{A}_{0}\subset H^1(\Omega)$. We have:

\vspace{0.2cm}

\noindent (i) We say that the family $\left\{\mathcal{A}_{\epsilon}\right\}_{\epsilon\in [0,\epsilon_0]}$ is $E$-upper semicontinuous in $H^1(\Omega_\epsilon)$ at $\epsilon=0$ if 
$$
dist\left(\mathcal{A}_{\epsilon},E_\epsilon \mathcal{A}_{0}\right):=\sup_{x_{\epsilon}\in \mathcal{A}_{\epsilon}} dist\left(x_{\epsilon},E_\epsilon \mathcal{A}_{0}\right)=\sup_{x_{\epsilon}\in \mathcal{A}_{\epsilon}} \inf_{x\in \mathcal{A}_{0}} \left\|x_{\epsilon}-E_\epsilon x\right\|_{H^1(\Omega_\epsilon)}\to 0, \quad \mbox{as $\epsilon\to 0$}.
$$
\noindent (ii) We say that the family $\left\{\mathcal{A}_{\epsilon}\right\}_{\epsilon\in [0,\epsilon_0]}$ is $E$-lower semicontinuous in $H^1(\Omega_\epsilon)$ at $\epsilon=0$ if 
$$
dist\left(E_\epsilon \mathcal{A}_{0},\mathcal{A}_{\epsilon}\right):=\sup_{x\in \mathcal{A}_{0}} dist\left(E_{\epsilon}x,\mathcal{A}_{\epsilon}\right)=\sup_{x\in \mathcal{A}_{0}} \inf_{x_{\epsilon}\in \mathcal{A}_{\epsilon}} \left\|x_{\epsilon}-E_{\epsilon}x\right\|_{H^1(\Omega_\epsilon)}\to 0, \quad \mbox{as $\epsilon\to 0$}.
$$
\noindent  (iii) We say that the family $\left\{\mathcal{A}_{\epsilon}\right\}_{\epsilon\in [0,\epsilon_0]}$ is $E$-continuous in $H^1(\Omega_\epsilon)$ at $\epsilon=0$ if $\left\{\mathcal{A}_{\epsilon}\right\}_{\epsilon\in [0,\epsilon_0]}$ is $E$-upper and $E$-lower semicontinuous in $H^1(\Omega_\epsilon)$ at $\epsilon=0$.
\end{definition}

 Let us remark that a similar definition can be stated for $E$-upper semicontinuity (lower semicontinuity and continuity) in $L^p(\Omega_\eps)$ or $W^{1,p}(\Omega_\eps)$.

\par\medskip  With this definition, we can state our main result of this paper, which shows the existence  of the attractors of \eqref{nbc} and \eqref{nbc_limite_gamma_F}, and their convergence  as the parameter $\eps\to 0$,

\begin{theorem}
\label{existence-upper-lower-attractors}
Assume that the family of domains $\{\Omega_\eps\}_{\eps\in [0,\eps_0]}$ and the nonlinearities $f$ and $g$ satisfy  the appropriate conditions stated below in Subsection \ref{perdominio} and Subsection \ref{equation-setting}. Then, for each $\eps\in (0,\eps_0]$, problem \eqref{nbc} has an attractor $\mathcal{A}_\eps\subset H^1(\Omega_\eps)$ and \eqref{nbc_limite_gamma_F} has an attractor $\mathcal{A}_0\subset H^1(\Omega)$. Moreover, we have:  

\vspace{0.2cm}

\noindent (i) The family of attractors $\{\mathcal{A}_{\epsilon}\}_{\epsilon \in [0,\epsilon_0]}$ is $E$-upper semicontinuous in $H^1(\Omega_\epsilon)$ at $\epsilon =0$.

\vspace{0.2cm}

\noindent (ii) If every equilibrium of the limit problem \eqref{nbc_limite_gamma_F} is hyperbolic, then the family of attractors $\{\mathcal{A}_{\epsilon}\}_{\epsilon \in [0,\epsilon_0]}$ is also $E$-lower semicontinuous in $H^1(\Omega_\epsilon)$ at $\epsilon =0$. In particular, the attractors are $E$-continous in $H^1(\Omega_\eps)$ at $\eps=0$. 
\end{theorem}

In order to prove this result, we will use the results from Carvalho-Piskarev \cite{CP}.  In this paper, the authors present an abstract result on the continuity of attractors in a similar setting as our case. Nevertheless the results from \cite{CP} cannot be applied directly to prove Theorem \ref{existence-upper-lower-attractors}, but they need to be adapted properly to our situation. This adaptation is not immediate and it deserves some rather technical results.  

\par\medskip 

One key point to adapt the results from \cite{CP} to our case is to extend appropriately the definition of extension operators and $E$-convergence to the case of fractional spaces with negative exponents. We dedicate Subsection  \ref{extensionoperators} to this issue. Another key point is to prove that the notion of $E$-convergence defined above is equivalent to the notion of convergence stablished in an abstract way in \cite{CP} and that the authors denote as $\cal{P}$-convergence.

\par\medskip 

Related to this work, we would like to emphasize that in \cite{pricila} the authors study the continuity of the family of attractors associated to semilinear parabolic problems with nonlinear Neumann boundary conditions, where the domain $\Omega \subset \mathbb{R}^2$ is a unit square and with perturbations of the domain which are small in the $C^1$-sense. It is important to note that the technique used in \cite{pricila} is different from the one used here. Actually, in this paper they follow the approach developed by D. Henry in  \cite{henry_perturbation} to deal with domain perturbations problems. This approach is based in considering the perturbed domains $\Omega_\eps=h_\eps(\Omega)$ where $h_\eps$ is a diffeomorphims which is close to the identity in $C^k$ norms for some $k\geq 1$.   Moreover, the technique of \cite{henry_perturbation}  was also applied in \cite{marconetoninho1,marconetoninho2}.  In these nice works, the perturbed domain approaches the unperturbed one in a $C^k$-sense for some $k\geq 1$,  which is a stronger convergence than the one we have in the present work. Also, because of this $C^1$ convergence the limiting boundary condition does not have the amplification factor $\gamma(x)$, as we have in our limiting equation \eqref{nbc_limite_gamma_F}. 
\par\medskip 

 Recently, in \cite{Pires1} the authors use compact convergence techniques to deal with the continuity of attractors of some reaction-diffusion equations under smooth perturbations of the domain subject to nonlinear Neumann boundary conditions. They also follow the setting  in \cite{henry_perturbation} and define a family of invertible linear operators $E_h:H^{1}(\Omega)\to H^{1}(\Omega_h)$ to compare the dynamics of perturbed and unperturbed problems in the same phase space, where $h:\Omega\to \mathbb{R}^N$ is a diffeomorphism which approaches
the inclusion $I_N :\Omega \to \mathbb{R}^N$ in the $\mathcal{C}^1$ topology and $\Omega_h=h(\Omega)$. Again, the perturbations are $C^1$ small and the factor $\gamma(x)$ do not appear in the limiting equation.  On the other hand, they obtain some rates on the convergence of the attractors following similar ideas as the ones from \cite{ArrietaSantamaria}.  
In \cite{Pires2} the authors treat the case of Dirichlet boundary condtions. 

\par\medskip 

This paper is organized as follow: in Section \ref{abstractsetting0} we give the precise hypotheses and definitions about the domain perturbation and nonlinearitites and we prove abstract results about embeddings and traces of fractional power space. We describe our approach to prove the continuity of the attractors, which consists in applying the setting of   \cite{CP}. Moreover, we extend the 
concepts of extension operators and convergence to
positive fractional power spaces and
we define extension operators and  convergence in spaces with negative exponents.  In Section \ref{P_continuity_attractors} we check the hypotheses of \cite{CP}, related to either compact convergence of the resolvent operators or convergence of the nonlinearities,
 to show the $\mathcal{P}$-continuity of the attractors. In Section \ref{continuity_attractors_H1} we prove that the notion of convergence established in an abstract way in the paper of \cite{CP} is the same as the notion of $E$-convergence we are using in this paper and we conclude
the $E$-continuity of the attractors.

\par\bigskip

\noindent {\bf Acknowledgements.}  

G.. S. Arag\~ao is partially supported by FAPESP 2020/14075-6, Brazil.

J. M. Arrieta is partially supported by grants PID2019-103860GB-I00, PID2022-137074NB-I00  and CEX2019-000904-S ``Severo Ochoa Programme for Centres of Excellence in R\&D'' the three of them from MICINN, Spain. Also by ``Grupo de Investigaci\'on 920894 - CADEDIF'', UCM, Spain.

%%%%%%%%%%%%%%%%%%%%%%%%%%%%%%%%%%%%%%%%%%%%%%%%%%%%%%%%%%%%%%%%%%%%%%%%%%%%%%%%%%%

%\section{Abstract setting}
\section{General setting and main results}
\label{abstractsetting0}

In this section we clarify the general setting of the problem. We provide hypotheses, definitions and results for the domain perturbation, nonlinearities,  abstract form,  and known  results  as well as a description of our approach and some results to implement  the appoach. 

More precisely, 
in Subsection \ref{perdominio} we describe in detail the domain perturbation that we are considering, providing the appropriate definitions and hypotheses. In Subsection \ref{equation-setting}   we write our problem in an abstract form and describe the hypothesis on the nonlinearities. In Subsection \ref{known-results} we review known results for this problem.  These  are  basically taken from  \cite{Arrieta95,AB0,AC}. In Subsection \ref{Carvalho-Piskarev} we describe our approach to prove the convergence of the attractors, which consists in applying the setting of Carvalho-Piskarev  \cite{CP} to problems (\ref{nbc}) and (\ref{nbc_limite_gamma_F}). 
In Subsection \ref{power space} we  study the relation between the fractional power spaces and other spaces like $L^p(\Omega_\eps)$ and $L^p(\partial \Omega_\eps)$ spaces. We pay special attention to the inclusion properties among these and other spaces and obtain that the embedding constants can be chosen uniformly for the whole family of domains $\{\Omega_\eps\}_{\eps\in [0,\eps_0]}$.  This is a very crucial step.  In Subsection \ref{extensionoperators} we extend the 
concepts of extension operators and convergence to
 fractional power spaces with positive and negative exponents.

%%%%%%%%%%%%%%%%%%%%%%%%%%%%%%%%%%%%%%%%%%%%%%%%%%%%%%%%%%%%%%%%%%%%%%%%%%%%%%%%%%%

\subsection{Setting of the perturbation of the domain}
\label{perdominio}

We consider a family of uniformly bounded smooth domains $\Omega_\epsilon \subset \mathbb{R}^N$, with $N\geq 2$ and $0\leq\epsilon\leq \epsilon_0$, for some $\epsilon_0>0$ fixed, and we regard $\Omega_\epsilon$ as a perturbation of the fixed domain $\Omega \equiv \Omega_0$. As well as in \cite{AB0}, we consider the following hypothesis

\begin{description}
\item[(H)]\textbf{(i)} For all $K\subset \Omega$, $K$ being compact, there exists $\epsilon(K)>0$ such that $K\subset \Omega_\epsilon$ for $0<\epsilon<\epsilon(K)$.

\textbf{(ii)} There exists a finite open cover $\{U_i\}_{i=0}^m$ of $\Omega$ such that $\overline{U}_0 \subset \Omega$, $\partial\Omega \subset \cup_{i=1}^m U_i$ and for each $i=1, \ldots, m$, there exists a Lipschitz diffeomorphism $\Phi_{i}: Q_N \to U_i$, where $Q_N=(-1,1)^N\subset \R^N$, such that
$$
\Phi_{i}(Q_{N-1} \times (-1,0)) = U_i \cap \Omega \qquad \mbox{and} \qquad \Phi_{i}(Q_{N-1} \times \{0\}) = U_i \cap \partial \Omega.
$$
We assume that $\overline{\Omega}_\epsilon \subset \cup_{i=0}^m U_i\equiv U$, and for each $i=1, \ldots, m$, there exists a Lipschitz  function  $\rho_{i,\eps}: Q_{N-1} \to (-1,1)$ such that $\rho_{i,\epsilon} \to 0$ as $\epsilon \to 0$, uniformly in $Q_{N-1}$, for each $i=1,\ldots,m$.

Moreover, we assume that $\Phi_{i}^{-1}(U_i \cap \partial \Omega_\epsilon)$ is the graph of $\rho_{i,\epsilon}$ this means 
$$
U_i \cap \partial \Omega_\epsilon=\Phi_{i}(\{(x',\rho_{i,\epsilon}(x')) \;:\; x'=(x_1,...,x_{N-1}) \in  Q_{N-1}\}).
$$
\end{description}

Note that if $\Omega \subset \Omega_\epsilon$, that is, $\Omega_\epsilon$ is an exterior perturbation of $\Omega$, then condition \textbf{(H)(i)} is satisfied.

\begin{figure}[h]
\begin{center}
\includegraphics[width=.7\linewidth,height=0.6\textheight,keepaspectratio]{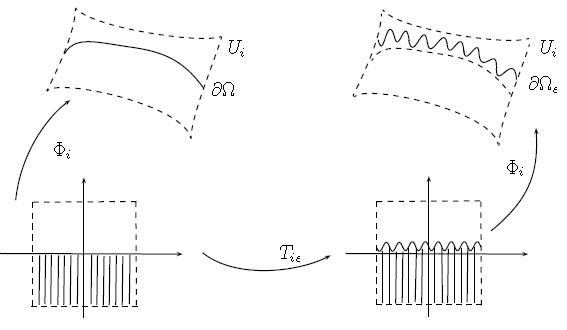}
\vspace*{0cm}
\caption{The parametrizations.}
\label{figure2}
\end{center}
\end{figure}

We consider the following  mappings $T_{i,\epsilon}:Q_N \to Q_N$ defined by
$$
T_{i,\epsilon}(x',s) =
\left\{
\begin{array}{ll}
(x',s+s\rho_{i,\epsilon}(x')+\rho_{i,\epsilon}(x')),&\hbox{ for $s \in (-1,0)$ and $x'\in Q_{N-1}$}\\
(x',s-s\rho_{i,\epsilon}(x')+\rho_{i,\epsilon}(x')),&\hbox{ for $s \in [0,1)$ and $x'\in Q_{N-1}$}.
\end{array}
\right.
$$
Also, 
$$
\Phi_{i,\epsilon}:=\Phi_{i}\circ T_{i,\epsilon}:Q_N \to U_i, 
$$
and we also denote by 
$$
\begin{array}{rl}
\phi_{i,\epsilon}:Q_{N-1}&\to U_i\cap\partial\Omega_\epsilon \\
x'& \mapsto \Phi_{i,\eps}(x',0)
\end{array} \qquad \mbox{and} \qquad
\begin{array}{rl}
\phi_{i}:Q_{N-1}&\to U_i\cap\partial\Omega \\
x'& \mapsto \Phi_{i}(x',0).
\end{array}
$$
Notice that $\phi_{i,\epsilon}$ and $\phi_{i}$ are local parametrizations of $\partial\Omega_\eps$ and $\partial\Omega$, respectively. Furthermore, observe that all the maps above are Lipschitz. Figure \ref{figure2} illustrates the parametrizations.

\vspace{6cm}

In order to give the hypothesis to deal with the deformation $\partial\Omega_\epsilon$ we need the following definition

\begin{definition}
\label{definition-jacobian}
Let $\eta: A \subset \R^{N-1} \to \R^N$ almost everywhere differentiable,  we define  
the $(N-1)$-dimensional Jacobian of $\eta$  as
%\begin{equation} 
%\label{defN-1Jacobian}
$$
J_{N-1}\eta \equiv \left|{\partial \eta \over \partial x_1}\wedge\ldots\wedge {\partial \eta\over \partial x_{N-1}}\right| = \sqrt{\sum_{j=1}^N (\mbox{det}({\rm{Jac}}\,\eta)_j)^2},
$$
%\end{equation}
where $v_1 \wedge\ldots\wedge v_{N-1}$ is the exterior product of the $(N-1)$ vectors $v_1,\ldots,v_{N-1} \in \mathbb{R}^{N}$
and $({\rm{Jac}}\,\eta)_j$ is the $(N-1)$-dimensional matrix obtained by deleting the $j$-th row of the Jacobian matrix of $\eta$.  
\end{definition}

We use $J_N$ for the absolute value of the $N$-dimensional Jacobian determinant. 

Now, to deal with the interaction between the nonlinear boundary condition and the oscillatory behavior of $\partial\Omega_\eps$  we consider the hypothesis

\begin{description}
\item[(F)]\textbf{(i)} $\| \nabla \rho_{i,\epsilon}\|_{L^\infty(Q_{N-1})}\leq C$, with $C>0$ independent of $\epsilon$, $i=1,\ldots,m$.

\textbf{(ii)} For each $i=1,\ldots,m$, there exists a function $\gamma_i \in L^\infty(Q_{N-1})$ such that
%\begin{equation} 
%\label{defgamma}
$$J_{N-1}\phi_{i,\eps} \stackrel{\epsilon \to 0}{-{\hspace{-2mm}}\rightharpoonup} \gamma_i \quad \mbox{in $L^1(Q_{N-1})$}.$$
%\end{equation}
\end{description}

%\stackrel{\epsilon \to 0}{-{\hspace{-2mm}}\rightharpoonup}
%\vspace{0.2cm}

Considering hypotheses {\bf(H)} and {\bf(F)}, let  $\gamma:  \partial\Omega \to \R$ be a function which measures the limit of the deformation of $\partial \Omega_\epsilon$ relatively to $\partial\Omega$. More precisely, we have

\begin{definition}
\label{definition-of-gamma}
For $x\in U_i\cap\partial\Omega$, let $(x',0) =\Phi_i^{-1} (x) \in Q_{N}$, we define $\gamma:  \partial\Omega \to \R$ as
$$
\gamma(x) = \frac{\gamma_i(x')}{J_{N-1}\phi_i(x')}.$$
\end{definition}

The function $\gamma$ is independent of $\phi_{i,\epsilon}$ and also on the choice of the charts $U_i$ and the maps $\Phi_i$, that is, $\gamma$ is independent of the parameterization
 chosen, and thefore it is unique. Moreover $\gamma \geq 1$. This was proved in \cite[Corollary 5.1]{AB0}. 

In order to deal also with non exterior perturbations, we  consider a family $K_\epsilon$ of smooth interior perturbations of $\Omega$ satisfying the  following hypothesis 

\begin{description}
\item[(I)]
For each $\epsilon \in (0,\epsilon_0]$, for some $\epsilon_0>0$ fixed, there exists $K_\epsilon \subset \Omega \cap\Omega_\epsilon$ such that $U_0 \subset K_{\epsilon_0} $, 
$K_{\epsilon_1} \subset K_{\epsilon_2}$ if $\epsilon_1 > \epsilon_2$, 
 $\theta_\epsilon :  \Omega \to K_\epsilon $ is a diffeomorphism such that ${ \theta_\epsilon}_{|K_{\epsilon_0}} $ is the identity in $K_{\epsilon_0}$, $\partial K_\epsilon$ is diffeomorphic to $\partial \Omega$ and $D \theta_\epsilon\ $ converges to $I$  in $L^\infty(\Omega, {\cal{L}}(\R^N))$ as $\epsilon \to 0$ and, for each $i=1,...,m$, there exists $\hat{\rho}_{i,\epsilon} : Q_{N-1} \to (-1,1)$ such that 
 $$
U_i\cap \partial K_\epsilon =\Phi_i  (\{(x^\prime, \hat{\rho}_{i,\epsilon}(x^\prime)): \; x^\prime =(x_1,...,x_{N-1}) \in Q_{N-1}\}).
 $$
%$$
%U_i\cap \partial K_\epsilon = \{\hat{ \phi}_{i,\epsilon}(x^\prime) = \Phi_i(x^\prime, \hat{\rho}_{i,\epsilon}(x^\prime)): \; x^\prime =(x_1,...,x_{N-1}) \in Q_{N-1}\}.
 %$$
 
We also suppose that $J_{N-1} \hat{ \phi}_{i,\epsilon} \stackrel{\epsilon \to 0}{-{\hspace{-2mm}}\rightharpoonup} 1$ in $L^1(Q_{N-1})$, where $\hat{ \phi}_{i,\epsilon}(x^\prime) = \Phi_i(x^\prime, \hat{\rho}_{i,\epsilon}(x^\prime))$ for $x^\prime \in Q_{N-1}$, and 
$$\int_{\Omega \setminus K_\epsilon} u = \int_0^\epsilon \int_{\partial K_\sigma} u \  d\sigma.$$
\end{description}

If $\Omega_\epsilon$ is an exterior perturbation of $\Omega$, then we can consider $K_\epsilon =\Omega$ and this hypothesis is satisfied.

We observe that if $\Omega$ is a $C^2$ domain then it satisfies this hypothesis, see \cite{AJCRB, henry_perturbation}.  
There are also some examples of $C^{0,1}$ domains that satisfies this hypothesis, for examples squares with particular perturbations, see \cite{pricila}.

%%%%%%%%%%%%%%%%%%%%%%%%%%%%%%%%%%%%%%%%%%%%%%%%%%%%%%%%%%%%%%%%%%%%%%%%%%%%%%

\subsection{Setting of the equations and hypotheses on the nonlinearities}
\label{equation-setting}

For $0\leq\eps\leq\eps_0$, consider the linear operator $A_\epsilon: D(A_\epsilon) \subset L^2(\Omega_\epsilon) \to
L^2(\Omega_\epsilon)$ defined by 
$$
A_\epsilon u_\epsilon = -\Delta u_\epsilon + u_\epsilon
$$
with domain 
$$
D(A_\epsilon) = \left\{u_\epsilon \in H^2(\Omega_\epsilon): \frac{\partial u_\epsilon}{\partial n_\epsilon}=0  \mbox{\ in\ } \partial\Omega_\eps\right\},
$$
where  we identify $\Omega\equiv\Omega_0$. Let us denote by $X^0_\epsilon = L^2(\Omega_\epsilon)$ and $X^1_\epsilon = D(A_\epsilon)$, endowed with the graph norm. Since this operator turns out to be sectorial in $X^0_\epsilon$, associated to it there is a scale of Banach spaces (the fractional power spaces)  $X^\alpha_\epsilon$, $\alpha \geq 0$, denoting the domain of the fractional power operators associated with $A_\epsilon$, that is, $X^\alpha_\epsilon:=D(A^{\alpha}_\epsilon)$ for $\alpha\geq 0$, where $X^\alpha_\epsilon$ endowed with the graph norm $\|x\|_{X^\alpha_\epsilon}=\|A^{\alpha}_\epsilon x\|_{X^0_\epsilon}$, $\alpha\geq 0$, with $X^{\half}_\epsilon = H^{1}(\Omega_\epsilon)$, see \cite{henry}. Since we are going to deal with a nonlinear boundary conditions problem, we also consider spaces of negative exponents by taking $X^{-\alpha}_\epsilon = (X^{\alpha}_\epsilon)^\prime$, for $\alpha >0$, with $X^{-\half}_\epsilon = (H^{1}(\Omega_\epsilon))^\prime$ denoted by $H^{-1}(\Omega_\epsilon)$. 

Denoting by $A_{\epsilon,\alpha}$, $\alpha \in \mathbb{R}$,  the realizations of $A_\epsilon$ in this scale we have that the operator $A_{\epsilon,-\half} \in \mathcal{L}(X^{\half}_\epsilon,X^{-\half}_\epsilon)$ is given by 
$$
\langle  A_{\epsilon,-\half}u_\epsilon, v_\epsilon\rangle = \int_{\Omega_\epsilon} (\nabla u_\epsilon \nabla v_\epsilon + u_\epsilon v_\epsilon), \quad \mbox{for $u_\epsilon, v_\epsilon \in X^{\frac{1}{2}}_{\epsilon}$}.
$$
With some abuse of notation we will identify all different realizations of this operator and we will write them all as $A_\epsilon$.  Moreover, we will denote $X^{\frac{\alpha}{2}}\equiv X^{\frac{\alpha}{2}}_{0}$ and $X^{-\frac{\alpha}{2}}\equiv X^{-\frac{\alpha}{2}}_{0}$.

\par\medskip

We can rewrite (\ref{nbc}) and (\ref{nbc_limite_gamma_F}) in an abstract form as 
\begin{equation}
\label{nbcabstract}
\left\{
\begin{array}{lll}
\dot{u}_\epsilon(t)+ A_{\epsilon} u_\epsilon(t) = h_\epsilon(u_\epsilon(t)), & t>0 \\
u_\epsilon(0)=u^0_\epsilon \in X^{\frac{1}{2}}_\epsilon,
\end{array} \right.
\end{equation}
where $h_\epsilon : X^{\frac{1}{2}}_\epsilon \to X^{-\frac{\alpha}{2}}_\epsilon$, with $0<\epsilon\leq\eps_0$ and  $\half <\alpha \leq 1$, is defined by
\begin{equation} 
\label{hepsilon}
\langle h_\epsilon(u_\epsilon), \psi_\epsilon\rangle = \int_{\Omega_\epsilon} f(x,u_\epsilon)\psi_\epsilon + \int_{\partial\Omega_\epsilon}g(x,u_\epsilon) \psi_\epsilon,  \quad \mbox{for $u_\epsilon \in X^{\frac{1}{2}}_\epsilon$ and $\psi_\epsilon \in X^{\frac{\alpha}{2}}_\epsilon$}.
\end{equation}
And $h_0: X^{\frac{1}{2}} \to X^{-\frac{\alpha}{2}}$, with $\half <\alpha \leq 1$, is defined by 
\begin{equation} 
\label{hzero}
\langle h_0(u), \psi\rangle = \int_{\Omega} f(x,u)\psi +\int_{\partial\Omega}\gamma(x) g(x,u) \psi,  \quad \mbox{for $u \in X^{\frac{1}{2}}$ and $\psi \in X^{\frac{\alpha}{2}}$}.
\end{equation}

\par\medskip  With respect to the nonlinearities $f:U\times \R\to \R$ and $g:U\times \R\to \R$ we will assume they are continuous in both variables,
$C^2$ in the second one
and satisfy
\begin{equation}
\label{bounded_f}
\left|f(x,u)\right|+\left|\partial_u f(x,u)\right|+\left|\partial_{uu} f(x,u)\right|\leq C,\quad \mbox{for all $x\in U$ and $u\in \mathbb{R}$},
\end{equation}
\begin{equation}
\label{bounded_g}
\left|g(x,u)\right|+\left|\partial_u g(x,u)\right|+\left|\partial_{uu} g(x,u)\right|\leq C,\quad \mbox{for all $x\in U$ and $u\in \mathbb{R}$}.
\end{equation}

\vspace{0.2cm}

\begin{remark}\label{remark-conditions-f-g}
Although conditions \eqref{bounded_f} and \eqref{bounded_g} do not look rather general, note that if 
$f$ and $g$ are functions satisfying appropriate growth and
sign conditions, we can prove  global existence and uniqueness of mild solutions of (\ref{nbc}) and (\ref{nbc_limite_gamma_F}) or \eqref{nbcabstract}, for $u^{0}_{\epsilon} \in X^{\frac{1}{2}}_{\epsilon}$, see  \cite[Theorem 1.2]{ACRcritical} and \cite[Theorems 2.1 and 2.2]{ACB}.  Assuming also  a dissipativeness 
condition
and proceeding in a similar way to
\cite[Proposition 3.2]{ACB}, we can obtain that the solutions of the equations
(\ref{nbc}) and (\ref{nbc_limite_gamma_F})  are bounded in $L^{\infty}(\Omega_{\epsilon})$, uniformly in $\epsilon$, so we may perform a cut off in the nonlinearities $f$ and $g$ in such a way that they become bounded with bounded derivatives up to second order without changing the solutions of the equations in a large bounded set in $L^\infty(\Omega_{\epsilon})$.  After these considerations, without loss of generality, we can assume that the nonlinearities satisfy \eqref{bounded_f} and \eqref{bounded_g}.
\end{remark}

%%%%%%%%%%%%%%%%%%%%%%%%%%%%%%%%%%%%%%%%%%%%%%%%%%%%%%%%%%%%%%%%%%%%%%%%%%%%%%

\subsection{Known results}
\label{known-results}

In this subsection we include some important definitions on $E$-convergence and we summarize some results on spectral  convergence, convergence of equilibria and their linearizations which  are taken mainly from \cite{Arrieta95,AB0,AC}.

\par\bigskip

\noindent {\bf Some definitions related to $E$-convergence.}  We include here some definitions that we will use throughout the paper.  We recall the definition of  the concept of compactness and of convergence of a family of  operators $T_\epsilon:Z_\eps\to W_\eps$ acting between two families of Banach spaces $W_\eps$ and $Z_\eps$, $\eps\in (0,\eps_0]$. For each of this family we have the appropriate family of ``extension'' operators $E_\eps^W:W\to W_\eps$ and $E_\eps^Z:Z\to Z_\eps$. Notice that with some abuse of notation we will denote by $E_\eps$ the extension operators for both families of Banach spaces.  
These concepts follow \cite{AB0}. We also refer to \cite{AB0,ACL1,ACL3,CCS} and references therein, for a detailed study of these notions  and its applications to differential equations. 
%Let $Y_\epsilon$, $\epsilon \in (0,\epsilon_0]$, and $Y$ be Banach spaces, we have

\begin{definition} \label{defEprecompact}
A sequence of elements $\{u_n\}_{n\in \mathbb{N}}$, with $u_n\in Z_{\epsilon_n}$ and $\eps_n\to 0$, is said to be $E$-precompact if for any subsequence $\{u_{n'}\}$ there exist a subsequence $\{u_{n''}\}$ and $u\in Z$ such that $u_{n''}\dto u$ as $n''\to \infty$. A family $\{u_\epsilon\}_{\epsilon\in (0,\epsilon_0]}$, $u_\epsilon\in Z_{\epsilon}$, is said to be $E$-precompact if each sequence $\{u_{\epsilon_n}\}_{n\in \mathbb{N}}$, with $\epsilon_n \to 0$, is $E$-precompact. 
\end{definition}

\begin{definition} \label{defconvof operators} We say that a family of
operators $T_\epsilon:Z_\epsilon\to W_{\epsilon}$, $\epsilon \in (0,\epsilon_0]$,
$E$-converges to $T:Z\to W$ as $\epsilon \to 0$, if $T_\epsilon
u_\epsilon \dto Tu\in W$ , whenever $u_\epsilon\dto u \in Z$. We denote
this by $T_\epsilon \ddto T$.
\end{definition}

\begin{definition} \label{defcompactconvoperators}
We say that a family of compact operators
$T_\epsilon:Z_\epsilon \to W_\epsilon$, $\epsilon \in (0,\epsilon_0]$, converges
compactly to a compact operator $T:Z\to W$ if for any  family
$\{u_\epsilon\}_{\epsilon\in (0,\epsilon_0]}$ with $\|u_\epsilon\|_{Z_{\epsilon}}$ bounded, the family $\{T_\epsilon u_\epsilon\}_{\eps\in (0,\eps_0]}$ is $E$-precompact and $T_\epsilon \ddto T$. We denote this by $T_\epsilon \ccto T$.
\end{definition}

\par\bigskip 

\noindent {\bf Spectral convergence.} The kind of domain perturbation that we are considering and that has been detailed in Subsection \ref{perdominio} guarantees the spectral convergence of the Laplace operator with homogeneous Neumann boundary conditions in $\Omega_\eps$ to the same operator in the limiting domain $\Omega$.  We refer to 
\cite{Arrieta95,AC} for this result.  In particular, this means that if we denote by $\mu_n^\eps$ the eigenvalues of the operator $-\Delta +I$ in $\Omega_\eps$ with homogeneous Neumann boundary conditions and by $\mu_n$ the eigenvalues in $\Omega$ also with Neumann boundary conditions, then we have $\mu_n^\eps\to \mu_n$ as $\eps\to 0$. A similar statement is obtained for the eigenfunctions.

\par\bigskip 

\noindent {\bf Equilibria behavior.}  The first step in the proof of the continuity of the attractors is to study the simplest elements from the attractor, the equilibrium solutions. The equilibrium solutions of (\ref{nbc}) and (\ref{nbc_limite_gamma_F}) are those solutions which are independent of time and therefore they are the solutions of the following nonlinear elliptic problems
\begin{equation}
\label{eql1}
\left\{
\begin{array}{lll}
-\Delta e_{\epsilon}+e_{\epsilon}=\displaystyle  f(x,e_{\epsilon}), &\hbox{in} \ \Omega_\eps  \\
\frac{\ts\partial e_{\epsilon}}{\ts\partial n_{\epsilon}}=g(x,e_{\epsilon}), & \hbox{on} \ \partial \Omega_\eps 
\end{array} \right.
\end{equation}
and
\begin{equation}
\label{eql2}
\left\{
\begin{array}{lll}
-\Delta e_{0}+e_{0}=f(x,e_{0}), &\hbox{in} \ \Omega  \\
\frac{\ts\partial e_{0}}{\ts\partial n}=\gamma(x)g(x,e_{0}), & \hbox{on} \ \partial \Omega. 
\end{array} \right.
\end{equation}

For each $\epsilon\in [0,\epsilon_{0}]$, we denote by $\mathcal{E}_{\epsilon}\subset X^{\frac{1}{2}}_{\epsilon}=H^1(\Omega_\eps)$ the set of solutions of (\ref{eql1}) and (\ref{eql2}).  Notice first that from \eqref{bounded_f} and \eqref{bounded_g} we easily obtain the uniform boundedness in $X^{\frac{1}{2}}_{\epsilon}$ of all equilibria, that is, there exists a $C>0$ independent of $\eps\in [0,\eps_0]$ so that  $\|e_\eps\|_{X^{\frac{1}{2}}_{\epsilon}}\leq C$ for all $e_\eps\in \mathcal{E}_{\epsilon}$ and for all $\eps\in [0,\eps_0]$.    Moreover, in \cite{AB0}, by using $E$-convergence  in $X^{\frac{1}{2}}_{\epsilon}$, it was studied the $E$-continuity of the equilibrium points. It was proved the $E$-upper semicontinuity of the family $\mathcal{E}_\epsilon$  at $\epsilon =0$ and, by assuming the hyperbolicity of the equilibrium points in $\mathcal{E}_0$  (and therefore there are only a finite number of them) it was proved the $E$-lower semicontinuity at  $\epsilon =0$. Moreover, the authors also proved that there exist $0<\delta<1$ and $M>0$ such that $\|e_\epsilon\|_{C^\delta(\Omega_\epsilon)} \leq M$ and if  $\|e_\epsilon - E_\epsilon e_0\|_{X^{\frac{1}{2}}_{\epsilon}} \to 0$ then  $\|e_\epsilon - E_\epsilon e_0\|_{C^\beta(\Omega_\epsilon)} \to 0$ as $\epsilon \to 0$, $0 < \beta < \delta$, see \cite[Proposition 5.3]{AB0}. 

\par\bigskip 

\noindent {\bf Spectral convergence of the linearizations around the equilibrium solutions.} The spectra of the linearization of (\ref{nbc}) around $e_\eps^*$ equilibrium solution of (\ref{nbc})
is given by the eigenvalue problem
\begin{equation}  \label{nbc-linearized}
\left\{
\begin{array}{ll}-\Delta w_{\epsilon}+w_{\epsilon}-\partial_u f(x,e_\eps^*)w_{\epsilon}=\lambda^{\epsilon} w_{\epsilon}, &\hbox{in} \ \Omega_\eps \\
{\ts\partial w_{\epsilon} \over \ts\partial n_{\epsilon}}+\partial_ug(x,e_\eps^*)w_{\epsilon}=0, &\hbox{on} \
\partial \Omega_\eps.
\end{array} \right.
\end{equation}
Similarly, if $e_0^*$ is an equilibrium solution  of (\ref{nbc_limite_gamma_F}),
then the spectra of its linearization is given by the eigenvalue
problem
\begin{equation}  \label{nbc-limit-linearized}
\left\{
\begin{array}{ll}-\Delta w_0+w_0-\partial_u f(x,e_0^*)w_0=\lambda^0 w_0, &\hbox{in} \ \Omega \\
{\ts\partial w_0 \over \ts\partial n}+\gamma(x) \partial_u g(x,e_0^*)w_0=0, &\hbox{on} \
\partial \Omega.
\end{array} \right.
\end{equation}

Notice that both problems, (\ref{nbc-linearized}) and
(\ref{nbc-limit-linearized}), are selfadjoint and of compact
resolvent. Hence, the eigenvalues of (\ref{nbc-linearized}) are
given by a sequence $\{\lambda_n^\eps\}_{n=1}^\infty$, ordered and
counting their multiplicity,  with $\lambda_n^\eps\to \infty$ as
$n\to \infty$. Similarly the eigenvalues of
(\ref{nbc-limit-linearized}) are also given by a sequence
$\{\lambda_n^0\}_{n=1}^\infty$ with $\lambda_n^0\to\infty$ as
$n\to \infty$.

In \cite[Theorem 2.2]{AB0}, the authors proved if $\|e_\eps^*-E_\eps
e_0^*\|_{X^{\frac{1}{2}}_{\epsilon}}\to 0$ then the eigenvalues and
eigenfunctions of (\ref{nbc-linearized}) converge to the
eigenvalues and eigenfunctions of (\ref{nbc-limit-linearized}).
That is, for each fixed $n\in \N$, $\lambda_n^\eps\to
\lambda_n^0$, as $\eps\to 0$. Moreover, if we denote by
$\{\varphi_n^\eps\}_{n=1}^\infty$ a set of orthonormal
eigenfunctions associated to $\{\lambda_n^\eps\}_{n=1}^\infty$,
then for each sequence $\eps_k\to 0$ there is another subsequence,
that we still denote by $\eps_k$, and a set of orthonormal
eigenfunctions $\{\varphi_n^0\}_{n=1}^\infty$ associated to
$\{\lambda_n^0\}_{n=1}^\infty$ such that, for all $n\in
\N$, we have $\|\varphi_n^{\eps_k}-E_\eps
\varphi_n^0\|_{X^{\frac{1}{2}}_{\epsilon_k}}\to 0$ as $\eps_k\to 0$. 

Notice also that if $\lambda_n^0$ is a simple eigenvalue, then
$\lambda_n^\eps$ is also simple for $\eps$ small enough and, via
subsequences, we always have that $\varphi_n^{\eps_k}\Eto \varphi_n^0$
or
$\varphi_n^{\eps_k}\Eto -\varphi_n^0$ as $\eps_k\to 0$.

\subsection{Our approach to prove the convergence of the attractors}
\label{Carvalho-Piskarev}
As we have mentioned in the introduction, an important ingredient in our proof is the results from   \cite{CP}.  We recall now the setting, hypotheses and main results in this paper.

The authors consider an abstract semilinear parabolic evolution equation of the type
$\dot{u}=Au+h(u)$  with  $A:D(A)\subset Y\to Y$ a closed linear operator with compact resolvent in the Banach space $Y$. Moreover,  $-A$ is a sectorial operator satisfying 

$$
\|(\lambda I-A)^{-1}\|_{\mathcal{L} (Y)}\leq \frac{M}{1+|\lambda|}, \quad \mbox{for all $\lambda\in \C$ with $Re(\lambda)\geq 0$}. 
$$
We denote by $Y^0=Y$,  $Y^1=D(-A)$ endowed with the graph norm and $Y^\beta=D((-A)^{\beta})$, for $0<\beta<1$, the fractional power spaces endowed with the graph norm $\|x\|_{Y^{\beta}}=\|(-A)^{\beta}x\|_{Y}$. We assume there exists  $0\leq\beta<1$ such that $h:Y^\beta\to Y$ is a globally Lipschitz, bounded and continuously Fr\'echet differentiable function. 

They also consider a family of semilinear parabolic problems which are regarded as a perturbation of the above equation. As a matter of fact they consider a sequence but to adapt their statements to our setting, we will consider a one parameter family of problems indexed by $\eps$. For each $0< \epsilon \leq \epsilon_0$, we express these problems as 
\begin{equation*}
\left\{\begin{array}{lll} 
\dot{u}_\epsilon  + A_\epsilon u_\epsilon = h_\epsilon(u_\epsilon), & t>0 \\
u_\epsilon(0) = u^{0}_{\epsilon} \in Y_\epsilon^\beta,
\end{array}\right.
\end{equation*}
where  $Y_\epsilon$ are Banach spaces, $A_\eps:D(A_\eps)\subset Y_\eps\to Y_\eps$ are operators as above satisfying 
$$
\|(\lambda I-A_\eps)^{-1}\|_{\mathcal{L} (Y_\eps)}\leq \frac{M}{1+|\lambda|}, \quad \mbox{for all $\lambda\in \C$ with $Re(\lambda)\geq 0$ and 
 for all $0< \eps\leq \eps_0$,}
$$
$Y^0_\eps=Y_\eps$, $Y_\eps^1=D(-A_\eps),$ $Y_\eps^\beta$ are the fractional power spaces and $h_\eps:Y_\eps^\beta \to Y_\eps$, where $0\leq \beta<1$, are globally Lipschitz, bounded and continuously Fr\'echet differentiable functions.

In order to relate the unperturbed problem and the perturbed ones they assume that there exist  linear bounded operators $p_\eps:Y\to Y_\eps$ satisfying 
\begin{equation}
\label{conv_p_eps}
\|p_\eps y\|_{Y_\eps}\to \|y\|_{Y}, \quad \mbox{as $\eps\to 0$ for any $y\in Y$.}  
\end{equation}
With this sequence of linear bounded operators, they construct the  operators $p_\eps^{\beta}:Y^\beta\to Y^\beta_\eps$ given by
$$p_\eps^\beta:= (-A_\eps)^{-\beta}p_\eps (-A)^{\beta},
$$ 
which are linear bounded operators satisfying  
$$
\|p_\eps^\beta y\|_{Y_\eps^\beta}\to \|y\|_{Y^\beta}, \quad \mbox{as $\eps\to 0$  for any $y\in Y^\beta$}. 
$$

The notions of convergences established in \cite{CP} are defined as follows.
\begin{definition}
\label{Definition_Concept_P_Convergence}
\noindent Let  $0\leq \beta< 1$, a family of elements $\{y_\epsilon\}_{\eps \in (0,\eps_0]}$, $y_\epsilon \in Y^{\beta}_{\epsilon}$, is said to be $\mathcal{P^{\beta}}$-convergent to $y\in Y^{\beta}$ if $\|y_\epsilon-p^{\beta}_\epsilon y\|_{Y^\beta_\eps}\to 0$ as $\epsilon \to 0$. We denote this convergence by $y_\epsilon {\stackrel{\mathcal{P^{\beta}}}{\longrightarrow}} y$. If $\beta=0$ then we simply write this as $y_\epsilon {\stackrel{\mathcal{P}}{\longrightarrow}} y$.
\end{definition}

\begin{definition} \label{def_PP_conv_operators} A family of bounded linear
operators $T_\epsilon: Y^{\beta}_{\eps}\to Y_{\eps}$, $\eps\in (0,\eps_0]$, is said to be $\mathcal{P}$-convergent to the bounded linear
operator $T:Y^{\beta} \to Y$ as $\epsilon \to 0$, if $T_\epsilon
y_\epsilon \Pto Ty \in Y$ whenever $y_\epsilon \bPto y \in Y^{\beta}$. We denote
this by $T_\epsilon \PPto T$.
\end{definition}

\par\medskip Under all above assumptions, the authors showed the existence of attractors $\mathcal{A}\subset Y^\beta$ and $\mathcal{A}_\eps\subset Y_\eps^\beta$, for each $0< \eps\leq \eps_0$, see \cite[Theorem 1.6]{CP}. 
Moreover, in order to show that the attractors $\mathcal{A}_\eps$ ``converge'' in certain sense to $\mathcal{A}$, they explicit the following two assumptions (see \cite[Section 1]{CP}):

\bigskip

$
{\bf [A1]}
\left\{
\begin{array}{l}
\bullet \quad \mbox{$A$ is a closed linear operator with compact resolvent and}\\ \quad \displaystyle  \|(\lambda I-A)^{-1}\|_{\mathcal{L}(Y)}\leq \frac{M}{1+|\lambda|}\,\, \hbox{ for all } Re(\lambda)\geq 0;\\ \bullet \quad \exists   M_2>0,\; w_2\in \R: \, \displaystyle \|(\lambda I-A_\eps)^{-1}\|_{\mathcal{L}(Y_\eps)}\leq \frac{M_2}{|\lambda-w_2|}\hbox{ for all }Re(\lambda)>w_2,\, \eps\in(0,\eps_0] ;\\
\bullet \quad \mbox{The region $\Delta_{cc}$ of compact convergence of the resolvents is non-empty and}\\
\quad \mbox{for $\lambda\in\Delta_{cc}$, the resolvents $(\lambda I-A_\eps)^{-1}$ compactly converge to $(\lambda I-A)^{-1}$};\\
%\Delta_{cc}\ne  \emptyset \hbox{ and for $\lambda\in\Delta_{cc}$, the resolvents }   (\lambda I-A_\eps)^{-1}\in {\mathcal{L}(Y)}\hbox{ are compact};\\
\bullet \quad h_\eps(y_\eps)\Pto h(y)\hbox{ whenever } y_\eps\bPto y.
\end{array}
\right .
$

\bigskip

$
{\bf [A2]}
\left\{
\begin{array}{l}
\bullet \quad h'_\eps(y_\eps)\PPto h'(y)\hbox{ whenever } y_\eps\bPto y;\\
\bullet \quad \hbox{If }\, y_\eps^*\bPto y^* \mbox{ then } \displaystyle \sup_{\eps\in [0,\eps_0]}\sup_{\|z_\eps\|_{Y_\eps^\beta}\leq \rho}\| h'_\eps(z_\eps+y_\eps^*)\|_{\mathcal{L}(Y^{\beta}_\eps,Y^{\beta})}<\infty, \mbox{ for some $\rho>0$}.
\end{array}
\right .
$

\bigskip

In \cite{CP} they authors proved that if {\bf{[A1]}} holds, then the family $\{\mathcal{A}_\eps\}_{\eps \in [0,\eps_0]}$ is $\mathcal{P}^{\beta}$-upper semicontinuous in $Y^{\beta}_{\eps}$ at  $\eps=0$ (see \cite[Theorem 5.4]{CP}). Moreover, if  {\bf{[A1]}} and {\bf{[A2]}} hold and each equilibrium of the limiting problem is hyperbolic, then the family $\{\mathcal{A}_\epsilon\}_{\eps \in [0,\eps_0]}$ is $\mathcal{P}^{\beta}$-lower semicontinuous in $Y^{\beta}_{\eps}$ at $\epsilon=0$ (see \cite[Theorem 5.15]{CP}).

We want to apply the results from \cite{CP} to our problem. But to do this, we need to state clearly our choice of space $Y_\eps$,  our choice of $\beta$ and show that hypotheses {\bf[A1]} and {\bf[A2]} hold.  Notice that since we have a nonlinearity acting on the boundary of the domain, we will need to choose $Y_\eps$ as a function space with negative exponent say $X_\eps^{-\frac{\alpha}{2}}$ with the property that the functions in the space $X^{\frac{\alpha}{2}}_\eps$ have traces at the boundary. This imposes some restrictions on $\alpha$. Moreover, since we want to obtain convergence of the solutions, attractors, etc. in $H^1(\Omega_\eps)$ we would like to choose $\beta$ so that $Y^{\beta}_\eps=X^{\frac{1}{2}}_\eps=H^1(\Omega_\eps)$, that is, $\beta=\frac{1+\alpha}{2}$.   That is, the nonlinearity $h_\eps:Y^\beta_\eps\to Y_\eps$, or equivalently  $h_\eps:X^{\frac{1}{2}}_\eps\to X^{-\frac{\alpha}{2}}_{\eps}$.

Notice also that we need an extension operator $p_\eps:Y\to Y_\eps$, that is, $p_\eps:X^{-\frac{\alpha}{2}}\to X^{-\frac{\alpha}{2}}_{\epsilon}$ and it is not completely straightforward how to choose these operators. Moreover, once this operators are choosen, following \cite{CP} we construct $p_\eps^\beta: H^1(\Omega)\to H^1(\Omega_\eps)$ 
and we will obtain statements on $\mathcal{P}^\beta$-convergence. But we would like to relate this convergence to the $E$-convergence defined in the introduction where the operator $E$ is a ``standard'' extension and restriction operator in $H^1$.  We will actually show that these two convergences are equivalent.

\subsection{Fractional power space, uniform equivalence and traces}
\label{power space}

Before getting into the proof of our result we need to clarify the relation between the fractional power spaces, the Bessel potential spaces and its dependence with respect to the domain. Notice that we are perturbing the boundary of the domain in a non very smooth way and we want to understand how the norms of the spaces and the norm of the inclusions depend on $\eps$. 

So,  let us denote by  $[Z,W]_\theta$ the $\theta$-complex interpolation space of Banach spaces $Z$ and $W$, for $0\leq \theta \leq 1$, see \cite{Amann1,triebel} for more details on interpolation theory.

Keeping the notation from Subsection \ref{equation-setting}, we start proving that the fractional power spaces $X^\alpha_\epsilon$, the Bessel potential spaces $H^{2\alpha}(\Omega_\epsilon)$  and the $\alpha$-complex interpolation spaces  indeed are the same spaces with equivalent norms for $\alpha$ in certain range of values. Moreover,  and very important, the constants appearing in the inequalities of  the equivalence of the norms can be chosen uniformly with respect to $\epsilon\in [0,\eps_0]$, for some values of $\alpha$. 
Notice that by using this, we will be able to prove that the constants of the embeddings of $X^\alpha_\epsilon$ into the space $L^p(\Omega_\epsilon)$ and of the trace operators into $L^q(\partial\Omega_\epsilon)$ (for appropriate $p$ and $q$) can be chosen uniformly in $\epsilon$. 

\par\medskip Let us start with the following 
\begin{lemma} 
\label{equivalentdefinitionsspacesnorms}
Let $\{\Omega_\epsilon\}_{\epsilon \in [0,\epsilon_0]}$ be a family of domains satisfying conditions {\bf{(H)}} and {\bf{(F)(i)}}. Then for $0<\alpha<\half$, we have 
$$
X^\alpha_\epsilon = [X^0_\epsilon, X^1_\epsilon]_\alpha = [X^0_\epsilon, X^{\half}_\epsilon]_{2\alpha}= [L^2(\Omega_\epsilon), H^1(\Omega_\epsilon)]_{2\alpha} = H^{2\alpha}(\Omega_\epsilon),
$$
with equivalent norms uniformly with respect to $\epsilon$. 
\end{lemma}
\proof 
As we have mentioned in Subsection \ref{known-results}, we have the spectral convergence of $A_\epsilon$ to $A_0$ as $\epsilon \to 0$.  In this case,  since $A_\epsilon$ are self-adjoint, $A_\epsilon=A_\epsilon^*  \geq c>0$ uniformly in $\epsilon$. Then,   by \cite[Example 4.7.3 (a)]{Amann1}, the purely imaginary powers are bounded by 1, that means, 
$$
\|A^{it}_\epsilon\|_{{\cal{L}}(L^2(\Omega_\epsilon))} \leq 1, \quad  \mbox{for all $t\in \R$ and $0\leq \epsilon \leq \epsilon_0$}.
$$
Therefore, by \cite[Subsection 1.15.3]{triebel}, the fractional power space $X^\alpha_\epsilon$ is characterized by the complex interpolation space
$$
X^\alpha_\epsilon = [X^0_\epsilon, X^1_\epsilon]_\alpha = [X^0_\epsilon, X^{\half}_\epsilon]_{2\alpha},
$$
with equivalent norms independent of $\epsilon$. Since $X^{0}_\epsilon=L^2(\Omega_\epsilon)$ and $X^{\half}_\epsilon=H^1(\Omega_\epsilon)$, then 
$$
X^\alpha_\epsilon = [X^0_\epsilon, X^1_\epsilon]_\alpha = [X^0_\epsilon, X^{\half}_\epsilon]_{2\alpha}= [L^2(\Omega_\epsilon), H^1(\Omega_\epsilon)]_{2\alpha} = H^{2\alpha}(\Omega_\epsilon), \quad \mbox{for $0<\alpha<\half$}.
$$
Observe that, by \cite[Theorem 16.1]{Yagi},  $X^\alpha_\epsilon =D(A^{\alpha}_\epsilon) = [X^0_\epsilon, X^1_\epsilon]_\alpha $ with isometry. 
\cqd

\begin{remark}
In view of this lemma, we can always consider that the norm in $H^\alpha(\Omega_\eps)$, $0<\alpha<1$, is given by the norm of the interpolation between $L^2(\Omega_\eps)$ and $H^1(\Omega_\eps)$.
 \end{remark}

Following the same notation as \cite[Lemma 4.1]{AB0}, we obtain

\begin{lemma}
\label{Halphaepsilon}
Let $\{S_\epsilon\}_{\epsilon \in [0,\epsilon_0]}$ be a family of Lipschitz bounded domains in $\R^N$. Assume there exists a family of Lipschitz, one-to-one mappings $F_\epsilon$ from $S_0$ onto $S_\epsilon$ such that the inverse is Lipschitz, $\|DF_\epsilon\|_{\mathcal{L}(L^\infty(S_0)^{N\times N})} \leq K$ and $\|DF^{-1}_\epsilon\|_{\mathcal{L}(L^\infty(S_\epsilon)^{N\times N})} \leq K$, with $K>0$ independent of $\epsilon$. 
Then $u \in H^\alpha(S_\epsilon)$ if, and only if, $u\circ F_\epsilon \in H^\alpha(S_0)$, for $0\leq \alpha\leq 1$. Moreover, there exist constants $C, D>0$ independent of $\epsilon$ such that
$$
C \|u\circ F_\epsilon\|_{H^\alpha(S_0)}\leq  \|u\|_{H^\alpha(S_\epsilon)}\leq D \|u\circ F_\epsilon\|_{H^\alpha(S_0)}.
$$
\end{lemma}
\proof 
By considering the pull-back operator $F_\epsilon^*$ (that is $F_\eps^*(u)=u\circ F_\eps$), we know by \cite[Lemma 4.1]{AB0} that 
$H^1(S_\epsilon) {\stackrel{F_\epsilon^*}{\longrightarrow}} H^1(S_0)$ and $L^2(S_\epsilon) {\stackrel{F_\epsilon^*}{\longrightarrow}} L^2(S_0)$, satisfying the inequalities   
$$ C\|u \circ  F_\epsilon\|_{H^1(S_0)} \leq \|u\|_{H^1(S_\epsilon)}  \leq D\|u \circ F_\epsilon\|_{H^1(S_0)}
$$
and 
$$ C\|u \circ  F_\epsilon\|_{L^2(S_0)} \leq \|u\|_{L^2(S_\epsilon)}  \leq D\|u \circ F_\epsilon\|_{L^2(S_0)}, 
$$
with $C,D$ independent of $\eps$. 
Since $H^\alpha(S_\epsilon)=[L^2(S_\epsilon), H^1(S_\epsilon)]_{\alpha}$  and $H^\alpha(S_0)=[L^2(S_0), H^1(S_0)]_{\alpha}$, for $0<\alpha<1$, then by definition of complex interpolation and its norms, we obtain 
$u \in H^\alpha(S_\epsilon)$ if, and only if, $u\circ F_\epsilon \in H^\alpha(S_0)$ and there exist constants $C,D>0$ independent of $\epsilon$ such that
$$
C \|u\circ F_\epsilon\|_{H^\alpha(S_0)}\leq  \|u\|_{H^\alpha(S_\epsilon)}\leq D \|u\circ F_\epsilon\|_{H^\alpha(S_0)}.
$$
\cqd

\begin{lemma} 
\label{embeddings-traces-uniformepsilon}
Let $\{\Omega_\epsilon\}_{\epsilon \in [0,\epsilon_0]}$ be a family of domains satisfying conditions {\bf{(H)}} and {\bf{(F)(i)}}. We have  
$$
X^\alpha_\epsilon \hookrightarrow L^p(\Omega_\epsilon), \quad \mbox{for $0< \alpha\leq \frac{1}{2}$ and $1\leq p\leq \frac{2N}{N-4\alpha}$},
$$ 
$$
X^\alpha_\epsilon \hookrightarrow X^\beta_\epsilon, \quad \mbox{for $0<\beta<\alpha< \half$}, \qquad \qquad \quad
$$
with embeddings constants independent of $\epsilon$.  Moreover, the norm of the trace operator 
$$
T: X^\alpha_\epsilon \longrightarrow L^q(\partial\Omega_\epsilon),\quad \mbox{for $\frac{1}{4}< \alpha\leq \frac{1}{2}$ and $1\leq q\leq \frac{2(N-1)}{N-4\alpha}$},
$$
can be chosen uniformly for all $\eps\in [0,\eps_0]$.
\end{lemma}
\proof First, the case $\alpha =\frac{1}{2}$ follows from \cite[Propositions 4.2 and 4.3]{AB0}. Now, we observe that $ X^\alpha_\epsilon \hookrightarrow L^p(\Omega_\epsilon), 1\leq p\leq 2$ is trivially satisfied. Now, observe that $$
\|u_\epsilon\|_{L^p(\Omega_\epsilon)}\leq \sum^{m}_{i=1}\|u_\epsilon\|_{L^p(\Omega_\epsilon\cap U_i)}\quad \mbox{and}\quad \|u_\epsilon\|_{L^p(\partial \Omega_\epsilon)}\leq \sum^{m}_{i=1}\|u_\epsilon\|_{L^p(\partial \Omega_\epsilon\cap U_i)}.
$$
In the Lemma \ref{Halphaepsilon}, for each $i=1,\ldots, m$, consider $\Omega_\epsilon \cap U_i$ and $\Phi_{i,\epsilon}:Q_N^- \to  \Omega_\epsilon \cap U_i$, where $Q_N^-=Q_{N-1} \times (-1,0)\subset Q_N$.
 Then, using Lemma \ref{Halphaepsilon} and the standard embeddings for a fixed domain, for $0< \alpha< \frac{1}{2}$ and $1\leq p\leq \frac{2N}{N-4\alpha }$,  we obtain
$$\|u_\epsilon\|_{L^p(\Omega_\epsilon\cap U_i)} \leq D  \|u_\epsilon \circ  \Phi_{i,\epsilon}\|_{L^p(Q^-_{N})} \leq \tilde{D} \|u_\epsilon \circ  \Phi_{i,\epsilon}\|_{H^{2\alpha}(Q^-_{N})} \leq \hat{D}\|u_\epsilon\|_{H^{2\alpha}(\Omega_\epsilon \cap U_i)}  \leq \hat{D}\|u_\epsilon\|_{H^{2\alpha}(\Omega_\epsilon)}.
$$
Now, using \cite[Lemma 4.1]{AB0} and the continuity of the trace operator for a fixed domain, we have
$$\|u_\epsilon\|_{L^q(\partial \Omega_\epsilon\cap U_i)} \leq D  \|u_\epsilon \circ  \Phi_{i,\epsilon}\|_{L^q(Q_{N-1})} \leq \tilde{D} \|u_\epsilon \circ  \Phi_{i,\epsilon}\|_{H^{2\alpha}(Q_{N}^{-})} \leq \hat{D}\|u_\epsilon \|_{H^{2\alpha}(\Omega_\epsilon \cap U_i)}  \leq \hat{D}\|u_\epsilon\|_{H^{2\alpha}(\Omega_\epsilon)},
$$
for $\frac{1}{4}< \alpha< \frac{1}{2}$ and $1\leq q\leq \frac{2(N-1)}{N-4\alpha }$. And the results follows from Lemma \ref{equivalentdefinitionsspacesnorms}.
\cqd

%%%%%%%%%%%%%%%%%%%%%%%%%%%%%%%%%%%%%%%%%%%%%%%%%%%%%%%%%%%%%%%%%%%%%%%%%%%%%%%%%%

\subsection{Extension operators and  convergence in spaces with negative exponents}
\label{extensionoperators}

As we have mentioned in the introduction, since we are dealing with the family of domains $\{\Omega_\eps\}_{\epsilon \in [0,\epsilon_0]}$, we need to devise a way to compare functions defined in the different domains and give a meaning to the fact that a family of functions $u_\eps$ defined in $\Omega_\eps$ converges to a function $u$ defined in $\Omega$. An effective way to accomplish this is with the concept of convergence mediated via an extension operator, as considered in \cite{CP,St1972a, St1972b, St1972c,V} and applied to some problems in \cite{AB0,AB1,ACL3, CCS}. Here, we will extend the 
concepts of extension operators and convergence to
positive fractional power spaces and
we will define extension operators and  convergence in spaces with negative exponents.

To define properly these concepts, we need to consider operators transforming functions defined in $\Omega$ to functions defined in $\Omega_\eps$ and also the other way around, operators transforming functions defined in $\Omega_\eps$ to functions defined in $\Omega$. Following the same notation as in Section \ref{introd},   we consider first the linear operator $E_\eps:Z_0\to Z_\eps$ given by \eqref{definitionE_epsilon} and satisfying \eqref{extensaobase}, where $Z_\eps=H^1(\Omega_\eps)$ or $L^p(\Omega_\eps)$ or $W^{1,p}(\Omega_\eps)$, for $\eps \geq 0$ and $1\leq p\leq \infty$.

On the other hand, we also need to define a family of operators $\hat{E}_\epsilon: Z_\epsilon \to Z_0$. At first sight one is tempted to follow a similar procedure as in the case of $E_\epsilon$, that is, consider the extension operator from $\Omega_\eps$ to $\R^N$
and then the restriction to $\Omega$. But since we may want to consider oscillations at the boundary in $\Omega_\eps$, the norm of the extension operators from $\Omega_\eps$ to $\R^N$ may not be very well controlled.  Therefore, we proceed in a different direction, following the ideas of \cite{AB1,AC}.

We go to hypothesis {\bf(I)} and consider the  family $K_\epsilon$ of interior smooth perturbation of $\Omega$  with $K_\epsilon \subset \Omega \cap \Omega_{\epsilon}$ and $\theta_\epsilon : \Omega \to K_\epsilon$ a diffeomorphism satisfying this hypothesis {\bf{(I)}}.  Using $K_\epsilon$ and $\theta_\epsilon$, we define  $$\hat{E}_\epsilon u_\epsilon = {u_\epsilon}_{ |K_\epsilon} \circ \theta_\epsilon,$$
that is, if $x\in \Omega$ then $(\hat{E}_\epsilon u_\epsilon)(x)=u_\eps(\theta_\eps(x))$.
If  $\Omega \subset \Omega_\epsilon$, that is, $\Omega_\epsilon$ is an exterior perturbation of $\Omega$,  then we can consider $K_\epsilon =\Omega$, $\theta_\eps$ the identity and $\hat{E}_\epsilon u_\epsilon = {u_\epsilon}_{|\Omega}$.

Using  now that $D(A_\epsilon^{\frac{\alpha}{2}})= X^{\frac{\alpha}{2}}_\epsilon=[L^2(\Omega_\epsilon),H^1(\Omega_\epsilon)]_\alpha$, $0<\alpha<1$, we can  extend the definitions of $E_\epsilon$ and $\hat{E}_\epsilon$ to the scale of positive fractional power spaces via interpolation. Indeed,  we obtain  
$$
E_\epsilon: X^{\frac{\alpha}{2}} \to X^{\frac{\alpha}{2}}_\epsilon \quad \mbox{and}\quad \hat{E}_\epsilon: X^{\frac{\alpha}{2}}_\epsilon \to X^{\frac{\alpha}{2}}, \quad \mbox{for $0 \leq \alpha \leq 1$},
$$
where $X^{\frac{\alpha}{2}} \equiv X^{\frac{\alpha}{2}}_{0}$ and they satisfy 
\begin{equation}\label{Eproperties}
\|E_\epsilon u\|_{X^{\frac{\alpha}{2}}_\epsilon} \leq k \|u\|_{X^{\frac{\alpha}{2}}}, \quad \|E_\epsilon u\|_{X^{\frac{\alpha}{2}}_\epsilon} \to \|u\|_{X^{\frac{\alpha}{2}}},\quad  \mbox{as $\epsilon \to 0$}, 
\end{equation}
and 
\begin{equation}\label{Ehatproperties}  
\|\hat{E}_\epsilon u_\epsilon\|_{X^{\frac{\alpha}{2}}} \leq  \hat{k} \|u_\epsilon\|_{X^{\frac{\alpha}{2}}_\epsilon}, \ \  \|\hat{E}_\epsilon \|_{\mathcal{L} (X_{\epsilon}^{\frac{\alpha}{2}}, X^{\frac{\alpha}{2}})} \to 1, \quad \mbox{as $\epsilon \to 0$},
\end{equation}
 where $k, \hat{k} > 0$ are  independent of  $\epsilon$. 

We note that for $x \in K_{\epsilon_0}$ and $0<\epsilon\leq \epsilon_0$, we have  $(\hat{E}_{\epsilon}E_{\epsilon}u)(x)=u(x)$  and $(E_{\epsilon}\hat{E}_{\epsilon}u_\epsilon)(x)=u_\epsilon(x)$, for all $u\in X^{\frac{\alpha}{2}}$ and $u_\epsilon \in X_\epsilon^{\frac{\alpha}{2}}$.  

\par\medskip 
In the following, we prove a lemma that allows us to compare 
$E_{\epsilon}\hat{E}_{\epsilon}u_\epsilon$  and $u_\epsilon$  in $X_\epsilon^{\frac{\alpha}{2}}$. 

\begin{lemma} 
\label{auxiliar-similarLemma4.2AB}
Let $\{\Omega_\epsilon\}_{\epsilon \in [0,\epsilon_0]}$ be a family of domains satisfying conditions 
{\bf{(H)}}, {\bf{(F)(i)}} and {\bf{(I)}}. Then there exists a function $c(\epsilon)$ with $c(\epsilon) \to 0$ as $\epsilon \to 0$ such that, for any $u_\epsilon \in H^1(\Omega_\epsilon)$, we have:

\vspace{0.2cm}

\noindent (i)
$\displaystyle \int_{\partial K_\sigma} |u_\epsilon - u_\epsilon \circ \theta_\epsilon| \leq c(\epsilon) \|u_\epsilon\|_{H^1(\Omega_\epsilon)}, \quad \mbox{for $\sigma > \epsilon$};
$

\vspace{0.2cm}

\noindent (ii)
$\displaystyle \int_{K_\epsilon} |u_\epsilon - u_\epsilon \circ \theta_\epsilon|  \leq c(\epsilon) \|u_\epsilon\|_{H^1(\Omega_\epsilon)}.$
\end{lemma}
 \proof (i)  By considering the parameterization $\phi_i$ and $\theta_\sigma$,  we have 
$$\int_{\partial K_\sigma} |u_\epsilon - u_\epsilon \circ \theta_\epsilon|  \leq 
\sum_{i=1}^m \int_{Q_{N-1} }  |u_\epsilon \circ \theta_\sigma \circ \phi_i(x')  - u_\epsilon \circ \theta_\epsilon \circ \theta_\sigma \circ \phi_i(x')|  J_{N-1}  (\theta_\sigma \circ \phi_i)(x')dx'.
$$

 Considering $\theta_\epsilon \circ \theta_\sigma = \theta_{\sigma(\epsilon)} $ and using similar arguments as \cite[Lemma 4.2]{AB0}, we obtain 
 \begin{eqnarray*} \int_{\partial K_\sigma} |u_\epsilon - u_\epsilon \circ \theta_\epsilon|  & \leq &
\sum_{i=1}^m \int_{Q_{N-1} }  |u_\epsilon \circ \Phi_i(x^\prime, \hat{\rho}_{i,\sigma} (x^\prime)) -   u_\epsilon \circ \Phi_i(x^\prime, \hat{\rho}_{i,\sigma(\epsilon)} (x^\prime))|  J_{N-1}  (\theta_\sigma \circ \phi_i) (x^\prime) dx^\prime\\
&\leq& \sum_{i=1}^m \int_{Q_{N-1} }\left|  \int_{ \hat{\rho}_{i,\sigma} (x^\prime)}^{\hat{\rho}_{i,\sigma(\epsilon) } (x^\prime)}  \frac{\partial  (u_\epsilon \circ \Phi_i)}{\partial x_N}(x^\prime, y) dx_N \right|  J_{N-1}  (\theta_\sigma \circ \phi_i) (x^\prime) dx^\prime\\
&\leq& \sum_{i=1}^m \|\hat{\rho}_{i,\sigma(\epsilon)}-\hat{\rho}_{i,\sigma}\|^{\frac{1}{2}}_{L^{\infty}(Q_{N-1})}  \\
&. & \int_{Q_{N-1}}\left( \int_{\hat{\rho}_{i,\sigma} (x^\prime)}^{\hat{\rho}_{i,\sigma(\epsilon)} (x^\prime)} \left| \frac{\partial  (u_\epsilon \circ \Phi_i)}{\partial x_N}(x^\prime, y)\right|^2 dx_N\right)^{\frac{1}{2}}   J_{N-1}  (\theta_\sigma \circ \phi_i) (x^\prime) dx^\prime
\\
&\leq&
c(\epsilon) \|u_\epsilon\|_{H^1(\Omega_\epsilon)},
\end{eqnarray*} 
with $c(\epsilon) \to 0$ as $\epsilon \to 0$ since $\hat{\rho}_{i,\sigma(\epsilon)}\to \hat{\rho}_{i,\sigma}$.

\vspace{0.2cm}

\noindent (ii)  Since ${\theta_\epsilon}_{|K_{\epsilon_0}} = id$, then  $\displaystyle \int_{K_{\epsilon_0}} |u_\epsilon - u_\epsilon \circ \theta_\epsilon| = 0$. Thus,
$$\int_{K_\epsilon} |u_\epsilon - u_\epsilon \circ \theta_\epsilon| = \int_{K_\epsilon\setminus K_{\epsilon_0}} |u_\epsilon - u_\epsilon \circ \theta_\epsilon|  = \int_{\epsilon}^{\epsilon_0} \int_{ \partial K_\sigma} |u_\epsilon - u_\epsilon \circ \theta_\epsilon| d\sigma  \leq c(\epsilon) \|u_\epsilon\|_{H^1(\Omega_\epsilon)}.$$
\cqd

\begin{lemma} 
\label{convergence_extension_ofthe_restriction}
Let $\{\Omega_\epsilon\}_{\epsilon \in [0,\epsilon_0]}$ be a family of domains satisfying conditions {\bf{(H)}}, {\bf{(F)(i)}} and {\bf{(I)}}. If $0<\alpha<1$ and $u_\epsilon \in X^{\frac{1}{2}}_{\epsilon}$ such that $\|u_\epsilon\|_{X^{\frac{1}{2}}_{\epsilon}} \leq K$, for some $K>0$ independent of $\epsilon$, then
$$
\|E_{\epsilon}\hat{E}_{\epsilon}u_{\epsilon} -u_{\epsilon}\|_{X^{\frac{\alpha}{2}}_{\epsilon}} \to 0, \quad \mbox{as  $\epsilon \to 0$}. 
$$
\end{lemma}
\proof
By \cite[Lemma 3.4]{AB1}, the interpolation's constants are uniformly bounded in $\epsilon$, and using Lemma \ref{equivalentdefinitionsspacesnorms}, we have
\begin{eqnarray*}
\|E_{\epsilon}\hat{E}_{\epsilon}u_{\epsilon} -u_{\epsilon}\|_{X^{\frac{\alpha}{2}}_{\epsilon}} \!\!\!&\leq&\!\!\! C\|E_{\epsilon}\hat{E}_{\epsilon}u_{\epsilon} -u_{\epsilon}\|_{H^{\alpha}(\Omega_{\epsilon})}\leq  \tilde{K}\|E_{\epsilon}\hat{E}_{\epsilon}u_{\epsilon} -u_{\epsilon}\|^{1-\alpha}_{L^2(\Omega_{\epsilon})}\|E_{\epsilon}\hat{E}_{\epsilon}u_{\epsilon} -u_{\epsilon}\|^{\alpha}_{H^1(\Omega_{\epsilon})}.
%\\ & \leq &
%2K\tilde{K} \left(\int_{\Omega_\epsilon}\!\! |E_{\epsilon}\hat{E}_{\epsilon}u_{\epsilon} -u_{\epsilon}|^2\right)^{\frac{1-\alpha}{2}} %\leq 2K\tilde{K} \left(\int_{\Omega_\epsilon\setminus \Omega}\!\!\!|E_{\epsilon}\hat{E}_{\epsilon}u_{\epsilon} -u_{\epsilon}|^2\right)^{\frac{1-\alpha}{2}},
\end{eqnarray*}
%where $K>0$ is independent of $\epsilon$. 
Since $E_{\epsilon}\hat{E}_{\epsilon}u_{\epsilon} = u_\epsilon\circ \theta_\epsilon$ in $K_\epsilon \subset \Omega$, then 
$$
\|E_{\epsilon}\hat{E}_{\epsilon}u_{\epsilon} -u_{\epsilon}\|^2_{L^2(\Omega_{\epsilon})}  = \int_{\Omega_\epsilon} |E_{\epsilon}\hat{E}_{\epsilon}u_{\epsilon} - u_\epsilon|^2  = \int_{\Omega_\epsilon\setminus K_\epsilon} |E_{\epsilon}\hat{E}_{\epsilon}u_{\epsilon} - u_\epsilon|^2  + \int_{K_\epsilon} |u_\epsilon \circ \theta_\epsilon - u_\epsilon|^2.  
$$ 
Thus, using Holder Inequality,  $|\Omega_\epsilon \setminus K_\epsilon| \to 0$ as $\epsilon \to 0$ and Lemma \ref{auxiliar-similarLemma4.2AB}, the result follows. 
\cqd

If  $\Omega \subset \Omega_\epsilon$, that is, $\Omega_\epsilon$ is an exterior perturbation of $\Omega$,  then $\hat{E}_{\epsilon}{E}_{\epsilon}u = u$. 
In the general case,  we obtain  
\begin{lemma} 
\label{convergence_restriction_ofthe_extension}
Let $\{\Omega_\epsilon\}_{\epsilon \in [0,\epsilon_0]}$ be a family of domains satisfying conditions {\bf{(H)}}, {\bf{(F)}(i)} and {\bf{(I)}}. If $0<\alpha<1$ and $u \in X^{\frac{1}{2}}$ such that $\|u\|_{X^{\frac{1}{2}}} \leq K$, for some $K>0$, then
$$
\|\hat{E}_{\epsilon}E_{\epsilon}u -u\|_{X^{\frac{\alpha}{2}}} \to 0, \quad \mbox{as  $\epsilon \to 0$}.
$$
\end{lemma}

\proof Using Lemma \ref{equivalentdefinitionsspacesnorms} and interpolation properties (see \cite[Theorem 1.11.3]{triebel} and  \cite[pag. 15]{Yagi}), we have
\begin{eqnarray*}
\|\hat{E}_{\epsilon}{E}_{\epsilon}u  -u \|_{X^{\frac{\alpha}{2}}} \!\!\!&\leq&\!\!\! C\|\hat{E}_{\epsilon}E_{\epsilon}u_{\epsilon} -u_{\epsilon}\|_{H^{\alpha}(\Omega)}\leq C\|\hat{E}_{\epsilon}{E}_{\epsilon}u  -u \|^{1-\alpha}_{L^2(\Omega)}\|\hat{E}_{\epsilon}{E}_{\epsilon}u  -u \|^{\alpha}_{H^1(\Omega)}. 
%\\ & \leq &
%2K\hat{K} \left(\int_{\Omega_\epsilon}\!\! |E_{\epsilon}\hat{E}_{\epsilon}u_{\epsilon} -u_{\epsilon}|^2\right)^{\frac{1-\alpha}{2}} \leq 2K\hat{K} \left(\int_{\Omega_\epsilon\setminus \Omega}\!\!\!|E_{\epsilon}\hat{E}_{\epsilon}u_{\epsilon} -u_{\epsilon}|^2\right)^{\frac{1-\alpha}{2}},
\end{eqnarray*}
Since $\hat{E}_\epsilon E_\epsilon u = u\circ\theta_\epsilon$ in $\Omega$ and ${\theta_\epsilon}_{|K_{\epsilon_0}} = id$, then 
$$
\int_{\Omega }\!|\hat{E}_{\epsilon}E_{\epsilon}u  -u |^2  =  \int_{\Omega} \!\!|u\circ\theta_\epsilon -u|^2 = \int_{\Omega\setminus K_{\epsilon_0} } \!\!|u\circ\theta_\epsilon -u|^2 \leq c(\epsilon)\|u\|^{2}_{H^{1}(\Omega)},
$$
with $c(\epsilon)\to 0$ as $\epsilon \to 0$, where in the last inequality we use the same arguments of Lemma \ref{auxiliar-similarLemma4.2AB}.

Therefore,
$$
\|\hat{E}_{\epsilon}{E}_{\epsilon}u  -u \|_{X^{\frac{\alpha}{2}}} \leq c(\epsilon, K)  \to 0, \quad \mbox{as  $\epsilon \to 0$}.
$$
This concludes the proof of the lemma. \cqd

The concepts of $E$-convergence and  $E$-weak convergence in positive fractional power spaces are defined as follows. 

\begin{definition}
\label{Definition_Concept_E_Convergence}
Let $0\leq \alpha\leq 1$, $u_\epsilon \in X^{\frac{\alpha}{2}}_\eps$ and $u \in X^{\frac{\alpha}{2}}$. We have:

\vspace{0.2cm}

\noindent (i) The family $\{u_\epsilon\}_{\eps \in (0,\eps_0]}$ $E$-converges to $u$ if $\|u_\epsilon-E_\epsilon u\|_{X^\frac{\alpha}{2}_\eps}\to 0$ as $\epsilon \to 0$. We denote this convergence by $u_\epsilon {\stackrel{E}{\longrightarrow}} u$;

\vspace{0.2cm}

\noindent (ii) The family $\{u_\epsilon\}_{\eps \in (0,\eps_0]}$ $E$-weak converges to $u$ if $(w_\epsilon,u_\epsilon)_{X^\frac{1}{2}_\eps}\to (w,u)_{X^\frac{1}{2}}$ as $\epsilon \to 0$, for any sequence $w_\epsilon {\stackrel{E}{\longrightarrow}} w$ in $X^\frac{1}{2}_\eps$, where $(\cdot,\cdot)_{X^\frac{1}{2}_\eps}$ and $(\cdot,\cdot)_{X^\frac{1}{2}}$  denote the inner product in $X^\frac{1}{2}_\eps$ and $X^\frac{1}{2}$, respectively. We denote this convergence by $u_\epsilon {\stackrel{E}{-{\hspace{-2mm}}\rightharpoonup}} u$.
\end{definition}

Once we have stablished the concept of $E$-convergence in $X^\frac{\alpha}{2}_\eps$ for $0\leq \alpha\leq 1$, we are going to extend appropriately this concept to spaces $X^{-\frac{\alpha}{2}}_\eps$ with $0< \alpha\leq 1$.
We define the following ``extension'' and ``restriction'' operators in spaces
with negative exponents by duality.

\begin{definition} \label{definition_extension_restrictions_dualspaces}
For $0 < \alpha \leq 1$, we define
\begin{equation}
\label{extensaodual}
\begin{array}{cccl}
E^*_{\epsilon}& : X^{-\frac{\alpha}{2}} &  \to & X^{-\frac{\alpha}{2}}_\epsilon
\\   
& \varphi & \mapsto & E^*_{\epsilon}\varphi: X^{\frac{\alpha}{2}}_\epsilon \; \to \; \mathbb{R} 
\end{array}
\end{equation}
where $\langle E^*_{\epsilon}\varphi,u_\epsilon \rangle = \langle \varphi,\hat{E}_\epsilon u_\epsilon\rangle$  for all $u_\epsilon \in X^{\frac{\alpha}{2}}_\epsilon$. On the other hand, we define 
\begin{equation}
\label{restricaodual}
\begin{array}{cccl} 
\hat{E}^*_{\epsilon}& : X^{-\frac{\alpha}{2}}_\epsilon &  \to & X^{-\frac{\alpha}{2}}\\   
& \varphi_\epsilon & \mapsto & \hat{E}^*_{\epsilon} \varphi_\epsilon : X^{\frac{\alpha}{2}} \; \to \; \mathbb{R} 
\end{array}
\end{equation}
where $\langle \hat{E}^*_{\epsilon}\varphi_\epsilon,u\rangle = \langle \varphi_\epsilon,E_\epsilon u\rangle$ for all $u \in X^{\frac{\alpha}{2}}$.
\end{definition}

We note that these operators satisfy 
$$
\|E^{*}_{\epsilon} \varphi\|_{X^{-\frac{\alpha}{2}}_\epsilon} \leq k^{*} \|\varphi\|_{X^{-\frac{\alpha}{2}}} \quad \mbox{and} \quad \|\hat{E}^{*}_{\epsilon} \varphi_\epsilon\|_{X^{-\frac{\alpha}{2}}} \leq  \hat{k}^{*} \|\varphi_\epsilon\|_{X^{-\frac{\alpha}{2}}_\epsilon},
$$ 
where $k^{*}, \hat{k}^{*} > 0$ are  independent of  $\epsilon$.

Again, if  $\Omega \subset \Omega_\epsilon$  then $\hat{E}^*_{\epsilon}E^*_{\epsilon}\varphi=\varphi$ for all $\varphi\in X^{-\frac{\alpha}{2}}$.  
In the general case,

\begin{lemma} 
\label{*convergence_restriction_ofthe_extension*}
Let $\{\Omega_\epsilon\}_{\epsilon \in [0,\epsilon_0]}$ be a family of domains satisfying conditions {\bf{(H)}}, {\bf{(F)}(i)} and {\bf{(I)}}. If $0<\alpha<1$ then
$$
\langle \hat{E}^*_{\epsilon}E^*_{\epsilon}\varphi,u\rangle  \to \langle \varphi,u\rangle, \quad \mbox{as  $\epsilon \to 0$},
$$
uniformly for $u$ in bounded sets of $X^\half$ and $\varphi$ in bounded sets of $X^{-\frac{\alpha}{2}}$.
\end{lemma}
\proof In fact, from Definition \ref{definition_extension_restrictions_dualspaces}, we have
\begin{eqnarray*}|\langle \hat{E}^*_{\epsilon}E^*_{\epsilon}\varphi,u\rangle  - \langle \varphi,u\rangle|  = |\langle E^*_{\epsilon}\varphi,E_{\epsilon}u\rangle  - \langle \varphi,u\rangle | = |\langle \varphi,\hat{E}_\epsilon E_{\epsilon}u\rangle  - \langle \varphi,u\rangle|\leq \|\varphi\|_{X^{-\frac{\alpha}{2}}} \|\hat{E}_\epsilon E_\epsilon u - u\|_{X^{\frac{\alpha}{2}}},
\end{eqnarray*} 
and the convergence follows using Lemma \ref{convergence_restriction_ofthe_extension}. 
\cqd

The following result proves that $E^*$ satisfies the condition as defined in \eqref{conv_p_eps} and that we will also see in Section \ref{EvsP}.

\begin{lemma} 
\label{conditionE*convergence} 
Let $\{\Omega_\epsilon\}_{\epsilon \in [0,\epsilon_0]}$  be a family of domains satisfying conditions {\bf{(H)}}, {\bf{(F)(i)}} and {\bf{(I)}}, then
$$\|E^*_\epsilon \varphi \|_{X^{-\frac{\alpha}{2}}_\epsilon} \to \|\varphi\|_{X^{-\frac{\alpha}{2}}}, \quad \mbox{as $\epsilon \to 0$ for any $\varphi \in X^{-\frac{\alpha}{2}}$}.
$$
\end{lemma}
\proof
In fact, for each $\varphi \in X^{-\frac{\alpha}{2}}$, we have % if $\Omega \subset \Omega_\epsilon$}, using the definition of $\|E_\epsilon^* \varphi\|_{X^{-\frac{\alpha}{2}}_\epsilon}$ we obtain  $\|E^*_\epsilon \varphi \|_{X^{-\frac{\alpha}{2}}_\epsilon} \leq \| \varphi\|_{X^{-\frac{\alpha}{2}}}$, for all $0<\epsilon \leq \epsilon_0$.
%{\large \color{blue} If $\Omega_\epsilon$ is not exterior to $\Omega$, .... MISSING!!!!!!
 \begin{eqnarray*} \|E_\epsilon^* \varphi\|_{X^{-\frac{\alpha}{2}}_\epsilon} &=& {\sup_{\begin{array}{ll}u_\epsilon \in X_\epsilon^\frac{\alpha}{2}\\ \|u_\epsilon\|_{X_\epsilon^\frac{\alpha}{2}}=1\end{array}}}
|\langle E_\epsilon^* \varphi, u_\epsilon\rangle|= {\sup_{\begin{array}{ll}u_\epsilon \in X_\epsilon^\frac{\alpha}{2}\\ \|u_\epsilon\|_{X_\epsilon^\frac{\alpha}{2}}=1\end{array}}}
|\langle \varphi, \hat{E}_\epsilon u_\epsilon\rangle|   \\
& \leq &   {\sup_{\begin{array}{ll}u_\epsilon \in X_\epsilon^\frac{\alpha}{2}\\ \|u_\epsilon\|_{X_\epsilon^\frac{\alpha}{2}}=1\end{array}}}  \|\varphi\|_{X^{-\frac{\alpha}{2}}} \| \hat{E}_\epsilon u_\epsilon\|_{X^{\frac{\alpha}{2}}}
\leq  \|\varphi\|_{X^{-\frac{\alpha}{2}}} \|\hat{E}_\epsilon \|_{{\cal{L}}(X_\epsilon^{\frac{\alpha}{2}},X^{\frac{\alpha}{2}})}
%\\
%& \leq &  \|\varphi\|_{X^{-\frac{\alpha}{2}}} (1+\eta_\epsilon)
\to 
\|\varphi\|_{X^{-\frac{\alpha}{2}}}, \quad \mbox{as $\epsilon \to 0$,}
 \end{eqnarray*} 
where we are using that $\|\hat{E}_\epsilon \|_{{\cal{L}}(X_\epsilon^{\frac{\alpha}{2}},X^{\frac{\alpha}{2}})}\to 1$, see \eqref{Ehatproperties}.

On the other hand, we note
$|\langle E^*_\epsilon\varphi,u_\epsilon\rangle|\leq\|E^*_\epsilon\varphi\|_{X^{-\frac{\alpha}{2}}_\epsilon}\|u_\epsilon\|_{X^{\frac{\alpha}{2}}_\epsilon}.$ %Moreover, $\langle E^*_\epsilon\varphi,u_\epsilon\rangle=\langle\varphi,\hat{E}_{\epsilon}u_\epsilon\rangle.$ 
Now, using the definition of $\| \varphi\|_{X^{-\frac{\alpha}{2}}}$ we obtain for each $\eta >0$  there exist $u \in X^\frac{\alpha}{2}$ and $\epsilon=\epsilon(\eta)$   such that 
$$\|\varphi\|_{X^{-\frac{\alpha}{2}}}-\eta<\frac{\langle \varphi,\hat{E}_\epsilon u_\epsilon\rangle}{\|u_\epsilon\|_{X^{\frac{\alpha}{2}}_\epsilon}} \leq\|\varphi\|_{X^{-\frac{\alpha}{2}}}.$$
Thus,
$$
\| \varphi\|_{X^{-\frac{\alpha}{2}}} -\eta < \frac{|\langle E^*_\epsilon \varphi,u_{\epsilon}\rangle|}{\|u_{\epsilon}\|_{X^{\frac{\alpha}{2}}_\epsilon}} \leq \|E_\epsilon^* \varphi\|_{X^{-\frac{\alpha}{2}}_\epsilon}, \quad \mbox{for $u_\epsilon = E_\epsilon u,\; \epsilon \leq \epsilon(\eta)$}.
$$ 
Therefore, 
$$
\|E^*_\epsilon \varphi  \|_{X^{-\frac{\alpha}{2}}_\epsilon} \to \|\varphi\|_{X^{-\frac{\alpha}{2}}}, \quad \mbox{as $\epsilon \to 0$.}
$$
\cqd

We also note that $\hat{E}^*_\epsilon \neq \hat{E}_\epsilon$ if $\alpha =0$, even in the case where the $\Omega_\epsilon$ is an exterior perturbation of $\Omega$. In fact, since $\langle \hat{E}^*_\epsilon u_\epsilon,v\rangle = \int_{\Omega_\epsilon}u_\epsilon E_\epsilon v$ and $\hat{E}_\epsilon u_\epsilon = {u_\epsilon}_{|\Omega}$ that means, $\hat{E}_\epsilon u_\epsilon$ as a functional in $L^2(\Omega)$ is given by $\langle \hat{E}_\epsilon u_\epsilon,v\rangle=\int_\Omega u_\epsilon v$.

\par\medskip

Now, we define the concepts of $E^*$-convergence and  $E^*$-weak convergence in negative fractional power spaces, 
using the operator $E^*_{\epsilon}$.

\begin{definition} 
\label{definition_E*epsilon}
Let $0< \alpha\leq 1$, $\varphi_\epsilon \in X^{-\frac{\alpha}{2}}_\epsilon$ and $\varphi \in X^{-\frac{\alpha}{2}}$. We have:

\vspace{0.2cm}

\noindent (i) The family $\{\varphi_\epsilon\}_{\eps \in (0,\eps_0]}$ $E^*$-converges to $\varphi$  if  
$\|\varphi_\epsilon - E^*_\epsilon \varphi\|_{X^{-\frac{\alpha}{2}}_\epsilon} \to 0$ as $\epsilon \to 0$. We denote this convergence by $\varphi_{\epsilon} {\stackrel{E^*}{\longrightarrow}} \varphi$;

\vspace{0.2cm}

\noindent (ii)  The family $\{\varphi_\epsilon\}_{\eps \in (0,\eps_0]}$ $E^*$-weak  converges to $\varphi$ if $\langle \varphi_\epsilon,u_\epsilon\rangle \to \langle \varphi,u\rangle$ whenever $u_\epsilon \dto u$ in  $X^{\frac{\alpha}{2}}_{\epsilon}$. We denote this convergence by $\varphi_\epsilon {\stackrel{E^*}{-{\hspace{-2mm}}\rightharpoonup}} \varphi$. 
\end{definition}

\begin{lemma} 
\label{lem_E_E} \label{Econvergence->E*convergence}
Let $\{\Omega_\epsilon\}_{\epsilon \in [0,\epsilon_0]}$ be a family of domains satisfying conditions {\bf{(H)}}, {\bf{(F)(i)}} and {\bf{(I)}}. If $v_\epsilon \in X^{\frac{1}{2}}_{\epsilon}$ and $v_0 \in X^{\frac{1}{2}}$ with  $v_{\epsilon} {\stackrel{E}{\longrightarrow}} v_0$ in $X^{\frac{1}{2}}_{\epsilon}$, then $
v_{\epsilon} {\stackrel{E^*}{\longrightarrow}} v_0$ in $X^{-\frac{\alpha}{2}}_{\epsilon}$, for $0<\alpha\leq 1$.
\end{lemma}
\proof Initially, for $\alpha=1$, we have
$$
\|v_\epsilon - E^*_{\epsilon} v_0\|_{X^{-\frac{1}{2}}_{\epsilon}}  =  {\ds{\sup_{\begin{array}{ll}u_\epsilon \in X^{\half}_{\epsilon}\\ \|u_\epsilon\|_{ X^{\half}_{\epsilon}}=1\end{array}}}}| \langle v_\epsilon -E^*_\epsilon v_0, u_\epsilon \rangle|.
$$
Since
\begin{equation}
\label{EQN1}
|\langle v_\epsilon -E^*_\epsilon v_0, u_\epsilon \rangle| = |\langle v_\epsilon , u_\epsilon \rangle -\langle  v_0 ,\hat{E}_\epsilon u_\epsilon\rangle| \leq |\langle v_\epsilon, u_\epsilon \rangle - \langle E_\epsilon v_0, u_\epsilon \rangle|+ | \langle E_\epsilon v_0, u_\epsilon \rangle - \langle v_0, \hat{E}_\epsilon u_\epsilon \rangle|.
\end{equation}
We have,
\begin{equation}
\label{EQN2}
\begin{array}{lll}
|\langle v_\epsilon, u_\epsilon \rangle - \langle E_\epsilon v_0, u_\epsilon \rangle| % & = & \left|\displaystyle \int_{\Omega_\epsilon} \nabla (v_\epsilon-E_\epsilon v_0)\nabla u_\epsilon+ \int_{\Omega_\epsilon} (v_\epsilon-E_\epsilon v_0)u_\epsilon\right| \\
%&\leq& (v_\epsilon-E_\epsilon v_0, u_\epsilon)_{H^1(\Omega_\epsilon)} 
\leq 
\|v_\epsilon-E_\epsilon v_0\|_{H^{1}(\Omega_\epsilon)} \|u_\epsilon\|_{H^{1}(\Omega_\epsilon)} \to 0, \quad \mbox{as $\epsilon \to 0$},
\end{array}
\end{equation}
uniformly in $u_\epsilon$ such that $\|u_\epsilon\|_{H^{1}(\Omega_\epsilon)} =1$. Moreover,
\begin{equation}
\label{EQN3}
\begin{array}{lll}
\quad |\langle E_\epsilon v_0, u_\epsilon \rangle-  \langle v_0, \hat{E}_\epsilon u_\epsilon \rangle|  = \left|\displaystyle \int_{\Omega_\epsilon}\!\! \nabla E_\epsilon v_0\nabla u_\epsilon +\int_{\Omega_\epsilon}\!\!  E_\epsilon v_0 u_\epsilon - \int_{\Omega}\!\! \nabla v_0\nabla \hat{E}_\epsilon u_\epsilon -\int_{\Omega}\!\! v_0 \hat{E}_\epsilon u_\epsilon \right|
\\
=
\left| \displaystyle \int_{\Omega_\epsilon \setminus K_\epsilon}\!\!\!\! (\nabla E_\epsilon v_0\nabla u_\epsilon + E_\epsilon v_0 u_\epsilon) + \int_{ K_\epsilon}\!\!\!\! (\nabla E_\epsilon v_0\nabla u_\epsilon  +  E_\epsilon v_0 u_\epsilon) - \displaystyle \int_{\Omega}\!\! (\nabla  v_0\nabla \hat{E}_\epsilon u_\epsilon + v_0 \hat{E}_\epsilon u_\epsilon) \right|. 
\end{array}
\end{equation}
Since 
\begin{equation}\label{parte0}  \left|  \displaystyle \int_{\Omega_\epsilon \setminus K_\epsilon}\!\!\!\! (\nabla E_\epsilon v_0\nabla u_\epsilon + E_\epsilon v_0 u_\epsilon) \right| \leq  \|u_\epsilon\|_{H^1(\Omega_\epsilon)} \|E_\epsilon v_0\|_{H^1(\Omega_\epsilon \setminus K_\epsilon)}\to 0, \quad \mbox{as $\epsilon \to 0$}, \quad 
\end{equation}
uniformly in $u_\epsilon$ such that $\|u_\epsilon\|_{H^{1}(\Omega_\epsilon)} =1$, and 
\begin{eqnarray*}
\begin{array}{lll}
& &\displaystyle \int_{ K_\epsilon}\!\!\!\! (\nabla E_\epsilon v_0\nabla u_\epsilon  +  E_\epsilon v_0 u_\epsilon) - \displaystyle \int_{\Omega}\!\! (\nabla  v_0\nabla \hat{E}_\epsilon u_\epsilon + v_0 \hat{E}_\epsilon u_\epsilon) \\
& = & \displaystyle \int_{ K_\epsilon}\!\!\!\! (\nabla E_\epsilon v_0\nabla u_\epsilon  +  E_\epsilon v_0 u_\epsilon)(1- J_N\theta_\epsilon^{-1} )  
+\int_{ K_\epsilon}\!\!\!\! (\nabla E_\epsilon v_0\nabla u_\epsilon  +  E_\epsilon v_0 u_\epsilon) J_N\theta_\epsilon^{-1} \\
& - & \displaystyle \int_{\Omega}\!\! \nabla  v_0\nabla \hat{E}_\epsilon u_\epsilon - \int_{\Omega}\!\!  v_0 \hat{E}_\epsilon u_\epsilon\\
& = & \displaystyle (I) +(II) - (III) - (IV). 
\end{array}
\end{eqnarray*}
We have,
$$
(III) = \int_{K_\epsilon}\!\! \nabla  v_0(\theta_\epsilon^{-1}(x)) \nabla  u_\epsilon(x) D\theta_\epsilon(\theta_\epsilon^{-1}(x)) J_N\theta_\epsilon^{-1}(x)dx
$$
and
$$
(IV) = \int_{K_\epsilon}\!\!  v_0(\theta_\epsilon^{-1} (x)) u_\epsilon(x) J_N\theta_\epsilon^{-1}(x) dx.
$$
Since $ J_N\theta_\epsilon^{-1} \to 1$  in $L^\infty(K_\epsilon)$ as $\epsilon \to 0$, then 
\begin{equation} \label{parteI} 
(I) = \int_{ K_\epsilon}\!\!\!\! (\nabla E_\epsilon v_0\nabla u_\epsilon  +  E_\epsilon v_0 u_\epsilon)(1- J_N\theta_\epsilon^{-1} )   \to 0, \quad \mbox{as $\epsilon \to 0$}.
\end{equation} 
Also, 
\begin{equation} \label{partesII-IV} 
\begin{array}{rcl} 
(II) -(III)-(IV) & = & \displaystyle \int_{K_\epsilon}\!\! [\nabla  v_0(x) \nabla  u_\epsilon(x)  - \nabla  v_0(\theta_\epsilon^{-1}(x))\nabla  u_\epsilon(x) D\theta_\epsilon(\theta_\epsilon^{-1}(x))] J_N\theta_\epsilon^{-1}(x)dx \\  
& + & \displaystyle \int_{K_\epsilon}\!\!  [v_0(x)- v_0(\theta_\epsilon^{-1} (x))] u_\epsilon(x) J_N\theta_\epsilon^{-1}(x) dx\to 0, \quad \mbox{as $\epsilon \to 0$}.
\end{array} 
\end{equation}

From \eqref{parte0}, \eqref{parteI} and \eqref{partesII-IV}, we have $\|v_\epsilon - E^*_{\epsilon} v_0\|_{X^{-\frac{1}{2}}_{\epsilon}} \to 0$, as $\epsilon \to 0$. 

Finally, using Lemma \ref{equivalentdefinitionsspacesnorms} and interpolation properties, we have
$$\|v_\epsilon -E^*_\epsilon v_0 \|_{X^{-\frac{\alpha}{2}}_{\epsilon}} \leq C \|v_\epsilon -E^*_\epsilon v_0 \|_{H^{-\alpha}(\Omega_\epsilon)} \leq C  \|v_\epsilon -E^*_\epsilon v_0 \|^{1-\alpha}_{H^{-1}(\Omega_{\epsilon})} \|v_\epsilon -E^*_\epsilon v_0 \|^{\alpha}_{L^2(\Omega_\epsilon)},$$
and since $\|v_\epsilon -E^*_\epsilon v_0 \|_{L^2(\Omega_\epsilon)}$ is uniformly bounded, we obtain the convergence in $X_\epsilon^{-\frac{\alpha}{2}}$. 
\cqd

\begin{lemma} 
\label{lemma4.3AB0} 
Let $\{\Omega_\epsilon\}_{\epsilon \in [0,\epsilon_0]}$ be a family of domains satisfying conditions {\bf{(H)}}, {\bf{(F)(i)}} and {\bf{(I)}}. Let $u_\epsilon \in X^{\frac{1}{2}}_{\epsilon}$ such that $\| u_\epsilon\|_ { X^{\frac{1}{2}}_{\epsilon}} \leq K$, for some $K>0$ independent of $\epsilon$. Then, there exist a subsequence denoted by $u_{\epsilon_k}$ and $u_0\in X^{\frac{1}{2}}$ such that $u_{\epsilon_k} \dwto u_0$ and $\hat{E}_{\epsilon_k} u_{\epsilon_k} \wto  u_0$  in $X^\half$.
\end{lemma}
\proof 
If $\| u_\epsilon\|_ { X^{\frac{1}{2}}_{\epsilon}} \leq K$ then $\| \hat{E}_\epsilon u_\epsilon \|_ { X^{\frac{1}{2}}} \leq \tilde{K}$. Thus, there exist a subsequence $\hat{E}_{\epsilon_k} u_{\epsilon_k}$ and $u_0\in  X^{\frac{1}{2}}$ such that $\hat{E}_{\epsilon_k} u_{\epsilon_k} \wto u_0$ in $ X^{\frac{1}{2}}$. 
%and $\hat{E}_{\epsilon_k} u_{\epsilon_k} \to u$ in $X^{\frac{\alpha}{2}}$, for $0<\alpha<1$.

Using \cite[Proposition 3.1]{AB0}, to prove that $u_{\epsilon_k} \dwto u_0$ it is enough to prove that 
$$
(E_{\epsilon_k}v,u_{\epsilon_k})_{ X^{\frac{1}{2}}_{\epsilon_k}}\to (v,u_0)_{ X^{\frac{1}{2}}}, \quad \mbox{for all $v\in  X^{\frac{1}{2}}$}.
$$

For each $v\in  X^{\frac{1}{2}}$, we have
$$
|(E_{\epsilon_k}v,u_{\epsilon_k})_{ X^{\frac{1}{2}}_{\epsilon_k}}- (v,u_0)_{ X^{\frac{1}{2}}}|\leq |(E_{\epsilon_k}v,u_{\epsilon_k})_{ X^{\frac{1}{2}}_{\epsilon_k}} - (v,\hat{E}_{\epsilon_k} u_{\epsilon_k})_{ X^{\frac{1}{2}}}|+|(v,\hat{E}_{\epsilon_k} u_{\epsilon_k})_{ X^{\frac{1}{2}}} - (v,u_0)_{ X^{\frac{1}{2}}}|.
$$

Since $\hat{E}_{\epsilon_k} u_{\epsilon_k} \wto u_0$ in $ X^{\frac{1}{2}}$ then $|(v,\hat{E}_{\epsilon_k} u_{\epsilon_k})_{ X^{\frac{1}{2}}} - (v,u_0)_{ X^{\frac{1}{2}}}| \to 0$. On the other hand,
\begin{equation}
\label{EQN3NOVO}
\begin{array}{lll}
|(E_{\epsilon_k}v,u_{\epsilon_k})_{ X^{\frac{1}{2}}_{\epsilon_k}} - (v,\hat{E}_{\epsilon_k} u_{\epsilon_k})_{ X^{\frac{1}{2}}}| = \left|\displaystyle \int_{\Omega_{\epsilon_k}}\!\! \nabla E_{\epsilon_k} v\nabla u_{\epsilon_k} +\int_{\Omega_{\epsilon_k}}\!\!  E_{\epsilon_k} v u_{\epsilon_k} - \int_{\Omega}\!\! \nabla v\nabla \hat{E}_{\epsilon_k} u_{\epsilon_k} -\int_{\Omega}\!\! v \hat{E}_{\epsilon_k} u_{\epsilon_k} \right|
\\
=
\left| \displaystyle \int_{\Omega_{\epsilon_k} \setminus K_{\epsilon_k}}\!\!\!\! (\nabla E_{\epsilon_k} v\nabla u_{\epsilon_k} + E_{\epsilon_k} v u_{\epsilon_k}) + \int_{ K_{\epsilon_k}}\!\!\!\! (\nabla E_{\epsilon_k} v\nabla u_{\epsilon_k}  +  E_{\epsilon_k} v u_{\epsilon_k}) - \displaystyle \int_{\Omega}\!\! (\nabla  v\nabla \hat{E}_{\epsilon_k} u_{\epsilon_k} + v \hat{E}_{\epsilon_k} u_{\epsilon_k}) \right|. 
\end{array}
\end{equation}
Notice that  (\ref{EQN3NOVO}) is the same as (\ref{EQN3}). Thus, the proof follows analogously to the Lemma \ref{lem_E_E}.
\cqd

%%%%%%%%%%%%%%%%%%%%%%%%%%%%%%%%%%%%%%%%%%%%%%%%%%%%%%%%%%%%%%%%%%%%%%%%%%%%%%%%%%

\section{$\mathcal{P}$-continuity of the attractors}
\label{P_continuity_attractors}

We have seen in  Subsection \ref{Carvalho-Piskarev} that the hypotheses {\bf[A1]} and {\bf[A2]} imply  the $\mathcal{P}^{\beta}$-continuity of the attractors in $Y^{\beta}_{\eps}$, for some $0\leq \beta<1$.  In this section we are going to define clearly the choice of our space $Y_\eps$ and $Y^\beta_\eps$ and prove that for this choice, the hypotheses \textbf{[A1]} and \textbf{[A2]} hold. This will imply the $\mathcal{P}^{\beta}$-continuity of the attractors of our problems (\ref{nbc}) and (\ref{nbc_limite_gamma_F}) or \eqref{nbcabstract}.   

If we fix a value $\alpha\in (0,1)$ and we choose  $Y_\eps=X^{-\frac{\alpha}{2}}_\eps$, $\eps\in [0,\eps_0]$, where $Y_0=Y$ and $X^{-\frac{\alpha}{2}}_{0}=X^{-\frac{\alpha}{2}}$, and consider them as the base spaces of the operators $A_\eps$. 

We define the operator $p_\eps:Y\to Y_\eps$, that is, $p_\eps:X^{-\frac{\alpha}{2}}\to X^{-\frac{\alpha}{2}}_{\epsilon}$ as $p_\eps =E^*_\epsilon$, where $E^*_\epsilon:X^{-\frac{\alpha}{2}}\to X^{-\frac{\alpha}{2}}_{\epsilon}$ is given by (\ref{extensaodual}). Moreover, choosing
$\beta=\frac{1+\alpha}{2}\in (\frac{1}{2},1)$ then we have $Y_\eps^{\frac{1+\alpha}{2}}=X^{\frac{1}{2}}_\eps=H^1(\Omega_\eps)$, $\epsilon \in[0,\epsilon_0]$,  and $p_\eps^{\frac{1+\alpha}{2}}:Y^{\frac{1+\alpha}{2}}\to Y^{\frac{1+\alpha}{2}}_\eps$ given by
$$
p_\eps^{\frac{1+\alpha}{2}} x=(A_\eps)^{-\frac{1+\alpha}{2}} E^*_\eps (A_0)^{\frac{1+\alpha}{2}} x, \quad \mbox{for $x\in Y^{\frac{1+\alpha}{2}}$.} 
$$
In particular, if $\alpha\in (\frac{1}{2},1)$ then $X^{\frac{\alpha}{2}}_{\eps}\equiv H^{\alpha}(\Omega_\eps)$ and therefore the trace operator from $X^{\frac{\alpha}{2}}_\eps$ to $L^2(\partial\Omega_\eps)$ is well defined and continuous for $\eps\in [0,\eps_0]$, see Lemmas \ref{equivalentdefinitionsspacesnorms} and \ref{embeddings-traces-uniformepsilon}. This is necessary in Subsection \ref{A2_nonlinearity} to deal with the nonlinearity at the boundary in the problems (\ref{nbc}) and (\ref{nbc_limite_gamma_F}).

We can distinguish within \textbf{[A1]} and \textbf{[A2]} hypotheses related to either resolvents or nonlinearities. So, in the next two subsections organized by those focus we are going to prove now several results which will conclude with a proof that hypotheses \textbf{[A1]} and \textbf{[A2]} hold with our notions of convergence and  operators defined in Subsection \ref{extensionoperators}. 
 Later, in Subsection \ref{EvsP},
 we will show the concepts of
$\mathcal{P}^{\frac{1+\alpha}{2}}$-convergence (see Definition \ref{Definition_Concept_P_Convergence}) and $E$-convergence (see Definition \ref{Definition_Concept_E_Convergence}) are equivalent and consequently we will conclude the
$E$-continuity of the attractors.

\subsection{Checking hypothesis {\bf[A1]} related to resolvent operators}
\label{A1_Resolvent}

 Initially, we note the operators $A_\eps$, $\eps \in [0,\eps_0]$, defined in Subsection \ref{equation-setting} are sectorial. 
  We will prove that  the resolvent operators are compact.   Moreover, we will prove the compact convergence of $A_\epsilon^{-1}$ to $A_0^{-1}$ (see Definition \ref{defcompactconvoperators}). 
  
\begin{lemma}
\label{compactlyconvergence}
Let $\{\Omega_\epsilon\}_{\epsilon \in [0,\epsilon_0]}$ be a family of domains satisfying conditions {\bf{(H)}}, {\bf{(F)(i)}} and  {\bf{(I)}}. If $0 <\alpha < 1$ then the family $\{ A_\epsilon^{-1} \in {\cal{L}}(X^{-\frac{\alpha}{2}}_{\epsilon},X^{\frac{1}{2}}_{\epsilon}): \epsilon \in (0, \epsilon_0]\}$  compactly converges to $A_0^{-1} \in {\cal{L}}(X^{-\frac{\alpha}{2}},X^{\frac{1}{2}})$, as $\epsilon \to 0$, that is, $A_\epsilon^{-1} {\stackrel{CC}{\longrightarrow}} A_0^{-1}$.
\end{lemma}
\proof
We are going to prove that: 

\vspace{0.2cm}

\noindent 1. $A_\epsilon^{-1}: X^{-\frac{\alpha}{2}}_{\epsilon} \to X^{\frac{1}{2}}_{\epsilon}$ is a compact operator, for $\epsilon \in [0, \epsilon_0]$ and $0<\alpha <1$;

\vspace{0.2cm}

\noindent 2.  $\{ A_\epsilon^{-1} h_\epsilon\}_{\epsilon \in (0, \epsilon_0]}$ is an $E$-precompact family whenever $\|h_\epsilon\|_{X^{-\frac{\alpha}{2}}_\epsilon}$ is bounded;

%\vspace{0.2cm}
 
\noindent 3. If $h_\epsilon {\stackrel{E^*}{\longrightarrow}} \ h_0$ then $A_\epsilon^{-1} h_\epsilon \dto A_0^{-1}h_0$.

\vspace{0.2cm}

Let us show each of the three points above.

\vspace{0.2cm}

\noindent 1.  Since $X^{1-\frac{\alpha}{2}}_\epsilon \hookrightarrow X^{\frac{1}{2}}_\epsilon$ is compact,  for $\alpha<1$, and $A_\epsilon^{-1} : X^{-\frac{\alpha}{2}}_\epsilon \to X^{1-\frac{\alpha}{2}}_\epsilon$ is continuous then $A_\epsilon^{-1} : X^{-\frac{\alpha}{2}}_\epsilon \to X^{\frac{1}{2}}_\epsilon$ is compact. 

\vspace{0.2cm}

\noindent  2. Let $\{h_\epsilon\}_{\epsilon \in (0, \epsilon_0]} $ be a family in $X^{-\frac{\alpha}{2}}_{\epsilon}$ such that $\|h_\epsilon\|_{X^{-\frac{\alpha}{2}}_{\epsilon}} \leq K$, for some $K>0$ independent of $\epsilon$, we will prove that there exists a  subsequence of  $\{A_\epsilon^{-1}h_\epsilon\}_{\epsilon \in (0, \epsilon_0]}$ which is $E$-convergent in $X^{\frac{1}{2}}_\epsilon$. 

For each $\epsilon \in (0, \epsilon_0]$, denote by $v_\epsilon=A_\epsilon^{-1}h_\epsilon \in X^{\frac{1}{2}}_\epsilon$, that is  $A_\epsilon v_\epsilon = h_\epsilon$. 
Note that $v_\eps\in X^{\frac{1}{2}}_\eps$ and $\|v_\epsilon\|^2_{X^{\frac{1}{2}}_\epsilon}$ is uniformly bounded. In fact, using Lemma \ref{embeddings-traces-uniformepsilon} to ensure that $X^{\frac{1}{2}}_\epsilon \hookrightarrow X^{\frac{\alpha}{2}}_\epsilon,$ with  embedding constant $c>0$ independent of $\epsilon$, we have \begin{equation}
\label{v_epsnorm} 
\begin{array} {lll}
\displaystyle \|v_\epsilon\|^2_{X^{\frac{1}{2}}_\epsilon} & = & \displaystyle \|v_\epsilon\|^2_{H^1(\Omega_\epsilon)}  = \int_{\Omega_\epsilon} |\nabla v_\epsilon |^2 + \int_{\Omega_\epsilon} |v_\epsilon|^2 = \langle A_\epsilon v_\epsilon, v_\epsilon\rangle = \langle h_\epsilon, v_\epsilon \rangle   \\
& \leq & \displaystyle \|h_\epsilon\|_{X^{-\frac{\alpha}{2}}_\epsilon} \|v_\epsilon\|_{X^{\frac{\alpha}{2}}_\epsilon} \leq c\|h_\epsilon\|_{X^{-\frac{\alpha}{2}}_\epsilon} \|v_\epsilon\|_{X^{\frac{1}{2}}_\epsilon}\leq cK \|v_\epsilon\|_{X^{\frac{1}{2}}_\epsilon}.
\end{array}
\end{equation}
Moreover, $\|\hat{E}_\epsilon v_\epsilon \|_{X^{\frac{1}{2}}} \leq C$, with $C>0$ independent of $\epsilon$. Thus,  by Lemma \ref{lemma4.3AB0}, we have that there exist a subsequence denoted by $v_{\epsilon_k}$  and $v_0 \in X^{\frac{1}{2}}$ such that 
\begin{eqnarray*}
\hat{E}_{\epsilon_k} v_{\epsilon_k}  \wto v_0  \quad \mbox{in  $X^{\frac{1}{2}}$},\qquad
\hat{E}_{\epsilon_k} v_{\epsilon_k}  \longrightarrow  v_0  \quad \mbox{in $X^{\frac{\alpha}{2}}$}\qquad \mbox{and} \qquad
v_{\epsilon_k}  \dwto v_0  \quad \mbox{ in $X^{\frac{1}{2}}_\epsilon$}.
\end{eqnarray*}

We need to prove that $A_{\epsilon_k}^{-1}h_{\epsilon_k} = v_{\epsilon_k} \dto v_0$. Since $v_{\epsilon_k} \dwto v_0$   then it is sufficient to prove that $\|v_{\epsilon_k}\|_{X^{\frac{1}{2}}_{\epsilon_k}} \to \|v_0\|_{X^{\frac{1}{2}}}$, as $k\to \infty$. 

Note that $\|h_\epsilon\|_{X^{-\frac{\alpha}{2}}_\epsilon} \leq K$ implies $\|\hat{E}_\epsilon^* h_\epsilon\|_{X^{-\frac{\alpha}{2}}} \leq \tilde{K}$. Using Banach-Alaouglu-Bourbaki Theorem, there exist $h_0 \in X^{-\frac{\alpha}{2}}$ and a subsequence $\hat{E}_{\epsilon_n}^* h_{\epsilon_n}$ such that, for all $\phi \in X^{\frac{\alpha}{2}}$,
\begin{equation} \label{hepsilon_convergence_on_Ephi}
\langle \hat{E}_{\epsilon_n}^* h_{\epsilon_n},\phi\rangle \to \langle h_0,\phi\rangle, \quad \mbox{as $n\to \infty$}.
\end{equation}

Now, we are going to prove that: 
\begin{equation}
\label{convergenceE*hat_epsheps}
\mbox{ If\ } \phi_{n}\in X^{\frac{\alpha}{2}} \mbox{ with\ } \phi_n \to \phi \mbox{ in  } X^{\frac{\alpha}{2}}, \mbox{ then  }
\langle \hat{E}_{\epsilon_n}^*h_{\epsilon_n},\phi_{n}\rangle \to \langle h_0,\phi\rangle, \quad \mbox{as $n\to \infty$}.
\end{equation}
In fact, using \eqref{hepsilon_convergence_on_Ephi} we obtain 
\begin{eqnarray*}
|\langle \hat{E}_{\epsilon_n}^*h_{\epsilon_n},\phi_{n}\rangle - \langle h_0,\phi\rangle| &\leq & |\langle\hat{E}_{\epsilon_n}^*h_{\epsilon_n},\phi_{n}\rangle - \langle\hat{E}_{\epsilon_n}^*h_{\epsilon_n},\phi\rangle| + |\langle\hat{E}_{\epsilon_n}^*h_{\epsilon_n},\phi\rangle - \langle h_0,\phi\rangle| \\ &\leq &
\|\hat{E}_{\epsilon_n}^*h_{\epsilon_n}\|_{X^{-\frac{\alpha}{2}}} \|\phi_{n} - \phi\|_{X^{\frac{\alpha}{2}}} + |\langle \hat{E}_{\epsilon_n}^*h_{\epsilon_n},\phi\rangle - \langle h_0,\phi\rangle| \to 0,\quad \mbox{as $n\to \infty$}.
\end{eqnarray*}

Taking subsequences and considering $\phi_{k}=\hat{E}_{\epsilon_k}v_{\epsilon_k} $ and $\phi=v_0$, we have $\phi_{k} \to \phi$ in $X^{\frac{\alpha}{2}}$. Thus, by \eqref{convergenceE*hat_epsheps} we obtain
$$
\langle \hat{E}_{\epsilon_k}^*h_{\epsilon_k},\hat{E}_{\epsilon_k}v_{\epsilon_k}\rangle \to \langle h_0,v_0\rangle, \quad \mbox{as $k\to \infty$}.
$$
Now, since $\langle \hat{E}^*_{\epsilon_k}h_{\epsilon_k},\hat{E}_{\epsilon_k}v_{\epsilon_k}\rangle =\langle h_{\epsilon_k},E_{\epsilon_k}\hat{E}_{\epsilon_k}v_{\epsilon_k}\rangle$ and, by Lemma \ref{convergence_extension_ofthe_restriction}, 
$\|E_{\epsilon_k}\hat{E}_{\epsilon_k}v_{\epsilon_k}- v_{\epsilon_k}\|_{X^{\frac{\alpha}{2}}_{\epsilon_k}} \to 0$, as $k \to \infty$, then  
$$
%\label{convh_eps}
\begin{array}{lll}
|\langle h_{\epsilon_k},v_{\epsilon_k}\rangle -\langle h_0,v_0\rangle| &\leq & |\langle h_{\epsilon_k},v_{\epsilon_k}\rangle-\langle h_{\epsilon_k},E_{\epsilon_k}\hat{E}_{\epsilon_k}v_{\epsilon_k}\rangle|+|\langle\hat{E}^*_{\epsilon_k}h_{\epsilon_k},\hat{E}_{\epsilon_k}v_{\epsilon_k}\rangle -\langle h_0,v_0\rangle|\\
&\leq & \|h_{\epsilon_k}\|_{X^{-\frac{\alpha}{2}}_{\epsilon_k}} \|E_{\epsilon_k}\hat{E}_{\epsilon_k}v_{\epsilon_k}- v_{\epsilon_k}\|_{X^{\frac{\alpha}{2}}_{\epsilon_k}} + |\langle \hat{E}^*_{\epsilon_k}h_{\epsilon_k},\hat{E}_{\epsilon_k}v_{\epsilon_k}\rangle-\langle h_0,v_0\rangle|  \to 0.
\end{array}
$$
Thus,
$$
\langle h_{\epsilon_k}, v_{\epsilon_k}\rangle \to \langle h_0,v_0 \rangle, \quad \mbox{as $k\to \infty$}.
$$

By \eqref{v_epsnorm} we have $\|v_{\epsilon_k}\|^2_{X^{\frac{1}{2}}_{\epsilon_k}} =  \langle h_{\epsilon_k}, v_{\epsilon_k}\rangle$,
then 
\begin{equation}
\label{convergencevepsilontoh0v0}
\|v_{\epsilon_k}\|^2_{X^{\frac{1}{2}}_{\epsilon_k}} \to \langle h_0,v_0 \rangle.
\end{equation}
Finally, we will prove that $\langle h_0,v_0 \rangle = \|v_0\|_{X^{\frac{1}{2}}}^2$ and $v_0=A_0^{-1} h_0$. We note
\begin{equation}
\label{v0norm}
\|v_0\|^2 _{X^{\frac{1}{2}}} = \|v_0\|^2 _{H^{1}(\Omega)} = \int_{\Omega} |\nabla v_0|^2 + \int_{\Omega} |v_0|^2 = \langle A_0v_0,v_0\rangle.
\end{equation}

On the one hand, $v_{\epsilon_k} \dwto v_0$, then
$$
(v_{\epsilon_k}, E_{\epsilon_k}w)_{X^{\frac{1}{2}}_{\epsilon_k}} \to  (v_{0}, w)_{X^{\frac{1}{2}}}, \quad \mbox{for all $w\in X^{\frac{1}{2}}$}.
$$
On the other hand, by \eqref{hepsilon_convergence_on_Ephi} we have
$$
(v_{\epsilon_k}, E_{\epsilon_k}w)_{X^{\frac{1}{2}}_{\epsilon_k}} = (A^{-1}_{\epsilon_k}h_{\epsilon_k}, E_{\epsilon_k}w)_{X^{\frac{1}{2}}_{\epsilon_k}} = \langle h_{\epsilon_k}, E_{\epsilon_k}w\rangle = \langle \hat{E}_{\epsilon_k}^* h_{\epsilon_k},w\rangle \to \langle h_{0}, w\rangle, \quad \mbox{for all $w\in X^{\frac{1}{2}}$}.
$$
Hence and the uniqueness of the limit, we have
$$
\langle h_0,w\rangle = (v_0,w)_{X^{\frac{1}{2}}}= \langle A_0v_0,w\rangle, \quad \mbox{ for all $w\in X^{\frac{1}{2}}$}.
$$
Thus,
\begin{equation} 
\label{v0solution}
v_0=A_0^{-1} h_0 \quad \mbox{or} \quad A_0 v_0=h_0.
\end{equation} 

By \eqref{convergencevepsilontoh0v0}, \eqref{v0norm} and \eqref{v0solution}, we obtain
$$
\|v_{\epsilon_k}\|^2_{X^{\frac{1}{2}}_{\epsilon_k}}\to \|v_{0}\|^2_{X^{\frac{1}{2}}}.
$$

\noindent  3. Since $h_\epsilon {\stackrel{E^*}{\longrightarrow}} h_0$ then $\|h_\epsilon \|_{X^{-\frac{\alpha}{2}}_\epsilon}\leq K$, for some $K>0$ independent of $\epsilon$. Thus, for any sequence $\epsilon_k \to 0$, we can extract another subsequence, which we denote
also as $\epsilon_k$, such that following the argument made above to prove item 2, we have $A_{\epsilon_k}^{-1} h_{\epsilon_k} \dto A_0^{-1}h_0$.
Since this has been proved for any sequence, we obtain the $E$-convergence for the
whole family, that is, $A_\epsilon^{-1} h_\epsilon \dto A_0^{-1}h_0$.  \cqd

Now, we verify that the resolvents $(\mu I+A_\eps)^{-1}$ compactly converge to $(\mu I+A_0)^{-1}$.

\begin{lemma}
\label{lemma_compactlyconvergence_limita_resolvente1}
Let $\{\Omega_\epsilon\}_{\epsilon \in [0,\epsilon_0]}$ be a family of domains satisfying conditions {\bf{(H)}}, {\bf{(F)(i)}} and  {\bf{(I)}}. If $0<\alpha<1$ and $\mu \in \rho(-A_0)$, then there exists $\epsilon_\mu >0$ such that $\mu \in \rho(-A_\epsilon)$, for all $\epsilon \in [0,\epsilon_\mu]$,  and there exists a constant $M_\mu>0$ independent of $\epsilon$ such that 
\begin{equation}
\label{NU_limita_resolvente}
\|(\mu I+A_\epsilon)^{-1}\|_{{\cal{L}}( X^{-\frac{\alpha}{2}}_\epsilon, X^{\frac{1}{2}}_\epsilon )} \leq M_\mu, \quad \mbox{for all $\epsilon \in [0,\epsilon_\mu]$}.
\end{equation}
%$$
%\sup_{\epsilon \in [0,\epsilon_\mu]}  \|(\mu+A_\epsilon)^{-1}\|_{{\cal{L}}( X^{-\frac{\alpha}{2}}_\epsilon, X^{\frac{1}{2}}_\epsilon )} <\infty.
%$$
Furthermore, $(\mu I + A_\epsilon)^{-1} {\stackrel{CC}{\longrightarrow}} (\mu I+A_0)^{-1}$.
\end{lemma}
\proof
Since $\mu \in \rho(-A_0)$ then
$$(\mu I+A_0)^{-1} = [(\mu A_0^{-1}+I)A_0)]^{-1} = A_0^{-1}(\mu A_0^{-1} + I)^{-1} = A_0^{-1}(I+\mu A_0^{-1})^{-1}.$$
Thus, $(I+\mu A_0^{-1})^{-1}$ exists and $(I+\mu A_0^{-1})$ is one-to-one and onto, thus $\ker(I+\mu A_0^{-1})=\{0\}$. 
By Lemma \ref{compactlyconvergence}, $A_\epsilon^{-1} {\stackrel{CC}{\longrightarrow}} A_0^{-1}$, then using \cite[Lemma 3.1]{AB0}, we obtain that there exist $\epsilon_\mu >0$ and  constants $\tilde{M}_\mu$, $C>0$ independents of $\epsilon$ such that  for all $\epsilon \in [0,\epsilon_\mu]$, $\|A_\epsilon^{-1}\|_{{\cal{L}}(X^{-\frac{\alpha}{2}}_\epsilon, X^{\frac{1}{2}}_\epsilon)} \leq C$, $\mu \in \rho(-A_\epsilon)$  and 
\begin{equation}
\label{bounded_resolvent_part}
\|(I+\mu A_\epsilon^{-1})^{-1}\|_{{\cal{L}}(X^{-\frac{\alpha}{2}}_\epsilon)} \leq \tilde{M}_\mu.
\end{equation}
Therefore, $A_\epsilon^{-1}(I+\mu A_\epsilon^{-1})^{-1}$ is well defined from $X^{-\frac{\alpha}{2}}_\epsilon$ to $X^{\frac{1}{2}}_\epsilon$, and it is bounded. In fact, 
$$
\|A_\epsilon^{-1}(I+\mu A_\epsilon^{-1})^{-1}\|_{{\cal{L}}(X^{-\frac{\alpha}{2}}_\epsilon,X^{\frac{1}{2}}_\epsilon)} \leq \|A_\epsilon^{-1}\|_{{\cal{L}}(X^{-\frac{\alpha}{2}}_\epsilon, X^{\frac{1}{2}}_\epsilon)} \|(I+\mu A_\epsilon^{-1})^{-1}\|_{{\cal{L}}(X^{-\frac{\alpha}{2}}_\epsilon)} \leq C \tilde{M}_\mu.
$$

Note that
$$
A_\epsilon^{-1}(I+\mu A_\epsilon^{-1})^{-1}=A_\epsilon^{-1}(\mu A_\epsilon^{-1}+I)^{-1}=[(\mu A_\epsilon^{-1}+I)A_\epsilon]^{-1} = (\mu I+A_\epsilon)^{-1}.
$$
Thus, there exist $\epsilon_\mu>0$ and $M_\mu >0$ independent of $\epsilon$ such that, for all $\epsilon \in [0,\epsilon_\mu]$, $\mu \in \rho(-A_\epsilon)$  and \eqref{NU_limita_resolvente} holds.
%$$
%\|(\mu I+A_\epsilon)^{-1}\|_{{\cal{L}}(X^{-\frac{\alpha}{2}}_\epsilon, X^{\frac{1}{2}}_\epsilon)}\leq M_\mu.
%$$ 

In order to show $(\mu I + A_\epsilon)^{-1} {\stackrel{CC}{\longrightarrow}} (\mu I+A_0)^{-1}$, we need to prove that:
 
\vspace{0.2cm}

\noindent 1. $(\mu I + A_\epsilon)^{-1}: X^{-\frac{\alpha}{2}}_{\epsilon} \to X^{\frac{1}{2}}_{\epsilon}$ is a compact operator, for $\epsilon \in [0, \epsilon_0]$ and $0<\alpha <1$;

\vspace{0.2cm}

\noindent 2. $\{ (\mu  I + A_\epsilon)^{-1} h_\epsilon\}_{\epsilon \in (0, \epsilon_0]}$ is an $E$-precompact family whenever $\|h_\epsilon\|_{X^{-\frac{\alpha}{2}}_\epsilon}$ is bounded;

%\vspace{0.2cm}
 
\noindent 3. If $h_\epsilon {\stackrel{E^*}{\longrightarrow}} h_0$ then $(\mu I + A_\epsilon)^{-1} h_\epsilon \dto (\mu I + A_0)^{-1} h_0$.

\vspace{0.2cm}

Let us show each of the three points above.

\vspace{0.2cm}

\noindent 1. The proof is similar to the item 2 that we will prove below.

\vspace{0.2cm}

\noindent 2. Let $\{h_\epsilon\}_{\epsilon \in (0, \epsilon_0]} $ be a family in $X^{-\frac{\alpha}{2}}_{\epsilon}$ such that $\|h_\epsilon\|_{X^{-\frac{\alpha}{2}}_{\epsilon}} \leq K$, for some $K>0$ independent of $\epsilon$. Note that
$$
(\mu I + A_\epsilon)^{-1} h_\epsilon = A_\epsilon^{-1}(I+\mu A_\epsilon^{-1})^{-1} h_\epsilon.
$$
Using \eqref{bounded_resolvent_part} we obtain $\{(I+\mu A_\epsilon^{-1})^{-1} h_\epsilon\}_{\epsilon \in (0, \epsilon_0]}$ is a family in $X^{-\frac{\alpha}{2}}_{\epsilon}$ which is uniformly bounded in $\epsilon$. Hence, from Lemma \ref{compactlyconvergence}, we get that $\{A_\epsilon^{-1} (I+\mu A_\epsilon^{-1})^{-1} h_\epsilon\}_{\epsilon \in (0, \epsilon_0]}$ is an $E$-precompact family in $X^{\frac{1}{2}}_\epsilon$. Consequently, $\{ (\mu I + A_\epsilon)^{-1} h_\epsilon\}_{\epsilon \in (0, \epsilon_0]}$ is an $E$-precompact family in $X^{\frac{1}{2}}_\epsilon$.

\vspace{0.2cm}

\noindent 3. If $h_\epsilon {\stackrel{E^*}{\longrightarrow}} h_0$ then $\|h_\epsilon\|_{X^{-\frac{\alpha}{2}}_\epsilon} \leq K$, for some $K>0$ independent of $\epsilon$. Now, from item 2. for any subsequence of $\{ (\mu I + A_\epsilon)^{-1} h_\epsilon\}_{\epsilon \in (0, \epsilon_0]}$ there exist a subsequence $\{ (\mu I + A_{\epsilon_k})^{-1} h_{\epsilon_k}\}_{k \in \mathbb{N}}$ and $y\in X^{\frac{1}{2}}$ such that $(\mu I + A_{\epsilon_k})^{-1}h_{\epsilon_k} {\stackrel{E}{\longrightarrow}} y$ in $X^{\frac{1}{2}}_{\epsilon_k}$. Taking $z_{\epsilon_k} = (I+\mu  A^{-1}_{\epsilon_k})^{-1} A^{-1}_{\epsilon_k}h_{\epsilon_k} $ we have
\begin{equation}
\label{res1}
z_{\epsilon_k} = (I+\mu  A^{-1}_{\epsilon_k})^{-1} A^{-1}_{\epsilon_k}h_{\epsilon_k} = (\mu I+ A_{\epsilon_k})^{-1}h_{\epsilon_k} {\stackrel{E}{\longrightarrow}} y, \quad \mbox{in $X^{\frac{1}{2}}_{\epsilon_k}$}.
\end{equation}

Using again Lemma \ref{compactlyconvergence}, we have 
\begin{equation}
\label{res2}
A^{-1}_{\epsilon_k}h_{\epsilon_k} {\stackrel{E}{\longrightarrow}} A^{-1}_{0}h_0, \quad \mbox{in $X^{\frac{1}{2}}_{\epsilon_k}$}.
\end{equation}
On the other hand, since $z_{\epsilon_k} {\stackrel{E}{\longrightarrow}} y$ in $X^{\frac{1}{2}}_{\epsilon_k}$, then by Lemma \ref{lem_E_E} $z_{\epsilon_k} {\stackrel{E^*}{\longrightarrow}} y$ in $X^{-\frac{\alpha}{2}}_{\epsilon_k}$. Using again Lemma \ref{compactlyconvergence}, we obtain
\begin{equation}
\label{res3}
A^{-1}_{\epsilon_k}z_{\epsilon_k} {\stackrel{E}{\longrightarrow}} A^{-1}_{0} y, \quad \mbox{in $X^{\frac{1}{2}}_{\epsilon_k}$}.
\end{equation}

Note that $A^{-1}_{\epsilon_k}h_{\epsilon_k} = (I+\mu A^{-1}_{\epsilon_k})z_{\epsilon_k}=z_{\epsilon_k}+\mu A^{-1}_{\epsilon_k} z_{\epsilon_k}$. Now, from \eqref{res1} and \eqref{res3}, we have 
\begin{equation}
\label{res4}
A^{-1}_{\epsilon_k}h_{\epsilon_k} {\stackrel{E}{\longrightarrow}} y+ \mu A^{-1}_{0}y=(I+\mu A^{-1}_{0})y, \quad \mbox{in $X^{\frac{1}{2}}_{\epsilon_k}$}.
\end{equation}

Using \eqref{res2}, \eqref{res4} and uniqueness of the limit, we have $A^{-1}_{0}h_0=(I+\mu A^{-1}_{0})y$. Thus, $y=(\mu I +A_0)^{-1}h_0$ and
$$
(\mu I + A_{\epsilon_k})^{-1} h_{\epsilon_k} \dto (\mu I + A_0)^{-1} h_0, \quad \mbox{in $X^{\frac{1}{2}}_{\epsilon_k}$}.
$$
In particular, $y$ is independent of the subsequence chosen. This implies that the whole sequence  $(\mu I + A_{\epsilon})^{-1} h_{\epsilon}$ $E$-converges to $y=(\mu I + A_0)^{-1} h_0$. Thus $(\mu I + A_\epsilon)^{-1} {\stackrel{CC}{\longrightarrow}} (\mu I+A_0)^{-1}$.
\cqd

Moreover, we can obtain that the boundedness \eqref{NU_limita_resolvente} in   Lemma \ref{lemma_compactlyconvergence_limita_resolvente1} is uniform.

\begin{lemma}
\label{lemma_limita_resolvente2}
Let $\{\Omega_\epsilon\}_{\epsilon \in [0,\epsilon_0]}$ be a family of domains satisfying conditions {\bf{(H)}}, {\bf{(F)(i)}} and  {\bf{(I)}}. Let $0<\alpha<1$ and $K$ be a compact set, $K \subset \rho(-A_0)$, then there exists $\epsilon_K>0$  such that $K \subset \rho(-A_\epsilon)$, for all  $\epsilon \in [0,\epsilon_K]$, and 
\begin{equation}
\label{uniformly_bounded_resolvent_epsilon1}
\sup_{\alpha\in(0,1)} \sup_{\epsilon \in [0,\epsilon_K]} \sup_{\mu \in K} \|(\mu I + A_\epsilon)^{-1}\|_{{\cal{L}}(X^{-\frac{\alpha}{2}}_\epsilon,X^{\frac{1}{2}}_\epsilon)} < \infty.
\end{equation}
\end{lemma}
\proof
From Lemma \ref{compactlyconvergence} we have that $A_\epsilon^{-1} {\stackrel{CC}{\longrightarrow}} A_0^{-1}$. Hence, by \cite[Proposition 3.3]{AB0}, there exists $\epsilon_K>0$ such that $K \subset \rho(-A_\epsilon)$, for all $\epsilon \in [0,\epsilon_K]$, and by \cite[Lemma 3.1]{AB0}, we obtain that there exists a constant  $C>0$ independent of $\epsilon$ such that $\|A_\epsilon^{-1}\|_{{\cal{L}}(X^{-\frac{\alpha}{2}}_\epsilon, X^{\frac{1}{2}}_\epsilon)} \leq C$. Since  $(\mu I + A_\epsilon)^{-1}=A_\epsilon^{-1} (I+\mu A_\epsilon^{-1})^{-1}$ then, in order to prove the uniformly boundedness \eqref{uniformly_bounded_resolvent_epsilon1} it is sufficient to prove 
\begin{equation}
\label{uniformly_bounded_resolvent_epsilon2}
\sup_{\alpha\in(0,1)} \sup_{\epsilon \in [0,\epsilon_K]} \sup_{\mu \in K} \|(I +\mu A_\epsilon^{-1})^{-1}\|_{{\cal{L}}(X^{-\frac{\alpha}{2}}_\epsilon)} < \infty.
\end{equation}

Suppose that \eqref{uniformly_bounded_resolvent_epsilon2} is not uniformly bounded, thus there exist sequences $\epsilon_n \to 0$ and $\mu_n \in K$ such that $\mu_n \to \bar{\mu}\in K$ and 
$$ 
\|(I +\mu_n A_{\epsilon_n}^{-1})^{-1}\|_{{\cal{L}}(X^{-\frac{\alpha}{2}}_{\epsilon_n})} \to \infty, \quad \mbox{as $n\to \infty$}.
$$
Since $\mu_n A_{\epsilon_n}^{-1} {\stackrel{CC}{\longrightarrow}} \bar{\mu} A_0^{-1}$ then using again \cite[Lemma 3.1]{AB0}, we obtain that there exist $\epsilon_0 >0$ and a constant $M>0$ such that $\|(I+\mu_n A_{\epsilon_n}^{-1})^{-1}\|_{{\cal{L}}(X^{-\frac{\alpha}{2}}_{\epsilon_n})} \leq M$, for all $\epsilon_n \in [0,\epsilon_0]$, which is a contradiction. 
\cqd

Therefore, the results above imply
 the conditions on the resolvent in the hypothesis   {\bf{[A1]}} hold.

%%%%%%%%%%%%%%%%%%

\subsection{Checking hypotheses {\bf[A1]} and {\bf[A2]} related to nonlinearities}
\label{A2_nonlinearity}

%%%%%%%%%%%%%%%%%%%%%%%%%%%%%%%%

In order to prove the $\mathcal{P}^{\frac{1+\alpha}{2}}$-continuity of the attractors, we need to prove the hypotheses  \textbf{[A1]} and \textbf{[A2]} in Subsection \ref{Carvalho-Piskarev} related to  nonlinearities. 

Initially, we will prove that the abstract nonlinearities given by \eqref{hepsilon} and \eqref{hzero} are bounded and globally Lipschitz.
 
\begin{lemma}
\label{limitacaoinfinito}
Suppose that $f$ and $g$ satisfy \eqref{bounded_f} and \eqref{bounded_g} and $\frac{1}{2}<\alpha \leq 1$. Then there exists $K>0$ independent of $\epsilon$ such that 
$$
\left\| h_{\epsilon}(u_\epsilon)\right\|_{X^{-\frac{\alpha}{2}}_\epsilon} \leq K, \quad \mbox{for all $u_\epsilon \in X^{\frac{1}{2}}_\epsilon$ and $0\leq \epsilon \leq \epsilon_{0}$}.
$$
\end{lemma}
\proof Initially, for each $\epsilon\in [0,\epsilon_0]$ we have
$$
\|h_{\epsilon}(u_\epsilon)\|_{X^{-\frac{\alpha}{2}}_\epsilon}=\sup_{\begin{array}{ccc}\psi_{\epsilon}\in X^{\frac{\alpha}{2}}_\epsilon\\ \|\psi_{\epsilon}\|_{X^{\frac{\alpha}{2}}_{\epsilon}}=1\end{array}}  |\langle h_{\epsilon}(u_{\epsilon}),\psi_{\epsilon}\rangle |.
$$

Since that $f$ and $g$ satisfy \eqref{bounded_f} and \eqref{bounded_g}, then there exists $C>0$ independent of $\epsilon$ such that
\begin{equation}
\label{limitacao_fg}
\| f(\cdot,u_{\epsilon}(\cdot))\|_{L^{\infty}(\Omega_{\epsilon})}\leq C \quad \mbox{and} \quad \| g(\cdot,u_{\epsilon}(\cdot))\|_{L^{\infty}(\Omega_{\epsilon})}\leq C, \quad \mbox{for all  $\epsilon \in [0,\epsilon_0]$}.
\end{equation}

For $\epsilon=0$, using (\ref{limitacao_fg}), $\gamma \in L^{\infty}(\partial \Omega)$, Cauchy-Schwartz inequality and Lemma \ref{embeddings-traces-uniformepsilon},  we get 
$$
\begin{array}{lll}
\displaystyle |\langle h_{0}(u_{0}),\psi_{0}\rangle | & \leq & \displaystyle \int_{\Omega}|f(x,u_0)||\psi_0|+\int_{\partial \Omega} |\gamma(x)||g(x,u_0)||\psi_0| \\
& \leq & \displaystyle C |\Omega|^{\frac{1}{2}} \| \psi_0\|_{L^{2}(\Omega)}+C\|\gamma\|_{L^{\infty}(\partial \Omega)} |\partial \Omega|^{\frac{1}{2}} \|\psi_0\|_{L^{2}(\partial \Omega)}\\
& \leq & \displaystyle (C_1 |\Omega|^{\frac{1}{2}}+C_2 |\partial \Omega|^{\frac{1}{2}} \|\gamma\|_{L^{\infty}(\partial \Omega)})  \| \psi_0\|_{X^{\frac{\alpha}{2}}}, 
\end{array}
$$
for all $\psi_0\in X^{\frac{\alpha}{2}}$.  Thus, there exists $K_1>0$ such that
\begin{equation}
\label{limitacao_fg_1}
\|h_{0}(u_0)\|_{X^{-\frac{\alpha}{2}}}\leq K_1.
\end{equation}

Now for $\epsilon \in (0,\epsilon_0]$, using (\ref{limitacao_fg}), Cauchy-Schwartz inequality,  Lemma \ref{embeddings-traces-uniformepsilon} and uniformly bounded of domains,  we get 
$$
\begin{array}{lll}
\displaystyle |\langle h_{\epsilon}(u_{\epsilon}),\psi_{\epsilon}\rangle | & \leq & \displaystyle \int_{\Omega_\epsilon}|f(x,u_\epsilon)||\psi_\epsilon|+\int_{\partial \Omega_\epsilon} |g(x,u_\epsilon)||\psi_\epsilon| \\
& \leq & \displaystyle C |\Omega_\epsilon|^{\frac{1}{2}} \| \psi_\epsilon\|_{L^{2}(\Omega_\epsilon)}+C |\partial \Omega_\epsilon|^{\frac{1}{2}} \|\psi_\epsilon\|_{L^{2}(\partial \Omega_\epsilon)}\\
& \leq & \displaystyle (\tilde{C}_1 |\Omega_\epsilon|^{\frac{1}{2}}+\tilde{C}_2 |\partial \Omega_\epsilon|^{\frac{1}{2}} )  \| \psi_\epsilon\|_{X^{\frac{\alpha}{2}}_\epsilon}, \end{array}
$$
for all $\psi_\epsilon\in X^{\frac{\alpha}{2}}_\epsilon$. 
 Thus, there exists $K_2>0$ independent of $\epsilon$ such that
\begin{equation}
\label{limitacao_fg_2}
\|h_{\epsilon}(u_\epsilon)\|_{X^{-\frac{\alpha}{2}}_\epsilon}\leq K_2.
\end{equation}

Therefore, the result follows from \eqref{limitacao_fg_1} and \eqref{limitacao_fg_2}. \cqd

\vspace{0.2cm}

\begin{lemma}
\label{lipschitz_h}
Suppose that $f$ and $g$ satisfy \eqref{bounded_f} and \eqref{bounded_g} and $\frac{1}{2}<\alpha \leq 1$. Then, for each $0\leq \eps\leq \eps_0$, the map $h_\epsilon : X^{\frac{1}{2}}_\epsilon \to X^{-\frac{\alpha}{2}}_\epsilon$ is globally Lipschitz, that is,
 there exists $L>0$ independent of $\epsilon$ such that 
$$
\left\| h_{\epsilon}(u_\epsilon)-h_{\epsilon}(v_\epsilon)\right\|_{X^{-\frac{\alpha}{2}}_\epsilon} \leq L \|u_{\epsilon}-v_{\epsilon} \|_{X^{\frac{1}{2}}_\epsilon}, \quad \mbox{for all $u_{\epsilon}, v_\epsilon \in X^{\frac{1}{2}}_\epsilon$ and $0\leq \epsilon \leq \epsilon_{0}$}.
$$
\end{lemma}
\proof Initially, for each $\epsilon\in [0,\epsilon_0]$ we have
$$
\|h_{\epsilon}(u_\epsilon)-h_{\epsilon}(v_\epsilon)\|_{X^{-\frac{\alpha}{2}}_\epsilon}=\sup_{\begin{array}{ccc}\psi_{\epsilon}\in X^{\frac{\alpha}{2}}_\epsilon\\ \|\psi_{\epsilon}\|_{X^{\frac{\alpha}{2}}_{\epsilon}}=1\end{array}}  |\langle h_{\epsilon}(u_\epsilon) - h_{\epsilon}(v_{\epsilon}),\psi_{\epsilon}\rangle |.
$$

For each $\eps\in [0,\eps_0]$, using Mean Value Theorem, we can write
$$
f(x,u_{\epsilon}(x))-f(x,v_{\epsilon}(x))=\partial_uf(x,\tilde{w}_\eps(x)u_{\epsilon}(x)+(1-\tilde{w}_{\eps}(x))v_{\epsilon}(x))[u_{\epsilon}(x)-v_{\epsilon}(x)], \quad x\in \Omega_\eps,
$$
and
$$
g(x,u_{\epsilon}(x))-g(x,v_{\epsilon}(x))=\partial_ug(x,\hat{w}_\eps(x) u_{\epsilon}(x)+(1-\hat{w}_\eps(x))v_{\epsilon}(x))[u_{\epsilon}(x)-v_{\epsilon}(x)], \quad x\in \partial \Omega_\eps,
$$
where $\tilde{w}_\epsilon(x)\in [0,1]$ for all $x \in \Omega_\epsilon$ and $\hat{w}_\epsilon(x)\in [0,1]$ for all $x \in \partial \Omega_\epsilon$.

For $\epsilon=0$, using  that $f$ and $g$ satisfy \eqref{bounded_f} and \eqref{bounded_g}, $\gamma \in L^{\infty}(\partial \Omega)$, Cauchy-Schwartz inequality and  Lemma \ref{embeddings-traces-uniformepsilon},  we get 
$$
\begin{array}{lll}
\displaystyle |\langle h_{0}(u_{0})-h_{0}(v_{0}),\psi_{0}\rangle | & \leq & \displaystyle \int_{\Omega}|f(x,u_0)-f(x,v_0)||\psi_0|+\int_{\partial \Omega} |\gamma(x)||g(x,u_0)-g(x,v_0)||\psi_0| \\
%& \leq & \displaystyle K_1\int_{\Omega}|u_0-v_0||\phi_0|+K_2\|\gamma\|_{L^{\infty}(\partial \Omega)} \int_{\partial \Omega} |u_0-v_0||\phi_0| \\
& \leq & \displaystyle C \| u_0-v_0\|_{L^{2}(\Omega)} \| \psi_0\|_{L^{2}(\Omega)}+C \|\gamma\|_{L^{\infty}(\partial \Omega)} \| u_0-v_0\|_{L^{2}(\partial \Omega)} \|\psi_0\|_{L^{2}(\partial \Omega)}\\
& \leq & \displaystyle (C_1 +C_2 \|\gamma\|_{L^{\infty}(\partial \Omega)}) \| u_0-v_0\|_{X^{\frac{1}{2}}}  \| \psi_0\|_{X^{\frac{\alpha}{2}}}, 
\end{array}
$$
for all $\psi_0\in X^{\frac{\alpha}{2}}$.
Thus, there exists $L_1>0$ such that
\begin{equation}
\label{lipschitziana_fg_1}
\|h_{0}(u_0)-h_{0}(v_0)\|_{X^{-\frac{\alpha}{2}}}\leq L_1 \| u_0-v_0\|_{X^{\frac{1}{2}}}.
\end{equation}

Now for $\epsilon \in (0,\epsilon_0]$, using that $f$ and $g$ satisfy \eqref{bounded_f} and \eqref{bounded_g}, Cauchy-Schwartz inequality and  Lemma \ref{embeddings-traces-uniformepsilon},  we get 
$$
\begin{array}{lll}
\displaystyle |\langle h_{\epsilon}(u_{\epsilon})-h_{\epsilon}(v_{\epsilon}),\psi_{\epsilon}\rangle | & \leq & \displaystyle \int_{\Omega_\epsilon}|f(x,u_\epsilon)-f(x,v_\epsilon)||\psi_\epsilon|+\int_{\partial \Omega_\epsilon} |g(x,u_\epsilon)-g(x,v_\epsilon)||\psi_\epsilon| \\
& \leq & \displaystyle C \| u_\epsilon-v_\epsilon\|_{L^{2}(\Omega_{\epsilon})} \| \psi_\epsilon\|_{L^{2}(\Omega_\epsilon)}+C \| u_\epsilon-v_\epsilon\|_{L^{2}(\partial \Omega_{\epsilon})} \|\psi_\epsilon\|_{L^{2}(\partial \Omega_\epsilon)}\\
& \leq & \displaystyle (\tilde{C}_1+\tilde{C}_2) \| u_\epsilon-v_\epsilon\|_{X^{\frac{1}{2}}_\epsilon}  \| \psi_\epsilon\|_{X^{\frac{\alpha}{2}}_\epsilon},  \end{array}
$$
for all $\psi_\epsilon\in X^{\frac{\alpha}{2}}_\epsilon$. 
 Thus, there exists $L_2>0$ independent of $\epsilon$ such that
\begin{equation}
\label{lipschitziana_fg_2}
\|h_{\epsilon}(u_\epsilon)-h_{\epsilon}(v_\epsilon)\|_{X^{-\frac{\alpha}{2}}_\epsilon}\leq L_2 \| u_\epsilon-v_\epsilon\|_{X^{\frac{1}{2}}_\epsilon}.
\end{equation}

Therefore, the result follows from \eqref{lipschitziana_fg_1} and \eqref{lipschitziana_fg_2}.
\cqd

To simplify the proof of results, from now on
 we will assume that the nonlinearities  
 depend only on $u_{\eps}$, $\eps \in [0,\eps_0]$, that is, we will consider
 $f(u_\eps)$ and $g(u_\eps)$ in the problems (\ref{nbc}) and (\ref{nbc_limite_gamma_F}) instead of $f(x,u_\eps)$ and $g(x,u_\eps)$.

We need to verify the $E^*$-convergence of the nonlinearities $h_\epsilon$ 
 given by \eqref{hepsilon} and \eqref{hzero}.  For this, we need the following result

\begin{lemma} \label{weakconvergenceissufficient}
Let $Q\subset \mathbb{R}^N$ be an open bounded set and $\psi_\epsilon, \ \psi \in L^p(Q)$, $0 < \epsilon \leq \epsilon_0$, such that  $\|\psi_\epsilon\|_{L^{p}(Q)} \leq K$ and $\|\psi\|_{L^p(Q)} \leq K$, for some $K>0$ independent of $\eps$. 
If $\psi_\epsilon \wto \psi$ in $L^{p}(Q)$ then 
$$
\left| \int_Q \vartheta_\eps (x) (\psi_\epsilon(x)-\psi(x))dx\right| = c(\epsilon, \vartheta_\eps) \to 0,
$$
uniformly for $\vartheta_\eps$ in compact sets of $L^{p^\prime}(Q)$, where $\frac{1}{p}+\frac{1}{p'}=1$. 
\end{lemma}

\proof  Suppose that the convergence is not uniform in compact sets of $L^{p^\prime}(Q)$.
Therefore, there exist a compact set $C$ of $L^{p^\prime}(Q)$, $\eta_0>0$  and  a sequence $\vartheta_\epsilon \in C$ such that
\begin{equation} \label{notuniform}
\left|\int_R \vartheta_\epsilon(x) (\psi_\epsilon(x)-\psi(x))dx \right| = c(\epsilon, \vartheta_\epsilon) \geq \eta_0, \  \mbox{ as } \epsilon \to 0.
\end{equation}
We are going to prove that \eqref{notuniform} is not possible. 

In fact, since $\vartheta_\epsilon \in C$, if necessary there exists a subsequence of $\vartheta_\epsilon$ and a function $\vartheta \in C$ such that $\vartheta_\epsilon \to \vartheta$ in $L^{p^\prime}(Q)$ as $\epsilon \to 0$. We have,
\begin{eqnarray*} 
\left|
\int_R \vartheta_\epsilon(x) (\psi_\epsilon(x)-\psi(x))dx \right|
\leq   \int_R |\vartheta_\epsilon(x)-\vartheta(x)|| \psi_\epsilon(x)-\psi(x)|dx  + \left| \int_R \vartheta(x) (\psi_\epsilon(x)-\psi(x))dx \right|.   
\end{eqnarray*}
By using Holder inequality, 
 \begin{eqnarray*} \int_R |\vartheta_\epsilon(x)-\vartheta(x)|| \psi_\epsilon(x)-\psi(x)|dx
 \leq   \|\vartheta_\epsilon -\vartheta\|_{L^{p^\prime}(Q)} \|\psi_\epsilon -\psi\|_{L^p(Q)} \leq 2K \|\vartheta_\epsilon-\vartheta\|_{L^{p^\prime}(Q)} \to 0.
\end{eqnarray*}
Since $\psi_\epsilon \wto \psi$ in $L^{p}(Q)$ then 
$$
\left| \int_R \vartheta(x) (\psi_\epsilon(x)-\psi(x))dx \right| \to 0,\quad \mbox{as $\epsilon \to 0$},
$$
and we obtain a contradiction to \eqref{notuniform}. \cqd

\begin{lemma} 
\label{convergence_nonlinearities}
Suppose that $f$ and $g$ satisfy \eqref{bounded_f} and \eqref{bounded_g}
and let $\frac{1}{2}<\alpha \leq 1,$ then %there exists a positive function $c(\epsilon) \to 0$ as $\epsilon \to 0$ such that 
$$
\|h_\epsilon (E_\epsilon u) - E^*_\epsilon h_0(u)\|_{X^{-\frac{\alpha}{2}}_{\epsilon}} \to 0,
$$
uniformly for $u$ in bounded sets of $X^{\frac{1}{2}}$. 
\end{lemma}
\proof We first consider $\alpha=1$. Since 
$$
\|h_\epsilon (E_\epsilon u) - E^*_\epsilon h_0(u)\|_{X^{-\frac{1}{2}}_\epsilon}= \sup _{\begin{array}{ccc}v_{\epsilon}\in X^{\frac{1}{2}}_\epsilon \\ \|v_\epsilon\|_{X^{\frac{1}{2}}_\epsilon}=1\end{array}}|\langle h_\epsilon (E_\epsilon u), v_\epsilon\rangle - \langle E^*_\epsilon h_0(u), v_\epsilon\rangle |.
$$

Let  $v_\epsilon \in X^{\frac{1}{2}}_\epsilon$ we have
\begin{equation}
\begin{array}{lll}
& &{\ds{|\langle h_\epsilon (E_\epsilon u), v_\epsilon\rangle - \langle E^*_\epsilon h_0(u), v_\epsilon\rangle | = |\langle h_\epsilon (E_\epsilon u), v_\epsilon\rangle - \langle  h_0(u), \hat{E}_\epsilon v_\epsilon\rangle | }}\\ 
% & = & {\ds{\left|\int_{\Omega_\epsilon} f(E_\epsilon w) v_\epsilon  + \int_{\partial\Omega_\epsilon} g(E_\epsilon w) v_\epsilon - \int_\Omega f(w) \hat{E}_\epsilon v_\epsilon - \int_{\partial\Omega} \gamma (x) g(w) \hat{E}_\epsilon v_\epsilon \right|}} \\  
& \leq & {\ds{\left|\int_{\Omega_\epsilon} f(E_\epsilon u) v_\epsilon -\int_\Omega f(u) \hat{E}_\epsilon v_\epsilon \right| 
 + \left| \int_{\partial\Omega_\epsilon} g(E_\epsilon u) v_\epsilon -  \int_{\partial\Omega} \gamma (x) g(u) \hat{E}_\epsilon v_\epsilon \right|}}. 
 \end{array}
\end{equation} 

Using \eqref{bounded_f}, Lemma \ref{auxiliar-similarLemma4.2AB},  $|\Omega_\epsilon\setminus K_\epsilon| \to 0$ and $|\Omega \setminus K_\epsilon| \to 0$ as $\epsilon \to 0$, we obtain 
$$ 
\begin{array}{lll}
 & & \displaystyle \left|\int_{\Omega_\epsilon} f(E_\epsilon u)v_\epsilon
  - \int_\Omega f(u)\hat{E}_\epsilon v_\epsilon\right|\\
  & \leq & \displaystyle\int_{\Omega_\epsilon\setminus K_\epsilon } |f(E_\epsilon u)v_\epsilon | + \int_{\Omega \setminus K_\epsilon } |f(u)\hat{E}_\epsilon v_\epsilon | + \int_{K_\epsilon}| f(u)(v_\epsilon-\hat{E}_\epsilon v_\epsilon)|\leq c_1(\epsilon) \|v_\epsilon\|_{X_\epsilon^{\frac{1}{2}}} ,
\end{array}
$$
with $c_1(\epsilon) \to 0$ as $\epsilon \to 0$.

Consider the finite cover $\{U_i\}_{i=0}^{m}$ such that $\overline{\Omega}_\epsilon \subset \cup_{i=0}^m U_i\equiv U$, the local parametrizations $\phi_{i,\epsilon}$ and $\phi_i$  of $\partial \Omega_\epsilon$ and $\partial \Omega$, respectively. By Definition \ref{definition-of-gamma}, we have $\gamma (\phi_{i}(x')) J_{N-1}\phi_i(x') = \gamma_i(x')$ for $x'\in Q_{N-1}$. Hence, 
\begin{eqnarray*} 
\begin{array}{llll}
& & \displaystyle \left| \int_{U_i\cap \partial\Omega_\epsilon} g(E_\epsilon u)  v_\epsilon -  \int_{U_i \cap \partial\Omega} \gamma (x) g(u) \hat{E}_\epsilon v_\epsilon \right| \\  
% & = & \displaystyle \left|\int_{Q_{N-1}} g((E_\epsilon w)(\phi_{i,\epsilon}(x'))) v_\epsilon(\phi_{i,\epsilon}(x')) J_{N-1}\phi_{i,\epsilon}(x')dx'  \right. \\
% &- & \displaystyle \left. \int_{Q_{N-1}} \gamma (\phi_{i}(x')) g(w(\phi_{i}(x'))) (\hat{E}_\epsilon v_\epsilon)(\phi_{i}(x'))J_{N-1}\phi_i(x') dx' \right| \\
%& = & \displaystyle \left|\int_{Q_{N-1}} g((E_\epsilon w)(\phi_{i,\epsilon}(x'))) v_\epsilon(\phi_{i,\epsilon}(x')) J_{N-1}\phi_{i,\epsilon}(x')dx' \right. \\
%& - &  \displaystyle \left. \int_{Q_{N-1}} g(w(\phi_{i}(x'))) (\hat{E}_\epsilon v_\epsilon)(\phi_{i}(x')) \tilde{\gamma}_i(x') dx' \right|.
& = & \displaystyle \left|\int_{Q_{N-1}} \!\!\!g((E_\epsilon u)(\phi_{i,\epsilon}(x'))) v_\epsilon(\phi_{i,\epsilon}(x')) J_{N-1}\phi_{i,\epsilon}(x')dx' \right. \\
  & - & \left. \displaystyle  \int_{Q_{N-1}}\!\!\!g(u(\phi_{i}(x'))) (\hat{E}_\epsilon v_\epsilon)(\phi_{i}(x')) \gamma_i(x') dx' \right|.
\end{array}
\end{eqnarray*}

Since  $(E_\epsilon u)_{|\Omega}=u$ and $\phi_i(x') \in U_i \cap \partial \Omega$, then $g((E_\epsilon u)(\phi_{i}(x'))) = g(u(\phi_i(x')))$ and 
\begin{eqnarray*}
\begin{array}{lll} 
& & \displaystyle \left|\int_{Q_{N-1}}g((E_\epsilon u)(\phi_{i,\epsilon}(x'))) v_\epsilon(\phi_{i,\epsilon}(x')) J_{N-1}\phi_{i,\epsilon}(x')dx' \right. 
\\
&- & \displaystyle \left. \int_{Q_{N-1}} g(u(\phi_{i}(x'))) (\hat{E}_\epsilon v_\epsilon)(\phi_{i}(x')) \gamma_i(x') dx'\right|
\\
& \leq & \displaystyle \int_{Q_{N-1}} |g((E_\epsilon u)(\phi_{i,\epsilon}(x')))-  g((E_\epsilon u)(\phi_{i}(x')))| |v_\epsilon(\phi_{i,\epsilon}(x'))| J_{N-1}\phi_{i,\epsilon}(x') dx'
\\ 
& + & \displaystyle \int_{Q_{N-1}}  |g(u(\phi_{i}(x')))||v_\epsilon(\phi_{i,\epsilon}(x'))-  (\hat{E}_\epsilon v_\epsilon)(\phi_{i}(x'))|J_{N-1}\phi_{i,\epsilon}(x') dx'
\\ 
& + & \left|\displaystyle \int_{Q_{N-1}} g(u(\phi_{i}(x'))) (\hat{E}_\epsilon v_\epsilon)(\phi_{i}(x'))( J_{N-1}\phi_{i,\epsilon}(x') -\gamma_i(x') )dx' \right|=\displaystyle (I) + (II)+(III).
\end{array}
\end{eqnarray*}

Using Mean Value Theorem, \eqref{bounded_g} and  \cite[Lemma 4.2]{AB0}, there exists $c_2(\epsilon)\to 0$ as $\epsilon \to 0$ such that
$$
\begin{array}{lll}
\displaystyle \|g((E_\epsilon u)(\phi_{i,\epsilon}))-   g((E_\epsilon u)(\phi_{i}))\|_{L^2(Q_{N-1})} \leq c_2(\epsilon)\| u\|_{H^1(U)},
\end{array}
$$
By similar arguments as   \cite[Lemma 4.2]{AB0}, there exists $c_3(\epsilon) \to 0$  as $\epsilon \to 0$ such that 
$$
\begin{array}{l}\displaystyle \| v_\epsilon(\phi_{i,\epsilon})-  (\hat{E}_\epsilon v_\epsilon)(\phi_{i})\|_{L^2(Q_{N-1})} \leq c_3(\epsilon)\|v_\epsilon\|_{H^1(U)}.
\end{array}
$$
Thus, using Cauchy-Schwartz inequality and \eqref{bounded_g}, we have
\begin{eqnarray*} 
(I) \leq c_2(\epsilon)\|J_{N-1}\phi_{i,\epsilon}\|_{L^\infty(Q_{N-1})} \| u\|_{H^1(U)}\|v_\epsilon(\phi_{i,\epsilon})\|_{L^2(Q_{N-1})}, 
\end{eqnarray*}
\begin{eqnarray*} 
(II) \leq \tilde{c}_3(\epsilon)\|J_{N-1}\phi_{i,\epsilon}\|_{L^\infty(Q_{N-1})} \|v_\epsilon\|_{H^1(U)},  
\end{eqnarray*}
with $\tilde{c}_3(\epsilon) \to 0$  as $\epsilon \to 0$.

Also, by hypothesis {\bf{(F)(ii)}}, $J_{N-1}\phi_{i,\epsilon} \stackrel{\epsilon \to 0}{-{\hspace{-2mm}}\rightharpoonup}  \gamma_i$ in $L^p(Q_{N-1})$, $1\leq p<\infty$, and  by Lemma \ref{weakconvergenceissufficient}, 
\begin{eqnarray*} 
(III) \leq  c_4(\epsilon, v_\epsilon ),\end{eqnarray*}
with $c_4(\epsilon,v_\epsilon) \to 0$ as $\epsilon \to 0$, uniformly for $v_\epsilon$ in bounded sets of $X^{\half}_\epsilon$.

Therefore, there exists $c(\epsilon)\to 0$ as $\epsilon \to0$ such that
\begin{equation}
\label{Conv1_heps_h0}
\|h_\epsilon (E_\epsilon u) - E^*_\epsilon h_0(u)\|_{X^{-\frac{1}{2}}_\epsilon}\leq c(\epsilon).
\end{equation}

Now, fix $\frac{1}{2}<\alpha_0<1$, by Lemma \ref{limitacaoinfinito} we have that there exists $K>0$ independent of $\epsilon$ such that
\begin{equation}
\label{Conv2_heps_h0}    
\|h_\epsilon (E_\epsilon u) - E^*_\epsilon h_0(u)\|_{X^{-\frac{\alpha_0}{2}}_\epsilon} \leq K.
\end{equation}
Then for any $\alpha$ such that $-1<-\alpha<-\alpha_0<-\frac{1}{2}$, using Lemma \ref{equivalentdefinitionsspacesnorms} and interpolation properties (see \cite[Theorem 1.11.3]{triebel} and  \cite[pag. 15]{Yagi}), \eqref{Conv1_heps_h0} and \eqref{Conv2_heps_h0}, we obtain
$$
\begin{array}{lll}
\displaystyle & & \|h_\epsilon (E_\epsilon u) - E^*_\epsilon h_0(u)\|_{X^{-\frac{\alpha}{2}}_{\epsilon}} \leq C\|h_\epsilon (E_\epsilon u) - E^*_\epsilon h_0(u)\|_{H^{-\alpha}(\Omega_\epsilon)}\\
\displaystyle & \leq & C\|h_\epsilon (E_\epsilon u) - E^*_\epsilon h_0(u)\|^{\theta}_{H^{-\alpha_0}(\Omega_\epsilon)}\|h_\epsilon (E_\epsilon u) - E^*_\epsilon h_0(u)\|^{1-\theta}_{H^{-1}(\Omega_\epsilon)} \leq c(\epsilon),
\end{array}
$$
for some $0<\theta<1$ and $c(\epsilon) \to 0$ as $\epsilon \to 0$.
\cqd

Now, we are in condition to
 verify completely
the hypothesis  {\bf{[A1]}}, that is,
$h_\epsilon(u_\epsilon) {\stackrel{E^*}{\longrightarrow}} \ h_0(u_0)$ whenever $u_\epsilon \dto u_0$.

\begin{lemma} \label{A1-3rdcondition} 
Suppose that $f$ and $g$ satisfy \eqref{bounded_f} and \eqref{bounded_g}
and let $\frac{1}{2}<\alpha \leq 1$. If $u_\epsilon \dto u_0$ in $X^{\frac{1}{2}}_{\eps}$,  then 
$$
\|h_\epsilon (u_\epsilon) - E^*_\epsilon h_0(u_0)\|_{X^{-\frac{\alpha}{2}}_{\epsilon}} \to 0,\quad \mbox{as $\epsilon \to 0$}.
$$
\end{lemma}
\proof By Lemmas  \ref{lipschitz_h} and  \ref{convergence_nonlinearities}, we obtain 
$$
\begin{array}{lll}
\|h_\epsilon (u_\epsilon) - E^*_\epsilon h_0(u_0)\|_{X^{-\frac{\alpha}{2}}_\epsilon} &\leq &\|h_\epsilon (u_\epsilon) - h_\epsilon (E_\epsilon u_0)\|_{X^{-\frac{\alpha}{2}}_\epsilon} + 
\|h_\epsilon (E_\epsilon u_0) - E^*_\epsilon h_0(u_0)\|_{X^{-\frac{\alpha}{2}}_\epsilon}\\
&\leq & L\|u_\epsilon - E_\epsilon u_0 \|_{X^{\frac{1}{2}}_\epsilon} + 
\|h_\epsilon (E_\epsilon u_0) - E^*_\epsilon h_0(u_0)\|_{X^{-\frac{\alpha}{2}}_\epsilon} \to 0, \quad \mbox{as $\epsilon \to 0$}.
\end{array}
$$
\cqd

Now, we will verify the properties in {\bf{[A2]}} related to the derivatives of 
 nonlinearities \eqref{hepsilon} and \eqref{hzero}. For this, we define 
$h'_\epsilon: X_\epsilon^\half\to \mathcal{L}(X_\epsilon^\half,X_\epsilon^{-\frac{\alpha}{2}})
$, with $0\leq \eps \leq \eps_0$ and $\half <\alpha \leq 1$, by
\begin{equation} 
\label{Dhepsilon}
\langle h'_\epsilon(e_{\epsilon})v_\epsilon, z_\epsilon\rangle = \int_{\Omega_\epsilon} f'(e_\epsilon)v_\epsilon z_\epsilon + \int_{\partial\Omega_\epsilon}g'(e_\epsilon) v_\epsilon z_\epsilon,  \quad \mbox{for $e_\eps,v_\epsilon \in X^{\frac{1}{2}}_\epsilon$ and $z_\epsilon \in X^{\frac{\alpha}{2}}_\epsilon$},
\end{equation}
and  
\begin{equation} 
\label{Dhzero}
\langle h'_0(e_{0})v, z \rangle = \int_{\Omega} f'(e_{0}) vz +\int_{\partial\Omega}\gamma(x) g'(e_{0}) vz,  \quad \mbox{for $e_0,v \in X^{\frac{1}{2}}$ and $z \in X^{\frac{\alpha}{2}}$}.
\end{equation}

First, by definition of $h^\prime_\epsilon$ in \eqref{Dhepsilon} and \eqref{Dhzero} and the hypotheses \eqref{bounded_f} and \eqref{bounded_g}, we can proceed
 in a similar way to Lemma \ref{limitacaoinfinito} and we can obtain 
\begin{equation} 
\label{A2-2ndcondition}
 \sup_{\epsilon \in [0,\epsilon_0]} \sup_{\|v_\epsilon\|_{X_\epsilon^\half} \leq \rho}\|h^\prime_\epsilon(v_\epsilon + e_\epsilon)\|_{\mathcal{L}(X_\epsilon^\half,X_\epsilon^{-\frac{\alpha}{2}})} \leq K,
\end{equation}
for some $\rho>0$ and $K>0$ independent of $\epsilon$. 

The following lemmas prove the $E^*$-convergence of  $h^\prime_\epsilon(e_\epsilon)$ to $h^\prime_0(e_0)$ whenever $e_\epsilon \dto e_0$. 
\begin{lemma} 
\label{E*convergencenonlinearities-linearization} 
Suppose that $f$ and $g$ satisfy \eqref{bounded_f} and \eqref{bounded_g} and let $\half <\alpha \leq 1$. If $e_\epsilon \dto e_0$ in $X^{\frac{1}{2}}_{\eps}$, then there exists a positive function $C(\epsilon)\to 0$ as $\epsilon \to 0$ such that  $$
\| h^\prime_{\epsilon}(e_\epsilon)(E_\epsilon v) - E^*_\epsilon h^\prime_0(e_0)(v)\|_{X^{-\frac{\alpha}{2}}_{\epsilon}} \leq C(\epsilon) \|v\|_{X^\half}, \quad \mbox{for all $v\in X^\half$}.
$$
\end{lemma}
\proof 
Initially, for each $\epsilon \in [0,\epsilon_0]$ and $\alpha=1$, we have
$$
\|h^\prime_{\epsilon}(e_\epsilon)(E_\epsilon v) - E^*_\epsilon h^\prime_0(e_0)(v)\|_{X^{-\frac{1}{2}}_{\epsilon}} =\sup_{\begin{array}{ccc}z_\epsilon \in X^{\frac{1}{2}}_{\epsilon}\\ \|z_\epsilon\|_{X^{\frac{1}{2}}_{\epsilon}}=1\end{array}}  | \langle h^\prime_{\epsilon}(e_\epsilon)(E_\epsilon v),z_\epsilon\rangle - \langle E^*_\epsilon h^\prime(e_0)(v),z_\epsilon\rangle|. 
$$
We note
\begin{eqnarray*}
& &|\langle h^\prime_{\epsilon}(e_\epsilon)(E_\epsilon v),z_\epsilon\rangle - \langle E^*_\epsilon h^\prime(e_0)(v),z_\epsilon\rangle|=|\langle h^\prime_{\epsilon}(e_\epsilon)(E_\epsilon v),z_\epsilon\rangle - \langle  h^\prime(e_0)(v),\hat{E}_\epsilon z_\epsilon\rangle|\\
%& = &  \left| \int_{\Omega_\epsilon}    f^\prime(e_\epsilon)(E_\epsilon v) z_\epsilon -  \int_{\partial\Omega_\epsilon}    g^\prime(e_\epsilon)(E_\epsilon v) z_\epsilon 
%-  \int_{\Omega}    f^\prime(e_0)(v) \hat{E}_\epsilon z_\epsilon +   \left. \int_{\partial\Omega} \gamma (\zeta) g^\prime(e_0)(v) \hat{E}_\epsilon z_\epsilon  \right| \\
& =  & \left|\int_{\Omega_\epsilon}   f^\prime(e_\epsilon)( E_\epsilon v)z_\epsilon  - \int_\Omega  f^\prime(e_0)( v)\hat{E}_\epsilon z_\epsilon  
 + \int_{\partial\Omega_\epsilon} g^\prime(e_\epsilon)( E_\epsilon v)z_\epsilon -  \int_{\partial\Omega}   \gamma(x)g^\prime(e_0)( v)\hat{E}_\epsilon z_\epsilon \right| \\
 & \leq &  \left|\int_{\Omega_\epsilon}   f^\prime(e_\epsilon)( E_\epsilon v)z_\epsilon  -\int_{\Omega_\epsilon}   f^\prime(E_\epsilon e_0)( E_\epsilon v)z_\epsilon \right|  
+ \left| \int_{\Omega_\epsilon}   f^\prime(E_\epsilon e_0)( E_\epsilon v)z_\epsilon - \int_\Omega  f^\prime(e_0)( v)\hat{E}_\epsilon z_\epsilon \right| \\ 
& + &\left| \int_{\partial\Omega_\epsilon} g^\prime(e_\epsilon)( E_\epsilon v)z_\epsilon -\int_{\partial\Omega_\epsilon}   g^\prime(E_\epsilon e_0)( E_\epsilon v)z_\epsilon \right| 
+ \left|\int_{\partial\Omega_\epsilon}   g^\prime(E_\epsilon e_0)( E_\epsilon v)z_\epsilon - \int_{\partial\Omega}   \gamma(x)g^\prime(e_0)( v)\hat{E}_\epsilon z_\epsilon \right|.
\end{eqnarray*}

Using Mean Value Theorem, we can write
\begin{eqnarray*}
f^\prime(e_\epsilon(x))-f^\prime((E_\epsilon e_0)(x)) =  f^{\prime\prime}(\tilde{e}_\epsilon(x)e_\epsilon(x) +(1-\tilde{e}_\epsilon(x)) (E_\epsilon e_0)(x))[e_\epsilon(x) - (E_\epsilon e_0)(x)], \quad x \in \Omega_\epsilon,
\end{eqnarray*}
and
\begin{eqnarray*}
g^\prime(e_\epsilon(x))-g^\prime((E_\epsilon e_0)(x)) = g^{\prime\prime}(\hat{e}_\epsilon(x)e_\epsilon(x) +(1-\hat{e}_\epsilon(x)) (E_\epsilon e_0)(x))[e_\epsilon(x) - (E_\epsilon e_0)(x)], \quad x \in \partial \Omega_\epsilon,
\end{eqnarray*}
where $\tilde{e}_\epsilon(x)\in [0,1]$ for all $x \in \Omega_\epsilon$ and $\hat{e}_\epsilon(x)\in [0,1]$ for all $x \in \partial \Omega_\epsilon$. Hence and using \eqref{bounded_f} and \eqref{bounded_g}, we get 
\begin{eqnarray*}
 & & |\langle  h^\prime_{\epsilon}(e_\epsilon)(E_\epsilon v),z_\epsilon\rangle - \langle E^*_\epsilon h^\prime(e_0)(v),z_\epsilon\rangle| \\
 & \leq &  C\int_{\Omega_\epsilon}|e_\epsilon - E_\epsilon e_0||E_\epsilon v||z_\epsilon|+\int_{\Omega_\epsilon\setminus K_\epsilon} |f^\prime(E_\epsilon e_0)|E_\epsilon v||z_\epsilon|   + \int_{\Omega \setminus K_\epsilon} |f^\prime( e_0)|v||\hat{E}_\epsilon z_\epsilon|  
\\&+&\int_{K_\epsilon\setminus K_{\epsilon_0} } |f^\prime( e_0)| v||z_\epsilon  - z_\epsilon\circ\theta_\epsilon|  
+  C\int_{\partial \Omega_\epsilon}|e_\epsilon - E_\epsilon e_0||E_\epsilon v||z_\epsilon| \\ &+& \left| \int_{\partial \Omega_\epsilon}   g^\prime(E_\epsilon e_0)( E_\epsilon v)z_\epsilon -  \int_{\partial\Omega}   \gamma(x)g^\prime(e_0)( v)\hat{E}_\epsilon z_\epsilon \right|= (I) +(II)+(III)+(IV)+(V)+(VI).
\end{eqnarray*}

For the last term (VI), we use the parametrizations $\phi_i$ and $\phi_{i,\epsilon}$ of the boundaries $\partial \Omega$ and $\partial \Omega_\epsilon$ respectively, then with similar arguments to Lemma \ref{convergence_nonlinearities} and considering that
$J_{N-1}\phi_{i,\epsilon} \stackrel{\epsilon \to 0}{-{\hspace{-2mm}}\rightharpoonup}  \gamma_i$ in $L^p(Q_{N-1})$, $1\leq p<\infty$. Applying Lemma \ref{auxiliar-similarLemma4.2AB} we get that term (IV)  goes to zero as $\epsilon\to 0$.

Using $e_\epsilon \dto e_0$, $|\Omega_\epsilon\setminus K_\epsilon| \to 0$ and $|\Omega \setminus K_\epsilon| \to 0$, we have that terms (I), (II), (III) and (V) goes to zero as $\epsilon\to 0$.  Thus, there exists a positive function $\tilde{C}(\epsilon)\to 0$ as $\epsilon \to 0$ such that
\begin{equation}
\label{1Conv_Derivadas_X_-1_2}
\|h^\prime_{\epsilon}(e_\epsilon)(E_\epsilon v) - E^*_\epsilon h^\prime_0(e_0)(v)\|_{X^{-\frac{1}{2}}_{\epsilon}}\leq \tilde{C}(\epsilon)\|v\|_{X^\half}. 
\end{equation}

Now, fix $\frac{1}{2}<\alpha_0<1$, by 
 similar arguments to \eqref{A2-2ndcondition},
we have that there exists $K>0$ independent of $\epsilon$ such that
\begin{equation}
\label{2Conv_Derivadas_X_-1_2}
\|h^\prime_{\epsilon}(e_\epsilon)(E_\epsilon v) - E^*_\epsilon h^\prime_0(e_0)(v)\|_{X^{-\frac{\alpha_0}{2}}_{\eps}}\leq K \|v\|_{X^\half}, \quad \mbox{for all $v\in X^\half$}.
\end{equation}
Finally, for any $\alpha$ such that $-1<-\alpha<-\alpha_0<-\frac{1}{2}$, using Lemma \ref{equivalentdefinitionsspacesnorms}, interpolation properties (see \cite[Theorem 1.11.3]{triebel} and  \cite[pag. 15]{Yagi}), \eqref{1Conv_Derivadas_X_-1_2} and \eqref{2Conv_Derivadas_X_-1_2}, we obtain
$$
\begin{array}{lll}
\displaystyle & & \|h^\prime_{\epsilon}(e_\epsilon)(E_\epsilon v) - E^*_\epsilon h^\prime_0(e_0)(v)\|_{X^{-\frac{\alpha}{2}}_{\epsilon}} \leq C\|h^\prime_{\epsilon}(e_\epsilon)(E_\epsilon v) - E^*_\epsilon h^\prime_0(e_0)(v)\|_{H^{-\alpha}(\Omega_\epsilon)}\\
\displaystyle & \leq & C\|h^\prime_{\epsilon}(e_\epsilon)(E_\epsilon v) - E^*_\epsilon h^\prime_0(e_0)(v)\|^{\theta}_{H^{-\alpha_0}(\Omega_\epsilon)}\|h^\prime_{\epsilon}(e_\epsilon)(E_\epsilon v) - E^*_\epsilon h^\prime_0(e_0)(v)\|^{1-\theta}_{H^{-1}(\Omega_\epsilon)}\leq C(\epsilon) \|v\|_{X^\half},
\end{array}
$$
for all $v\in X^\half$ and
for some $0<\theta<1$ and $C(\epsilon) \to 0$ as $\epsilon \to 0$.
\cqd

Now, we are in condition to
 verify completely
the hypothesis  {\bf{[A2]}}, that is, $ h^\prime_{\epsilon}(e_\epsilon)  \ddto h^\prime_0(e_0)$, 
whenever $e_\epsilon\dto e_0$.

\begin{lemma} % firstcondition_in_[A2]
\label{iii-compactconvergence}
Suppose that $f$ and $g$ satisfy \eqref{bounded_f} and \eqref{bounded_g} and let $\half <\alpha \leq 1$. If  $e_\epsilon \dto e_0$ and $v_\epsilon \dto v$ in $X^{\frac{1}{2}}_{\eps}$, then 
$$\| h^\prime_{\epsilon}(e_\epsilon)(v_\epsilon) - E^*_\epsilon h^\prime_0(e_0)(v)\|_{X^{-\frac{\alpha}{2}}_{\epsilon}}  \to 0,\quad \mbox{as }\epsilon \to 0. $$
\end{lemma}
\proof Using \eqref{A2-2ndcondition}, similar arguments to Lemma \ref{lipschitz_h},  and Lemma \ref{E*convergencenonlinearities-linearization}, we obtain 
\begin{eqnarray*} 
& &\| h^\prime_{\epsilon}(e_\epsilon)(v_\epsilon) - E^*_\epsilon h^\prime_0(e_0)(v)\|_{X^{-\frac{\alpha}{2}}_{\epsilon}} \\ &\leq & 
\| h^\prime_{\epsilon}(e_\epsilon)(v_\epsilon) -  h_\epsilon ^\prime(e_\epsilon)(E_\epsilon v)\|_{X^{-\frac{\alpha}{2}}_{\epsilon}} + \| h_\epsilon^\prime(e_\epsilon)(E_\epsilon v) - E^*_\epsilon h^\prime_0(e_0)(v)\|_{X^{-\frac{\alpha}{2}}_{\epsilon}} \\ &\leq & 
L\| v_\epsilon -  E_\epsilon v\|_{X^{\frac{1}{2}}_{\epsilon}} + \| h_\epsilon^\prime(e_\epsilon)(E_\epsilon v) - E^*_\epsilon h^\prime_0(e_0)(v)\|_{X^{-\frac{\alpha}{2}}_{\epsilon}} \to 0,\quad \mbox{as $\epsilon \to 0$},
\end{eqnarray*}
for some $L>0$ independent of $\epsilon$.
\cqd

\begin{remark}
The results of Subsections \ref{A1_Resolvent} and \ref{A2_nonlinearity} imply that the hypotheses {\bf[A1]} and {\bf[A2]} hold. Consequently, there exists a family of attractors $\left\{\mathcal{A}_{\epsilon}\right\}_{\epsilon\in [0,\epsilon_0]}$ for our problems (\ref{nbc}) and (\ref{nbc_limite_gamma_F}) or \eqref{nbcabstract} in $Y^{\frac{1+\alpha}{2}}_{\eps}$ and we obtain the
$\mathcal{P}^{\frac{1+\alpha}{2}}$-continuity of this family in $Y^{\frac{1+\alpha}{2}}_{\eps}$ at $\eps=0$.   
\end{remark}

%%%%%%%%%%%%%%%%%%%%%
\section{$E$-Continuity of the attractors}
\label{continuity_attractors_H1}

We have shown in  Section \ref{P_continuity_attractors} the $\mathcal{P}^{\frac{1+\alpha}{2}}$-continuity of the attractors. For this we have
applied the abstract results from \cite{CP}. Now, we want to obtain the $E$-continuity of the attractors. In order to accomplish this, we need to show the equivalence between the concepts of
$\mathcal{P}^{\frac{1+\alpha}{2}}$-convergence and $E$-convergence.

%%%%%%%%%%%%%%%%%%%%%%%%%%%%%%%%%%%%%%%%%%%%%%%%%%%%%%%%
%%%%%%%%%%

\subsection{$E$-convergence vs $\mathcal{P}$-convergence}
\label{EvsP}

The purpose of this subsection is to prove that the notion of $\mathcal{P}$-convergence established in an abstract way in the paper of \cite{CP} is the same as the notion of $E$-convergence we are using in this paper.

Based on the notations in Subsection \ref{Carvalho-Piskarev}, 
if we fix a value $\alpha\in (0,1)$ and we choose  $Y_\eps=X^{-\frac{\alpha}{2}}_\eps$, $\eps\in [0,\eps_0]$, where $Y_0=Y$ and $X^{-\frac{\alpha}{2}}_{0}=X^{-\frac{\alpha}{2}}$, and consider them as the base spaces of the operators $A_\eps$. Analogously to Section \ref{P_continuity_attractors}, we define the operator $p_\eps:Y\to Y_\eps$, that is, $p_\eps:X^{-\frac{\alpha}{2}}\to X^{-\frac{\alpha}{2}}_{\epsilon}$ as $p_\eps =E^*_\epsilon$, where $E^*_\epsilon:X^{-\frac{\alpha}{2}}\to X^{-\frac{\alpha}{2}}_{\epsilon}$ is given by (\ref{extensaodual}), then we have $Y_\eps^{\frac{1+\alpha}{2}}=X^{\frac{1}{2}}_\eps=H^1(\Omega_\eps)$, $\epsilon \in[0,\epsilon_0]$,  and $p_\eps^{\frac{1+\alpha}{2}}:Y^{\frac{1+\alpha}{2}}\to Y^{\frac{1+\alpha}{2}}_\eps$ given by
$p_\eps^{\frac{1+\alpha}{2}} x=(A_\eps)^{-\frac{1+\alpha}{2}} E^*_\eps (A_0)^{\frac{1+\alpha}{2}} x$, for $x\in Y^{\frac{1+\alpha}{2}}$.

So we have two concepts of convergence, one given by Definition \ref{Definition_Concept_E_Convergence} using the family of linear operators $E_\epsilon:X^{\frac{1}{2}}\to X^{\frac{1}{2}}_{\epsilon}$  in \eqref{definitionE_epsilon} and another  given by Definition \ref{Definition_Concept_P_Convergence} using the family of  operators $p_\eps^{\frac{1+\alpha}{2}}:X^{\frac{1}{2}}\to X^{\frac{1}{2}}_{\epsilon}$. 
%also  have the $E$-convergence, where  $E_\eps: Y_0^{(1+\alpha)/2} \to Y_\epsilon^{(1+\alpha)/2}$ is defined in \ref{definitionE_epsilon}. 
We are going to prove that these two concepts of convergence are equivalents.

\begin{lemma}  \label{equivalence_E_Palpha}
Let $\{\Omega_\epsilon\}_{\epsilon \in [0,\epsilon_0]}$ be a family of domains satisfying conditions {\bf{(H)}}, {\bf{(F)}(i)} and {\bf{(I)}}. If $u_\epsilon \in X^{\half}_\epsilon$ and $u \in X^\half$,  then $u_\epsilon\dto u$ if, and only if, $u_\epsilon {\stackrel{\mathcal{P}^{\frac{1+\alpha}{2}}}{\longrightarrow}} u$. 
\end{lemma} 
\proof
We first prove that given $u \in X^\half = H^1(\Omega)$, the sequence $\{p_\eps^{\frac{1+\alpha}{2}}\ u\}$  $E$-converges to  $u$ in $X^{\half}_\epsilon= H^{1}(\Omega_{\epsilon})=Y_\epsilon^{\frac{1+\alpha}{2}} $,  as $\epsilon \to 0$, that means, 
\begin{equation}
\label{Equation_1_Dia}
\|p_\eps^{\frac{1+\alpha}{2}}\ u - E_\epsilon u\|_{X^{\half}_\epsilon} \to 0 , \mbox{ \ as \ } \epsilon \to 0. 
\end{equation}
In fact, let $u \in X^{\half}$ and $v\in X^{-\frac{\alpha}{2}}$ such that $u=A_0^{-\frac{1+\alpha}{2}}v$, then taking $v_\epsilon = E^*_\epsilon v$ and $u_\epsilon=A_\epsilon^{-\frac{1+\alpha}{2}}v_\epsilon$, we have $u_\epsilon = p_\eps^{\frac{1+\alpha}{2}}\ u$. 
%We are going to estimate $$\|u_\epsilon - E_\epsilon u\|_{H^1(\Omega_\epsilon)}$$

Now, let $\{(\lambda_i,\varphi_i)\}_{i=1}^{\infty}$ the eigenvalues and eigenfunctions of $A_0: X^{\half} \to X^{-\half}$
and $\{(\lambda_{i,\epsilon},\varphi_{i, \epsilon}) \}_{i=1}^{\infty}$ the eigenvalues and eigenfunctions of $A_\epsilon:  X^{\half}_\epsilon \to X^{-\half}_\epsilon $. Let $u \in X^{\half}$ then 
$$
u=\sum_{i=1}^{\infty} (u,\varphi_i)_{X^{\half}} \varphi_i \quad \mbox{and} \quad \|u\|_{X^{\half}}^2= \sum_{i=1}^{\infty} (u,\varphi_i)^2_{X^{\half}} \lambda_i < \infty.
$$
So, $v= A^{\frac{1+\alpha}{2}}_{0} u = \displaystyle \sum_{i=1}^{\infty} (u,\varphi_i)_{X^{\half}} A^{\frac{1+\alpha}{2}}_{0}\varphi_i =  \sum_{i=1}^{\infty} (u,\varphi_i)_{X^{\half}} \lambda_i^{\frac{1+\alpha}{2}}\varphi_i.$ On the other hand, 
$$
E^*_\epsilon v = v_\epsilon =  \displaystyle \sum_{i=1}^{\infty} (v_\epsilon,\varphi_{i,\epsilon} )_{X^{-\frac{\alpha}{2}}_{\epsilon}} \varphi_{i,\epsilon} = \sum_{i=1}^{\infty} v_{i,\epsilon} \varphi_{i,\epsilon},
$$
where, by definition of $E^*_\epsilon$,  
$$
\begin{array}{lll}
v_{i,\epsilon} & = & \displaystyle (v_\epsilon,\varphi_{i,\epsilon} )_{X^{-\frac{\alpha}{2}}_{\epsilon}}  =  (E^*_\epsilon v,\varphi_{i,\epsilon} )_{X^{-\frac{\alpha}{2}}_{\epsilon}} = ( v,\hat{E}_\epsilon\varphi_{i,\epsilon} )_{X^{-\frac{\alpha}{2}}} 
\\
& = &  \displaystyle ( v,\hat{E}_\epsilon\varphi_{i,\epsilon} - \varphi_i )_{X^{-\frac{\alpha}{2}}} + ( v, \varphi_{i})_{X^{-\frac{\alpha}{2}}}.
\end{array}
$$

From Lemma \ref{compactlyconvergence} we have that $A_\epsilon^{-1} {\stackrel{CC}{\longrightarrow}} A_0^{-1}$. Hence,
using \cite[Proposition 3.3]{AB0}, we have $\lambda_{i,\epsilon} \to \lambda_i$ and $\varphi_{i,\epsilon} \dto \varphi_i$  in $X^{\frac{1}{2}}_{\epsilon}$. So, from Lemma \ref{Econvergence->E*convergence} we obtain that $\varphi_{i,\epsilon} {\stackrel{E^*}{\longrightarrow}} \varphi_i$ in $X^{-\frac{\alpha}{2}}_{\epsilon}$. Moreover, as consequence of the Lemma \ref{convergence_restriction_ofthe_extension},  $\hat{E}_\epsilon \varphi_{i,\epsilon}$  converges to $\varphi_i$ in $X^{\frac{\alpha}{2}}$, as $\epsilon \to 0$. Thus, 
$$
v_{i,\epsilon} = (v_\epsilon,\varphi_{i,\epsilon} )_{X^{-\frac{\alpha}{2}}_{\epsilon}}  \to ( v, \varphi_{i} )_{X^{-\frac{\alpha}{2}}} = v_i, \quad \mbox{as $\epsilon \to 0$}.
$$

Also, by Lemma  \ref{conditionE*convergence},   $\|E^*_\epsilon v\|_{X^{-\frac{\alpha}{2}}_\epsilon} \to \| v\|_{X^{-\frac{\alpha}{2}}}$. Thus, 
\begin{equation}
\label{Eq_N1_Aplic_Prop3.1_Prop3.2}
\|u_\epsilon\|_{X^\half_\epsilon} = \|v_\epsilon\|_{X^{-\frac{\alpha}{2}}_\epsilon} =  \|E^*_\epsilon v \|_{X^{-\frac{\alpha}{2}}_\epsilon} \to \| v\|_{X^{-\frac{\alpha}{2}}} = \|u\|_{X^\half}, \quad \mbox{as $\epsilon \to 0$}.
\end{equation}

Now, we prove that $u_\epsilon = A_\epsilon^{- \frac{1+\alpha}{2}} v_\epsilon = \displaystyle \sum_{i=1}^{\infty}v_{i,\epsilon} \lambda_{i,\epsilon}^{-\frac{1+\alpha}{2}}\varphi_{i,\epsilon} $ $E$-converges to $u= \displaystyle \sum_{i=1}^{\infty}(u,\varphi_i)_{X^\frac{1}{2}} \varphi_{i}$, that means,
\eqref{Equation_1_Dia} holds. For this, considering \eqref{Eq_N1_Aplic_Prop3.1_Prop3.2}  and \cite[Propositions 3.1 and 3.2]{AB0}, it is sufficient to prove $u_{\epsilon} \dwto u$, that is, 
$$
(u_\epsilon,E_\epsilon w)_{X^\half_\epsilon} \to (u,w)_{X^\half}, \quad \mbox{for all $w\in X^\half$}. 
$$
In fact,
$$(u_\epsilon, E_\epsilon w)_{X^\half_\epsilon}  =   \sum_{i=1}^{\infty} v_{i,\epsilon} \lambda_{i,\epsilon}^{-\frac{1+\alpha}{2}}(\varphi_{i,\epsilon} , E_\epsilon w)_{X^{\half}_{\epsilon}} \to   \sum_{i=1}^{\infty} v_{i} \lambda_{i}^{-\frac{1+\alpha}{2}}(\varphi_{i} ,w)_{X^\half} = (u,w)_{X^\half}, \quad \mbox{as $\epsilon \to 0$}.$$

Finally, let $\{w_\epsilon\}$ be a family in $X^{\frac{1}{2}}_{\epsilon}$ and $w \in X^{\frac{1}{2}}$  such that $w_\epsilon\dto w$, then using \eqref{Equation_1_Dia}, we obtain
$$\|w_\epsilon - p_\epsilon^{\frac{1+\alpha}{2}}w\|_{X^{\frac{1}{2}}_{\epsilon}}  \leq \|w_\epsilon - E_\epsilon w\|_{X^{\frac{1}{2}}_{\epsilon}} + \|E_\epsilon w - p_\epsilon^{\frac{1+\alpha}{2}}w\|_{X^{\frac{1}{2}}_{\epsilon}}  \to 0, \quad \mbox{as $\epsilon \to 0$}.$$
Therefore, $w_\epsilon {\stackrel{\mathcal{P}^{\frac{1+\alpha}{2}}}{\longrightarrow}} w$. Using similar arguments, if $w_\epsilon {\stackrel{\mathcal{P}^{\frac{1+\alpha}{2}}}{\longrightarrow}} w$
then $w_\epsilon\dto w$. \cqd

\begin{remark}  
Notice that once we have established the equivalence of $\mathcal{P}^{\frac{1+\alpha}{2}}$-convergence and $E$-convergence in Lemma \ref{equivalence_E_Palpha} and since we have proven in  Subsection \ref{P_continuity_attractors} 
  the $\mathcal{P}^{\frac{1+\alpha}{2}}$-continuity of the attractors  in  $Y_\eps^{\frac{1+\alpha}{2}}=X^{\frac{1}{2}}_\eps=H^1(\Omega_\eps)$, then we obtain the $E$-continuity of the attractors in $H^1(\Omega_\eps)$. Therefore, the main result given by Theorem \ref{existence-upper-lower-attractors} is proved.
\end{remark}

%%%%%%%%%%%%%%%%%%%%%%%%%%%%%%%%%%%%%%%%%%%%%%%%%%%%%%%%%%%%%%%%%%%%%%%%%%%%%%%%

\newpage

\end{document}